\documentclass[11pt,twoside,openany,reqno]{amsart}

\usepackage{mathrsfs}
\usepackage{latexsym}
\usepackage{amsmath}
\usepackage{amssymb}
\usepackage{amsthm}
\usepackage{amscd}
\usepackage{amsfonts}
\usepackage{epsfig}
\usepackage{multicol}
\usepackage{CJK}
\usepackage{url}
\usepackage{bbm}
\usepackage{dutchcal}
\usepackage{arydshln}
\usepackage{amsmath}\numberwithin{equation}{section}

\def\gs{\begin{aligned}}
\def\egs{\end{aligned}}

\usepackage[OT2,T1]{fontenc}
\DeclareSymbolFont{cyrletters}{OT2}{wncyr}{m}{n}
\DeclareMathSymbol{\Sha}{\mathalpha}{cyrletters}{"58}

\theoremstyle{plain}
\newtheorem{thm}{Theorem}[section]
\newtheorem{definition}[thm]{Definition}
\newtheorem{theorem}[thm]{Theorem}

\newtheorem{conjecture}[thm]{Conjecture}
\newtheorem{corollary}[thm]{Corollary}

\newcommand{\pf}{\noindent\begin {proof}}
\newcommand{\epf}{\end{proof}}

\begin{document}

\title{The rank of 2-Selmer group associate to $\theta$-congruent numbers }

\author{Tao Wei, Xuejun Guo}
\dedicatory{Department of Mathematics, Nanjing University, Nanjing 210093, China
\\ 
weitao@smail.nju.edu.cn \ \ \ \ \ guoxj@nju.edu.cn}
\thanks{The authors are supported by National Nature Science Foundation of China (Nos. 11971226, 12231009).}

\date{}

\maketitle

\noindent

\begin{abstract} We study the parity of rank of $2$-${\rm Selmer}$ groups associated to $\pi/3$ and $2\pi/3$-congruent numbers. Our second result gives some positive densities about $\pi/3$ and $2\pi/3$ non-congruent numbers which can support the even part of Goldfeld's conjecture. We give some necessary conditions such that $n$ is non $\pi/3$-congruent number for elliptic curves $E_n$ whose Shafarevich-Tate group is non-trivial. In the last section, we show that for $n=pq\equiv 5(resp. \ 11)\pmod{24}$, the density of non $\pi/3$($resp.$ $2\pi/3$)-congruent numbers is at least 75\%, where $p,q$ are primes.

\end{abstract}

\medskip


\section{\bf Introduction}

Recall a positive integer $n$ is congruent number if it is the area of a right triangle with rational side lengths. The congruent numbers conjecture says, for $n\equiv1,2,3(\text{mod }8)$, $n$ is congruent; for $n\equiv5,6,7(\text{mod }8)$, almost all $n$ are non-congruent. So far, many breakthroughs have been made towards this conjecture. Tian-Yuan-Zhang \cite{TYZ17} and A. Smith \cite{S16} show that congruent number have positive density for $n\equiv1,2,3\pmod{8}$; A. Smith \cite{S17} proved the non-congruent part of this conjecture, i.e., for $n\equiv1,2,3\pmod{8}$, the density of $n$ being non-congruent is $1$. Their result depend on considering a certain quadratic twist family of elliptic curves are called "congruent number elliptic curves"(CNEC): $$A_n:ny^2=x^3-x.$$

Fujiwara \cite{F96} introduced a general definition: a positive integer $n$ is called a $\theta$-congruent number if there exists a rational $\theta$-triangle with area $n\sqrt{r^2-s^2}$, where $\theta$ be a real number with $0<\theta<\pi$, cos$\theta=s/r$ is rational and $r,s\in\mathbb{Z},\ {\rm gcd}(r,s)=1$, a rational $\theta$-triangle is a triangle with rational sides and an angle $\theta$. Note that congruent number can be seen as $\pi/2$-congruent number.

Another famous problem is "tiling a triangle"\cite{B10}. We call triangle can be $n$-tiled or $n$ is a tiling numbers if a triangle can be tiled into $n$ congruent triangles. Hibino and Kan in \cite{HK01} \cite{Ka00}, He-Hu-Tian in \cite{HHT21} builds infinity many tiling number by using Heegner points.  Fujiwara \cite{F96} show if $n$ is tiling number if and only if $n$ is $\pi/3$ or $2\pi/3$-congruent numbers.  

The $\theta$-congruent number problems can also be studied by elliptic curves: the elliptic curves associated to $\theta$ is defined by

$$E_{\theta}:y^2=x(x+(r+s))(x-(r-s)),$$
where cos$\theta=s/r$, $r,s\in\mathbb{Z},\ {\rm gcd}(r,s)=1$.
The quadratic twist family of $E_{\theta}$ is $E_{n,\theta}:y^2=x(x+n(r+s))(x-n(r-s))$ for $n\in\mathbb{Z}$ and $n$ is square free. 

There are some important criteria for $\theta$-congruent numbers by Fujiwara. 

\begin{theorem}\cite{F96}
	Let $n$ be any square free natural number, $0<\theta<\pi$. Then
	
$1)$ $n$ is $\theta$-congruent if and only if $E_{n,\theta}$ has a rational point of order greater
than $2$.

$2)$ For $n\nmid 6$, $n$ is $\theta$-congruent if and only if the Modell-Weil group $E_{n,\theta}(\mathbb{Q})$ has a positive
rank.
\end{theorem}

Let $L(s,E)$ denote the Hasse-Weil L-function of elliptic curve $E$ over $\mathbb{Q}$, $E(\mathbb{Q})$ the Modrell-Weil group of $E$ over $\mathbb{Q}$.
The $\phi$-Selmer group of $E$ is the subgroup of $H^1(\mathbb{Q},E[\phi])$ defined by $${\rm Ker}\left(H^1(\mathbb{Q},E[\phi])\rightarrow \prod_{p}H^1(\mathbb{Q}_p,E)\right).$$
The Shafarevich-Tate group $\Sha(E)$ is defined by
$$\Sha(E)={\rm Ker}\left( H^1(\mathbb{Q},E)\rightarrow\prod_{p}H^1(\mathbb{Q}_p,E)\right) .$$

We have the fundamental short exact sequence which is called the first 2-descents for elliptic curves \cite{C98}:
$$0\rightarrow E(\mathbb{Q})/2E(\mathbb{Q}) \rightarrow {\rm Sel}_2(E)\rightarrow \Sha(E)[2]\rightarrow 0.$$

It is worth noting that, the Selmer groups are very helpful to study Moderell-Weil groups. For example, Alexander Smith\cite{S17} proved the density of $\pi/2$-congruent numbers is $0$ for $n\equiv1,2,3\pmod{8}$ by virtue his theorem about ${\rm Sel}_{2^\infty}(A_n)$ as following:

\begin{theorem}\cite[Corollary 1.2]{S17}
Take $E/\mathbb{Q}$ to be an elliptic curve with full rational 2-torsion. Assume that $E$ has no rational cyclic subgroup of order four. Then, for any
$N>1$, we have

$$\#\{1\leq d \leq N|\ corank{\rm Sel}_{2^\infty}(A_n)\geq \ 2\}=o(N).$$
\end{theorem}

 In this paper, our main research objects are $\pi/3$ and $2\pi/3$ congruent number elliptic curve, i.e, the "tiling number elliptic curve"(TNEC) and its $2$-Selmer group. Significantly, TNEC and CNEC have essential differences: CNEC have complex multiplication with full rational 2-torsion, then we have distribution of $2^\infty$-Selmer group \cite{S17}; TNEC haven't complex multiplication with 4-torsion, we haven't had a conjecture about distribution of $2^\infty$-Selmer group until now. 
 
 For each square-free positive integer number $n$, Monsky [HBM94, Appendix] proved the rank of 2-Selmer groups associated to $E_{n,\pi/2}$ is even if and only if $n\equiv1,2,3\pmod{8}$. Youcef Mokrani \cite{M20} skillfully used the tool he called "adapting Monsky matrix" to deal some problems with $\pi/3$ and $2\pi/3$-congruent number.

In this paper, at first, we get a similar result with Monsky about pairty of the rank of 2-Selmer groups associated to $E_{n,\pi/3}$ and $E_{n,2\pi/3}$ as following.

\begin{theorem}
	Let $s_2(E_{n,\theta})$ be the rank of $2$-Selmer groups of $E_{n,\theta}$, where $n\in\mathbb{Z}^+$ is squarefree. We have:
	\begin{align*}
		1)\ &s_2(E_{n,\pi/3})\text{ is even if and only if }n\equiv 1,2,3,5,7,9,14,15,19\pmod{24},\\
		&s_2(E_{n,\pi/3})\text{ is odd if and only if }n\equiv 6,10,11,13,17,18,21,22,23\pmod{24}.\\
		2)\ &s_2(E_{n,2\pi/3})\text{ is even if and only if }n\equiv 1,2,3,6,7,11,13,14,18\pmod{24},\\
		&s_2(E_{n,2\pi/3})\text{ is odd if and only if }n\equiv 5,9,10,15,17,19,21,22,23\pmod{24}.
	\end{align*}
	
\end{theorem}

In 1979, Goldfeld proposed his famous conjecture \cite{G79} as the following:

\begin{conjecture}
	Let $E$ be an elliptic curve over $\mathbb{Q}$.	Then, for a density one subset of square-free integers $n$ with $\epsilon(E_n)=+1$ $(resp.\ \epsilon(E_n)=-1)):$
	$$ord_{s=1}L(s,E_n)=0,\ \ \ (resp.\ ord_{s=1}L(s,E_n)=1).$$
	
\end{conjecture}

Where the root number $\epsilon(E_n)=+1(resp. -1)$ means the analytic rank is even(odd). We call the even part or odd part of Goldfeld's conjecture for root number is $+1$ or $-1$, respectively.

In this paper, we determine the tiling numbers by consider the $2^\infty$ part of class group of quadratic field $\mathbb{Q}(\sqrt{-n})$ like Burungale and Tian \cite{BAT22},  Ouyang \cite{OZ14} \cite{OZ15},  Youcef \cite{M20}, and deduce a positive density to support even part of  Goldfeld's conjecture:

\begin{corollary}
	Let $n\equiv 3,7,15,19\pmod{24}$ be a positive integer number and the $4$-rank $r_4$ of class group ${\rm Cl}(K)$ for quadratic field $K=\mathbb{Q}(\sqrt{-n})$ be $0$. Then $n$ is a non $\pi/3$-congruent number. 
\end{corollary}

\begin{corollary}
	Let $n\equiv 2,3,6,11,14,18\pmod{24}$ be a positive integer number, the $4$-rank $r_4$ of class group ${\rm Cl}(K)$ for quadratic field $K=\mathbb{Q}(\sqrt{-n})$ be $0$. Then $n$ is a non $2\pi/3$-congruent number. 
\end{corollary}

{\bf Remark}: The density $4$-rank $r_4=0$ of class group ${\rm Cl}(K)$ for quadratic field $K=\mathbb{Q}(\sqrt{-n})$ is equal to $\prod_{i=1}^{\infty}(1-\frac{1}{2^i})\approx28.87\%$. The case of non $\pi/3$-congruent number for $n\equiv 3,7\pmod{24}$ had been proved by Burunagle and Tian in Theorem 6.4 of \cite{BAT22}.  

Moreover, for the cases of Shafarevich-Tate groups being non-trivial, like Wang\cite{WZ16}, we have a infinite but zero-density result as following:

\begin{corollary}
	Let $n= p_1p_2...p_t\equiv  19\pmod{24}$. Assume the $4$-rank $r_4$ of class group ${\rm Cl}(K)$ for quadratic field $K=\mathbb{Q}(\sqrt{-n})$ is $1$, i.e., ${\rm Cl}(K)[2]\cap2{\rm Cl}(K)=\{1,[(d,\frac{-n+\sqrt{-n}}{2})]\}$, for some $d|n$.	 If 
	
	$1)$ $n$ have only one prime factor $p_t\equiv3\pmod{4}$ and $r_8=1-\left[\frac{6}{d}\right]$; 
	or
	
	$2)$ all prime factors of $d$ congruent to $1\pmod{8}$and $r_8=0$,
	
\noindent then $n$ is a non $\pi/3$-congruent number and $\Sha(E_n)[2]$ is non-trivial.
\end{corollary}

Actually, we give a better necessary condition than Corollary 1.9 for Shafarevich-Tate groups is non-trivial(this result can be saw in Theorem 6.5). However, the expression is very complicated. We guess it contains some information on positive densities.

In particular, Qin in \cite{Q22} discuss some special cases of congruent numbers for 2 prime factor i.e., $n=pq$, $p.q$ are prime, like:

\begin{theorem} \cite[Theorem 6.2]{Q22}
	Let $p\equiv5,\ q\equiv7\pmod{8}$ be two primes. If $[\frac{q}{p}]=1$ or $[\frac{q}{p}]=0$ and the class number $h(-pq)\equiv 4\pmod{8}$, then $L(A_{pq})\neq 0$, therefore rk$A_{pq}(\mathbb{Q})=0$ and $pq$ is not a congruent number.
\end{theorem}

 In the last section, we analogous give some examples of non $\pi/3$ or non $2\pi/3$-congruent numbers for $n=pq$ to support even Goldfeld's conjuncture:

\begin{theorem}
	For $n=pq\equiv 5\pmod{24}$, $p,q$ are prime, $[\frac{p}{q}]=1$ or $[\frac{p}{q}]=0$ and $[\frac{\beta}{q}]=1$ where $\beta$ satisfy $\beta^2=p\pmod{q}$ and $q\mid c+\beta a$, then $n$ is a non $\pi/3$-congruent number. Moreover, the density of $n$ is non $\pi/3$-congruent number is at least $75\%$. 
\end{theorem}

\begin{theorem}
	For $n=pq\equiv 11\pmod{24}$, $p,q$ are prime, $[\frac{p}{q}]=1$ or $[\frac{p}{q}]=0$ and $[\frac{\beta}{q}]=1$ where $\beta$ satisfy $\beta^2=p\pmod{q}$ and $q\mid c+\beta a$, then $n$ is a non $2\pi/3$-congruent number. Moreover, the density of $n$ is non $2\pi/3$-congruent number is at least $75\%$. 
\end{theorem}

\section{\bf Congruent elliptic curves and 2-Selmer groups}
For a real number $0<\theta<\pi$ stratifies cos $\theta=s/r$, $r,s\in\mathbb{Z},\ {\rm gcd}(r,s)=1$, $E_{\theta}$ be the elliptic curves defined by

$$y^2=x(x+(r+s))(x-(r-s)).$$
The quadratic twist family of $E_{\theta}$ is $E_{n,\theta}:y^2=x(x+n(r+s))(x-n(r-s))$ for $n\in\mathbb{Z}$ and square free.  In this paper, we consider the case of $\theta=\pi/3$ or $2\pi/3$, i.e., $s=\pm 1,\ r=2$. We denote $E_{n,\pi/3}$ by $E_{n}$ for simplifying notation,
Obviously, $E_{n,\pi/3}$ and $E_{-n,2\pi/3}$ represent the same elliptic curve. Hence we only  explore $E_{n,\pi/3}$ and denote $E_{n,\pi/3}$ by $E_n$ for simplifying notation.

In this paper, we let $\widetilde{n}=|\frac{n}{{\rm gcd}(6,n)}|=\prod_{i=1}^{t}p_i$, $p_i\neq2,3$. By a classical result (see \cite{S}), we know that ${\rm Sel}_2(E)$(that is the 2-Selmer group of $E:y^2=x(x-e_1)(x-e_2)$) can be written as
 $$\{\Lambda=(b_1,b_2)\in \mathbb{Q}(S,2)^2 | C_{\Lambda}(\mathbb{A})\neq\emptyset\},$$
where $e_1,e_2\in\mathbb{Z}$, $\mathbb{A}$ is the Ad$\grave{e}$le ring of $\mathbb{Q}$, $\mathbb{Q}(S,2)^2$ is the subgroup of $\mathbb{Q}^{\times}/(\mathbb{Q}^{\times})^2$ supported on $S$ which is the set of primes dividing $2e_1e_2(e_1-e_2)\infty$ and $C_{\Lambda}$ is a genus one curve in $\mathbb{P}^3$ defined by 

$$\left\{
\begin{aligned}
	H_1:&\ b_2 u_2^2-b_1b_2 u_3^2=(e_2-e_1)t^2, \\
	H_2:&\ b_1 u_1^2-b_1b_2 u_3^2=e_2t^2, \\
	H_3:&\ b_1 u_1^2-b_2 u_2^2=e_1t^2.
\end{aligned}
\right.
$$

For elliptic curves $E_n$, i.e., $e_1=n$, $e_2=-3n$, 2-torsion points $O,$ $(0,0),$ $(n,0),$ $(-3n,0)$ correspond to $(b_1,b_2)=(1,1),\ (-3,-n),\ (n,1)$, $(-3n,-n)$, respectively. For $\theta=\pi/3$ or $2\pi/3$, $b_1,b_2$ have the form by \cite{M20}:

\begin{equation}\tag{$\star$}
	 \begin{aligned}
	 	b_1=&(-1)^{\gamma_1}2^{\gamma_2}3^{\gamma_3}\prod_{i=1}^{t}p_i^{x_i}, \\
	 	b_2=&(-1)^{\xi_1}2^{\xi_2}3^{\xi_3}\prod_{i=1}^{t}p_i^{y_i},
	 \end{aligned}
\end{equation}
with $\gamma_i,\xi_i,x_i,y_i\in\{0,1\}.$ Then we can represent any pair $(b_1,b_2)$ by an unique vector $v\in\mathbb{F}_2^{2t+6}$ with

$$v^{T}=(\xi_1,\xi_2,\xi_3,\gamma_1,\gamma_2,\gamma_3,y_1,...,y_t,x_1,...,x_t).$$

\begin{theorem}$(\text{Hensel's Lemma})$\cite{Milne}
 Let $x$ be the closed point of $X=Spec\ A$, where $A$ is a local ring with maximal ideal $\mathcal{m}$ and residue field $k$. The following are equivalent:

$1)\ A$ is Henselian;

$2)$ let $f_1,...,f_n\in A[T_1,...,T_n]$; if there exists an $a=(a_1,...,a_n)\in k^n$ such that $\overline{f_i}(a)=0,i=1,2,...,n$, and ${\rm det}((\partial \overline{f_i}/\partial T_j)(a))\neq 0$, then there exists a $b\in A^n$ such that $\overline{b}=a$ and $f_i(b)=0,i=1,2,...,n$.
	
\end{theorem}

As we all know, every complete local ring $A$ is Henselian and $C_{\Lambda}(\mathbb{A})\neq \emptyset$ if and only if $C_{\Lambda}(\mathbb{Q}_S)\neq \emptyset$. Then we can calculate the ${\rm Sel}_2(E_n)$. An example in \cite{HHT21} show the condition for $(6,n)=1$ and $\widetilde{n}\equiv1({\rm mod}\ 8)$. In this case,

$C_{\Lambda}(\mathbb{R})\neq \emptyset$ if and only if 
\begin{equation*}
	\begin{cases}
		b_1b_2>0, &\text{if}\ n>0\\
		b_2>0, &\text{if}\ n<0 
	\end{cases},
\end{equation*}

For each $p|n$, $C_{\Lambda}(\mathbb{Q}_p)\neq \emptyset$ if and only if

\begin{equation*}
	\begin{cases}
			\left( \frac{b_1}{p}\right)=\left( \frac{b_2}{p}\right)=1, &\text{if}\ p \nmid b_1b_2\\
				\left( \frac{nb_1}{p}\right)=\left( \frac{b_2}{p}\right)=1, &\text{if}\ p\ |\ b_1,p\nmid b_2 \\
							\left( \frac{-3b_1}{p}\right)=\left( \frac{-nb_2}{p}\right)=1, &\text{if}\ p\nmid b_1,p\ |\  b_2\\
				\left( \frac{-3nb_1}{p}\right)=\left( \frac{-nb_2}{p}\right)=1, &\text{if}\ p\ |\ b_1,p\ |\  b_2 
		\end{cases},
\end{equation*}

$C_{\Lambda}(\mathbb{Q}_3)\neq \emptyset$ if and only if $3\nmid b_2$ and
\begin{equation*}
	\begin{cases}
		\left( \frac{b_2}{3}\right)=1, &\text{if}\ 3 \nmid b_1\\
		\left( \frac{-nb_2}{3}\right)=1, &\text{if}\ 3\ |\ b_1 
	\end{cases},
\end{equation*}

$C_{\Lambda}(\mathbb{Q}_2)\neq \emptyset$ if and only if $2\nmid b_2$, $(b_1\ {\rm mod}\ 8,b_2\ {\rm mod}\ 4)=(1,1)$ or $(5,3)$; or $2\ |\ b_2$, $(b_1\ {\rm mod}\ 8,b_2\ {\rm mod}\ 8)=(7,6)$ or $(3,2)$.

\section{\bf Narrow class group and R$\acute{e}$dei Matrix}

Let $d\neq 1$ be a squarefree integer, $K=\mathbb{Q}(\sqrt{d})$ be the corresponding quadratic field, $D=d$ if $d\equiv1(\text{mod 4})$ or $D=4d$ for others, and ${\rm Cl}^+(d):={\rm Cl}^+(\mathcal{O}_K)$ the narrow class group of $K$, i.e., the quotient of the fractional ideal group $I$ by the subgroup of principal ideals $(x)=x\mathcal{O}_K$ with generator of positive norm $N(x)$. If $K$ is real quadratic field and the norm of fundamental unit $\epsilon_d$ is $1$, ${\rm Cl}^+(d)$ has twice the size of ${\rm Cl}(\mathcal{O}_K)$. In other cases, ${\rm Cl}^+(d)={\rm Cl}(\mathcal{O}_K)$. For $k\in\mathbb{Z}^+$, we denote by $r_{2^k}(d)$ the $2^k$-rank of ${\rm Cl}^+(d)$, i.e., $r_{2^k}(d):={\rm dim}_{\mathbb{F}_2}2^{k-1}{\rm Cl}^+(d)/2^k{\rm Cl}^+(d)$.

Gauss's genus theory shows that $r_2(d)=t-1$ with $t$ the number of prime divisors of $D$. R$\acute{e}$dei\cite{RR33}, Waterhouse\cite{W73}, and Kolster\cite{K05} show $r_{2^{k+1}}(d)$ can be expressed by the $2$-rank of R$\acute{e}$dei Matrix $R^{(k)}(k\geq1)$ whose entries are certain Hilbert symbols.

\begin{theorem}
	$r_{2^{k+1}}(d)=t-1-rk(R^{(k)}(d)).$
\end{theorem}

For $k=1$, the R$\acute{e}$dei Matrix of quadratic field $K=\mathbb{Q}(\sqrt{d})$ is 

$$R^{(1)}(d)=\begin{pmatrix}
	[p_1,d]_{p_1} & [p_2,d]_{p_1} & \cdots & [p_r,d]_{p_1} \\
	[p_1,d]_{p_2} & [p_2,d]_{p_2} & \cdots & [p_r,d]_{p_2} \\
	\vdots & \vdots & \ddots & \vdots \\
	[p_1,d]_{p_r} & [p_2,d]_{p_r} & \cdots & [p_r,d]_{p_r} 
\end{pmatrix},$$
where $\{p_1,p_2,\cdots,p_r\}$ is the all prime numbers ramify in $K$, and $[\cdot,\cdot]_p$ is the additive Hilbert symbol(that is, $[\cdot,\cdot]_p=\frac{1}{2}(1-(\cdot,\cdot)_p)$). Then $R^{(1)}(d)$ is a $\mathbb{F}_2$-matrix. In this paper, we denote $R^{(1)}(d)$ by $R(d)$ for simplification.

\section{\bf Monsky Matrix}
 The additive Legendre symbol $[\frac{n}{m}]:=\frac{1}{2}(1-\left(\frac{n}{m}\right))\in \mathbb{F}_2$ is additive, i.e. $[\frac{ab}{n}]=[
 \frac{a}{n}]+[\frac{b}{n}]$.The law of quadratic reciprocity of additive Legendre symbol is  $[\frac{p}{q}]+[\frac{q}{p}]=[\frac{-1}{p}][\frac{-1}{q}]$ for primes $p,q$. It is worth noting that the equivalence condition of $C_{\Lambda}(\mathbb{A})\neq \emptyset$ can be rephrased as linear equation in $\mathbb{F}_2$ by the additive Legendre symbol. 
 
 For a matrix $U$, let $r(U)$ be the sum of the rows of $U$, and $c(U)$ be the sum of the columns of $U$. Note that $r_{d}:=\left([\frac{d}{p_1}],...,[\frac{d}{p_t}]\right)\in \mathbb{F}_2^t$ and $D_d:=diag(r_{d})\in M_{t\times t}(\mathbb{F}_2)$ for an integer $d$ coprime to $\widetilde{n}$. It is easy to see that $r(D_d)=r_d$ and $c(D_d)=r_d^T$. Let $O$ be any zero matrix or vector. Some symbols are same with Monsky in \cite[Appendix]{HBM94}.
 
 Let $A=(a_{ij})\in M_{t\times t}(\mathbb{F}_2)$ be the R$\acute{e}$dei Matrix $R(-\widetilde{n})$ associated to $\widetilde{n}$ where $a_{ij}=[\frac{p_j}{p_i}]$ for $i\neq j$ and $a_{ii}=\sum\limits_{j\neq i}a_{ij}=[\frac{\widetilde{n}p_i}{p_i}]$. We have $r(A)=(1+[\frac{-1}{\widetilde{n}}])r_{-1}$, $c(A)=0$, and 
 $$A+A^T=D_{-1}+r_{-1}^Tr_{-1}.$$
 
  In this paper, let $v_p(d)$ denote the $p$-adic valuation of $d$, for $p\neq \infty$ and  $v_\infty(d):=sign(d)$ the signum function.
 
 Now we define the Monsky matrix:
 
 \begin{definition}
 	For elliptic curve $E$ over number field $K$, we let full Monsky matrix be a matrix $M$ over $F_2$ associated to $E$ if we have a bijection
 	$$\phi:{\rm Sel}_2(E)\rightarrow\{v\in\mathbb{F}_2^{2s}|Mv=0\}$$
 	given by $(b_1,b_2)\rightarrow(v_{p_1}(d_1),...,v_{p_s}(d_1),v_{p_1}(d_2),...,v_{p_s}(d_2))$, where $p_i\in S$, $s=\# S$, and $S$ is a finite prime place set.
 \end{definition}

In this paper, we will show:

\begin{theorem}
 For any non-zero integer $n$, the full Monsky matrix $M_n$ associated to $E_n$ is exist.
\end{theorem}

And we list all Monsky matrices $M_n$ in Section 6. For example, let $n>0$, $(6,n)=1$ and  $\widetilde{n}\equiv 1\pmod 8$, the full Monsky matrix 

$$M_n=\begin{pmatrix}
	1 & 0 & 0 & 1 & 0 & 0 & O & O\\
	0 & 1 & 0 & 1 & 0 & 1 & O & r_{-1}\\
	0 & 0 & 0 & 0 & 1 & 0 & O & O\\
	1 & 1 & 1 & 0 & 0 & 1 & r_{-1} & r_2\\
	1 & 1 & 0 & 0 & 0 & 1+[\frac{-3}{\widetilde{n}}] & r_{-3} & O\\
	0 & 0 & 1 & 0 & 0 & 0 & O & O\\
	0 & 0 & 0 & r_{-1}^T & r_{2}^T & r_{3}^T & D_{-3} & A\\
	r_{-1}^T & r_{2}^T & r_{3}^T & 0 & 0 & 0 & A+D_{-1} & O
\end{pmatrix}.$$

The rank $s_2(E_n)$ of 2-Selmer group ${\rm Sel}_2(E_n)$ is same as the dimension of solution space $V$ of $M_nv=0$. Then we have $s_2(E_n)=2n+6-rk(M_n)$.

\section{\bf Cassels Pairing}

Cassels\cite{C62} constructed a pairing $\left \langle \cdot,\cdot \right \rangle$ on $\Sha(E/K)$ for elliptic curves $E$ over $K$ in 1962.
\begin{theorem}
	Let $K$ be an algebraic number field, $E$ an abelian variety of dimension $1$  and $\Sha$ the corresponding Tate-Shafarevich group of classes of principal homogeneous spaces for $E$ defined over $K$ which are everywhere locally trivial. Then there exists a naturally defined skew-symmetric form $\left \langle \Lambda,\Pi \right \rangle$ defined for $\Lambda,\Pi \in \Sha$ and taking values in the group $\mathbb{Q}/\mathbb{Z}$ of the rationals modulo $1$ which has the following property:
	
	Let $m$ be a natural number and suppose that $\left \langle \Lambda,\Pi \right \rangle=0$ for all $\Lambda\in \Sha$ such that $m\Lambda=0$; then $\Pi=m\Pi^{\prime}$ for some $\Pi^{\prime}\in \Sha$.
\end{theorem}

Moreover, on the $\mathbb{F}_2$-vector space ${\rm Sel}_2(E)/E(\mathbb{Q})[2]$, the Cassels pairing $\left \langle \cdot,\cdot \right \rangle$ which is a natural lifting of on $\Sha(E/K)[2]$ is computable by Cassels[C98]. In this paper, we use "addition" instead of "multiplication" in the process of calculating Cassels pairing, i.e., acting a group isomorphism  $\{\pm1\}\rightarrow \mathbb{F}_2$ by mapping $1$ to $0$ and map $-1$ to $1$ for value. By Cassels, we have $\left \langle \cdot,\cdot \right \rangle=\sum_p\left \langle \cdot,\cdot \right \rangle_p$ where almost all local Cassels pairings $\left \langle \cdot,\cdot \right \rangle_p$ are trivial. 

By this theorem, if $\left \langle \Lambda,\Pi \right \rangle\neq0$ for $\Lambda,\Pi\in {\rm Sel}_2(E)/E(\mathbb{Q})[2]$, we have $\Lambda,\Pi\in \Sha(E/K)[2]$. We can see Zhangjie Wang\cite{WZ16} for more details about calculating Cassels pairing.

\section{\bf Main result}

 Let $\eta=n/\widetilde{n}$. We use $\dot{m}:=\frac{m}{{\rm gcd}(m,p)}$ for simplification when working in $\mathbb{Z}/q\mathbb{Z}$ where $q$ is a power of prime $p$ and $m$ is an integer. In this section, we write all Monsky matrices for each square-free $n\in\mathbb{Z}$ module $24$, and calculate the parity of ${\rm Sel}_2(E_n)$, respectively. 
 
 In our proof, we use a fact coming linear algebra that a skew-symmetric matrix have even rank and a density result for $2^k$-rank of class group of quadratic field as following:

\begin{theorem}\cite{FK07}
	Let $k$ be a non negative integer. Then 

$1)$ The density of negative fundamental discriminant $D$ such that $r_4(D)=k$ is equal to
$$\frac{\prod\limits_{i=1}^{\infty}(1-2^{-i})}{2^{k^2}\prod\limits_{i=1}^{k}(1-2^{-i})^2},$$

$2)$ The density of positive fundamental discriminant $D$ such that $r_4(D)=k$ is equal to
$$\frac{\prod\limits_{i=1}^{\infty}(1-2^{-i})}{2^{k(k+1)}\prod\limits_{i=1}^{k}(1-2^{-i})\prod\limits_{i=1}^{k+1}(1-2^{-i})}.$$
\end{theorem}

For $k=0$, the density of negative fundamental discriminant is equal to $\prod_{i=1}^{\infty}(1-\frac{1}{2^i})\approx28.87\%$; the density of positive fundamental discriminant is equal to $\frac{1}{2}\prod_{i=1}^{\infty}(1-\frac{1}{2^i})\approx14.43\%$.

\subsection{\bf $\eta=1$} Let 
$$A_1=\begin{pmatrix}
	1 & 0 & 0 & 1 & 0 & 0 & O & O\\
	0 & 1 & 0 & 1 & 0 & 1 & O & r_{-1}\\
	0 & 0 & 0 & 0 & 1 & 0 & O & O\\
	1 & 1 & 1 & 0 & 0 & 1 & r_{-1} & r_2\\
	1 & 1 & 0 & 0 & 0 & 1+[\frac{-3}{\widetilde{n}}] & r_{-3} & O\\
	0 & 0 & 1 & 0 & 0 & 0 & O & O\\
	O & O & O & r_{-1}^T & r_{2}^T & r_{3}^T & D_{-3} & A\\
	r_{-1}^T & r_{2}^T & r_{3}^T & O & O & O & A+D_{-1} & O
\end{pmatrix},$$

and 
$$A_2=\begin{pmatrix}
	1 & 0 & 0 & 1 & 0 & 0 & O & O\\
	0 & 0 & 0 & 0 & 1 & 0 & O & O\\
	0 & 1 & 0 & 0 & 0 & 0 & O & O\\
	0 & 0 & 0 & 1 & 0 & 1 & O & r_{-1}\\
	1 & 1 & 0 & 0 & 0 & 1+[\frac{-3}{\widetilde{n}}] & r_{-3} & O\\
	0 & 0 & 1 & 0 & 0 & 0 & O & O\\
	O & O & O & r_{-1}^T & r_{2}^T & r_{3}^T & D_{-3} & A\\
	r_{-1}^T & r_{2}^T & r_{3}^T & O & O & O & A+D_{-1} & O
\end{pmatrix},$$

and
$$A_3=\begin{pmatrix}
	1 & 0 & 0 & 1 & 0 & 0 & O & O\\
	0 & 1 & 0 & 0 & 0 & 0 & O & O\\
	0 & 0 & 0 & 0 & 1 & 0 & O & O\\
	1 & 0 & 1 & 0 & 0 & 0 & r_{-1} & O\\
	1 & 1 & 0 & 0 & 0 & 1+[\frac{-3}{\widetilde{n}}] & r_{-3} & O\\
	0 & 0 & 1 & 0 & 0 & 0 & O & O\\
	O & O & O & r_{-1}^T & r_{2}^T & r_{3}^T & D_{-3} & A\\
	r_{-1}^T & r_{2}^T & r_{3}^T & O & O & O & A+D_{-1} & O
\end{pmatrix}.$$

\begin{theorem}
	For $n=\widetilde{n}=p_1...p_t,p_i \neq 2,3$, then we have Monsky matrix 
\begin{equation*}
	M_n=\begin{cases}
		A_1, &\text{if } \widetilde{n}\equiv 1({\rm mod }\ 8),\\
		A_2, &\text{if } \widetilde{n}\equiv 5({\rm mod }\ 8), \\
		A_3, &\text{if } \widetilde{n}\equiv 3({\rm mod }\ 4).
	\end{cases}
\end{equation*} 
\end{theorem}

\emph{\bf Proof:} We list all the equivalent conditions and transform them to additional Legendre symbols to build system of linear equations. 

We know that $C_{\Lambda}(\mathbb{A})\neq \emptyset$ if and only if $C_{\Lambda}(\mathbb{Q}_S)\neq \emptyset$, i.e., 

1) $C_{\Lambda}(\mathbb{R})\neq \emptyset$ if and only if $b_1b_2>0$ that is equivalent to $\gamma_1+\xi_1=0$;

2) $C_{\Lambda}(\mathbb{Q}_{p_i})\neq \emptyset$ for $i=1,..,t$ if and only if 
\begin{equation*}
	\begin{cases}
		\left( \frac{\dot{b}_1}{p_i}\right)=\left( \frac{\dot{b}_2}{p_i}\right)=1, &\text{if}\ p_i \nmid b_1b_2\\
		\left( \frac{\dot{n}\dot{b}_1}{p_i}\right)=\left( \frac{\dot{b}_2}{p_i}\right)=1, &\text{if}\ p_i\ |\ b_1,p_i\nmid b_2 \\
		\left( \frac{-3\dot{b}_1}{p_i}\right)=\left( \frac{-\dot{n}\dot{b}_2}{p_i}\right)=1, &\text{if}\ p_i\nmid b_1,p_i\ |\  b_2\\
		\left( \frac{-3\dot{n}\dot{b}_1}{p_i}\right)=\left( \frac{-\dot{n}\dot{b}_2}{p_i}\right)=1, &\text{if}\ p_i\ |\ b_1,p_i\ |\  b_2 
	\end{cases}
\end{equation*}

by Hensel's Lemma. The conditions are equivalent to
\begin{equation*}
	\begin{cases}
		[\frac{\dot{b}_1}{p_i}]+[\frac{-3}{p_i}]y_i+[\frac{\widetilde{n}p_i}{p_i}]x_i=0 \\
		[\frac{\dot{b}_2}{p_i}]+[\frac{-\widetilde{n}p_i}{p_i}]y_i=0 
	\end{cases};
\end{equation*}
In other words, we have
\begin{equation*}
	\begin{cases}
[\frac{-1}{p_i}]\gamma_1+[\frac{2}{p_i}]\gamma_2+[\frac{3}{p_i}]\gamma_3+[\frac{-3}{p_i}]y_i+\sum\limits_{j\neq i}[\frac{p_j}{p_i}]x_j+[\frac{\widetilde{n}p_i}{p_i}]x_i=0 \\
[\frac{-1}{p_i}]\xi_1+[\frac{2}{p_i}]\xi_2+[\frac{3}{p_i}]\xi_3+\sum\limits_{j\neq i}[\frac{p_j}{p_i}]y_j+[\frac{-\widetilde{n}p_i}{p_i}]y_i=0 
	\end{cases};
\end{equation*}

3) $C_{\Lambda}(\mathbb{Q}_3)\neq \emptyset$ if and only if $3\nmid b_2$ and
\begin{equation*}
	\begin{cases}
		\left( \frac{\dot{b}_2}{3}\right)=1, &\text{if}\ 3 \nmid b_1\\
		\left( \frac{-\dot{n}\dot{b}_2}{3}\right)=1, &\text{if}\ 3\ |\ b_1 
	\end{cases}
\end{equation*}

by Hensel's Lemma. The conditions are equivalent to

\begin{equation*}
	\begin{cases}
		(1+[\frac{-3}{\widetilde{n}}])\gamma_3+[\frac{-3}{\dot{b}_2}]=0\\
	\xi_3=0
	\end{cases};
\end{equation*}
In other words, we have 
\begin{equation*}
	\begin{cases}
		\xi_1+\xi_2+(1+[\frac{-3}{\widetilde{n}}])\gamma_3+\sum\limits_{i}[\frac{-3}{p_i}]y_i=0\\
		\xi_3=0
	\end{cases};
\end{equation*}
4) $C_{\Lambda}(\mathbb{Q}_2)\neq \emptyset$. If  $\widetilde{n}\equiv 1\pmod{8}$, $C_{\Lambda}(\mathbb{Q}_2)\neq \emptyset$ if and only if (i):$2\nmid b_1b_2$ and $(\dot{b}_1,\dot{b}_2)\equiv (1,1),(1,5),(5,7),(5,3)({\rm mod }\ 8,{\rm mod }\ 8)$ or (ii):$2\nmid b_1,2\ |\ b_2$ and $(\dot{b}_1,\dot{b}_2)\equiv (7,3),$ $(3,1)({\rm mod }\ 8,{\rm mod }\ 4)$;
 if $\widetilde{n}\equiv 5\pmod{8}$, $C_{\Lambda}(\mathbb{Q}_2)\neq \emptyset$ if and only if $2\nmid b_1b_2$ and $\dot{b}_1\equiv 1({\rm mod }\ 4)$; 
 if $\widetilde{n}\equiv 3\pmod{4}$, $C_{\Lambda}(\mathbb{Q}_2)\neq \emptyset$ if and only if $2\nmid b_1b_2$ and $\dot{b}_2\equiv 1({\rm mod }\ 4)$ by Hensel's Lemma. The conditions are equivalent to
 
 \begin{equation*}
 	\begin{cases}
 		\xi_2+\gamma_1+\gamma_3+\sum\limits_{i}[\frac{-1}{p_i}]x_i=0, \ \ \ \ \gamma_2=0,\\
 		\xi_1+\xi_2+\xi_3+\gamma_3+\sum\limits_{i}[\frac{-1}{p_i}]y_i+\sum\limits_{i}[\frac{2}{p_i}]x_i=0
 	\end{cases} \text{if } \widetilde{n}\equiv 1\pmod{8},
 \end{equation*}
and
 \begin{equation*}
	\begin{cases}
		\xi_2=0, \ \ \ \ \gamma_2=0,\\
		\gamma_1+\gamma_3+\sum\limits_{i}[\frac{-1}{p_i}]x_i=0
	\end{cases} \text{if } \widetilde{n}\equiv 5\pmod{8},
\end{equation*}
and
 \begin{equation*}
	\begin{cases}
		\xi_2=0, \ \ \ \ \gamma_2=0,\\
		\xi_1+\xi_3+\sum\limits_{i}[\frac{-1}{p_i}]y_i=0
	\end{cases} \text{if } \widetilde{n}\equiv 3\pmod{4}.
\end{equation*}
The Monsky matrix can be naturally derived from the above system of linear equations.
$\Box$

Then we have a corollary about pairty of 2-rank as following:

\begin{corollary}
	For $n=\widetilde{n}=p_1...p_t,p_i\neq2,3$, the rank of $2$-Selmer group $s_2(E_n)$ is even (resp. odd) if and only if $n\equiv1,5,7,19(resp.\ 11,13,17,23)\pmod{24}$. 
\end{corollary}

\emph{\bf Proof:}
Note that $s_2(E_n)=2t+6-r(M_n)$, then we only need consider the rank of $M_n$ which is denoted by $r(M_n)$.

For $\widetilde{n}\equiv1(\text{mod 8})$, we have:
\begin{align*}
	r(A_1)&=r\begin{pmatrix}
		1 & 0 & 0 & 1 & 0 & 0 & O & O\\
		0 & 1 & 0 & 1 & 0 & 1 & O & r_{-1}\\
		0 & 0 & 0 & 0 & 1 & 0 & O & O\\
		1 & 1 & 1 & 0 & 0 & 1 & r_{-1} & r_2\\
		1 & 1 & 0 & 0 & 0 & 1+[\frac{-3}{\widetilde{n}}] & r_{-3} & O\\
		0 & 0 & 1 & 0 & 0 & 0 & O & O\\
		O & O & O & r_{-1}^T & r_{2}^T & r_{3}^T & D_{-3} & A\\
		r_{-1}^T & r_{2}^T & r_{3}^T & O & O & O & A+D_{-1} & O
	\end{pmatrix} \\
&=r\begin{pmatrix}
	1 & 0 & 0 & 1 & 0 & 0 & O & O\\
	0 & 1 & 0 & 1 & 0 & 1 & O & r_{-1}\\
	0 & 0 & 0 & 0 & 1 & 0 & O & O\\
	0 & 0 & 0 & 0 & 0 & 0 & r_{-1} & r_{-2}\\
	0 & 0 & 0 & 0 & 0 & [\frac{-3}{\widetilde{n}}] & r_{-3} & r_{-1}\\
	0 & 0 & 1 & 0 & 0 & 0 & O & O\\
	O & r_{-1}^T & O & O & O & r_{-3}^T & D_{-3} & A+r_{-1}^Tr_{-1}\\
	O & r_{-2}^T & O & O & O & r_{-1}^T & A+D_{-1} & r_{-1}^Tr_{-1}
\end{pmatrix} \\ &\text{  by row transformations.} \\
&=r\begin{pmatrix}
	1 & 0 & 0 & 0 & 0 & 0 & O & O\\
	0 & 0 & 0 & 1 & 0 & 0 & O & O\\
	0 & 0 & 0 & 0 & 1 & 0 & O & O\\
	0 & 0 & 0 & 0 & 0 & 0 & r_{-1} & r_{-2}\\
	0 & 0 & 0 & 0 & 0 & [\frac{-3}{\widetilde{n}}] & r_{-3} & r_{-1}\\
	0 & 0 & 1 & 0 & 0 & 0 & O & O\\
	O & r_{-1}^T & O & O & O & r_{-3}^T & D_{-3} & A+r_{-1}^Tr_{-1}\\
	O & r_{-2}^T & O & O & O & r_{-1}^T & A+D_{-1} & r_{-1}^Tr_{-1}
\end{pmatrix} \\ &\text{  by column transformations.} \\
&=4+r\begin{pmatrix}
	[\frac{-3}{\widetilde{n}}] & 0  & r_{-3} & r_{-1}\\
	0 & 0  & r_{-1} & r_{-2}\\
	r_{-3}^T & r_{-1}^T & D_{-3} & A^T+D_{-1}\\
	r_{-1}^T &r_{-2}^T & A+D_{-1} & r_{-1}^Tr_{-1}
\end{pmatrix} \\
 &\text{  by } A+A^T+r_{-1}^Tr_{-1}+D_{-1}=O. 
\end{align*}

If $[\frac{-3}{\widetilde{n}}]=0$, i.e. $n\equiv1({\rm mod }\ 24)$, $$r(A_1)=4+r\begin{pmatrix}
	0 & 0  & r_{-3} & r_{-1}\\
	0 & 0  & r_{-1} & r_{-2}\\
	r_{-3}^T & r_{-1}^T & D_{-3}+r_{-3}^Tr_{-3} & A^T+D_{-1}+r_{-3}^Tr_{-1}\\
	r_{-1}^T &r_{-2}^T & A+D_{-1}+r_{-1}^Tr_{-3} & O
\end{pmatrix} $$ is even; and if $[\frac{-3}{\widetilde{n}}]=1$, i.e. $n\equiv17({\rm mod }\ 24)$, $$r(A_1)=5+r\begin{pmatrix}
 0  & r_{-1} & r_{-2}\\
 r_{-1}^T & D_{-3}+r_{-3}^Tr_{-3} & A^T+D_{-1}+r_{-3}^Tr_{-1}\\
 r_{-2}^T & A+D_{-1}+r_{-1}^Tr_{-3} & O
\end{pmatrix} $$ is odd.

For $\widetilde{n}\equiv5(\text{mod 8})$, we have:
\begin{align*}
	r(A_2)&=r\begin{pmatrix}
		1 & 0 & 0 & 1 & 0 & 0 & O & O\\
		0 & 0 & 0 & 0 & 1 & 0 & O & O\\
		0 & 1 & 0 & 0 & 0 & 0 & O & O\\
		0 & 0 & 0 & 1 & 0 & 1 & O & r_{-1}\\
		1 & 1 & 0 & 0 & 0 & 1+[\frac{-3}{\widetilde{n}}] & r_{-3} & O\\
		0 & 0 & 1 & 0 & 0 & 0 & O & O\\
		O & O & O & r_{-1}^T & r_{2}^T & r_{3}^T & D_{-3} & A\\
		r_{-1}^T & r_{2}^T & r_{3}^T & O & O & O & A+D_{-1} & O
	\end{pmatrix} \\
	&=r\begin{pmatrix}
		1 & 0 & 0 & 1 & 0 & 0 & O & O\\
		0 & 0 & 0 & 0 & 1 & 0 & O & O\\
		0 & 1 & 0 & 0 & 0 & 0 & O & O\\
		0 & 0 & 0 & 1 & 0 & 1 & O & r_{-1}\\
		0 & 0 & 0 & 0 & 0 & [\frac{-3}{\widetilde{n}}] & r_{-3} & r_{-1}\\
		0 & 0 & 1 & 0 & 0 & 0 & O & O\\
		O & O & O & O & O & r_{-3}^T & D_{-3} & A+r_{-1}^Tr_{-1}\\
		O & O & O & O & O & r_{-1}^T & A+D_{-1} & r_{-1}^Tr_{-1}
	\end{pmatrix} \\ &\text{  by row transformations.} \\
	&=r\begin{pmatrix}
		1 & 0 & 0 & 0 & 0 & 0 & O & O\\
		0 & 0 & 0 & 0 & 1 & 0 & O & O\\
		0 & 1 & 0 & 0 & 0 & 0 & O & O\\
		0 & 0 & 0 & 1 & 0 & 0 & O & O\\
		0 & 0 & 0 & 0 & 0 & [\frac{-3}{\widetilde{n}}] & r_{-3} & r_{-1}\\
		0 & 0 & 1 & 0 & 0 & 0 & O & O\\
		O & O & O & O & O & r_{-3}^T & D_{-3} & A+r_{-1}^Tr_{-1}\\
		O & O & O & O & O & r_{-1}^T & A+D_{-1} & r_{-1}^Tr_{-1}
	\end{pmatrix} \\ &\text{  by column transformations.} \\
	&=5+r\begin{pmatrix}
		[\frac{-3}{\widetilde{n}}]  & r_{-3} & r_{-1}\\
		r_{-3}^T & D_{-3} & A^T+D_{-1}\\
		r_{-1}^T & A+D_{-1} & r_{-1}^Tr_{-1}
	\end{pmatrix} \\
	&\text{  by } A+A^T+r_{-1}^Tr_{-1}+D_{-1}=O. 
\end{align*}

If $[\frac{-3}{\widetilde{n}}]=0$, i.e. $n\equiv13\pmod{24}$, $$r(A_2)=5+r\begin{pmatrix}
	0  & r_{-3} & r_{-1}\\
	r_{-3}^T & D_{-3}+r_{-3}^Tr_{-3} & A^T+D_{-1}+r_{-3}^Tr_{-1}\\
	r_{-1}^T & A+D_{-1}+r_{-1}^Tr_{-3} & O
\end{pmatrix} $$ is odd; and if $[\frac{-3}{\widetilde{n}}]=1$, i.e. $n\equiv5\pmod{24}$, $$r(A_2)=6+r\begin{pmatrix} D_{-3}+r_{-3}^Tr_{-3} & A^T+D_{-1}+r_{-3}^Tr_{-1}\\
A+D_{-1}+r_{-1}^Tr_{-3} & O
\end{pmatrix} $$ is even.

For $\widetilde{n}\equiv3\pmod{4}$, we have:
\begin{align*}
	r(A_3)&=r\begin{pmatrix}
		1 & 0 & 0 & 1 & 0 & 0 & O & O\\
		0 & 1 & 0 & 0 & 0 & 0 & O & O\\
		0 & 0 & 0 & 0 & 1 & 0 & O & O\\
		1 & 0 & 1 & 0 & 0 & 0 & r_{-1} & O\\
		1 & 1 & 0 & 0 & 0 & 1+[\frac{-3}{\widetilde{n}}] & r_{-3} & O\\
		0 & 0 & 1 & 0 & 0 & 0 & O & O\\
		O & O & O & r_{-1}^T & r_{2}^T & r_{3}^T & D_{-3} & A\\
		r_{-1}^T & r_{2}^T & r_{3}^T & O & O & O & A+D_{-1} & O
	\end{pmatrix} \\
	&=r\begin{pmatrix}
	1 & 0 & 0 & 1 & 0 & 0 & O & O\\
	0 & 0 & 0 & 0 & 1 & 0 & O & O\\
	0 & 1 & 0 & 0 & 0 & 0 & O & O\\
	0 & 0 & 0 & 1 & 0 & 0 & r_{-1} & O\\
	0 & 0 & 0 & 0 & 0 & 1+[\frac{-3}{\widetilde{n}}] & r_{3} & O\\
	0 & 0 & 1 & 0 & 0 & 0 & O & O\\
	O & O & O & O & O & r_{3}^T & D_{-3}+r_{-1}^Tr_{-1} & A\\
	O & O & O & O & O & O & A+D_{-1}+r_{-1}^Tr_{-1} & O
\end{pmatrix} \\ &\text{  by row transformations.} \\
	&=r\begin{pmatrix}
	1 & 0 & 0 & 0 & 0 & 0 & O & O\\
	0 & 0 & 0 & 0 & 1 & 0 & O & O\\
	0 & 1 & 0 & 0 & 0 & 0 & O & O\\
	0 & 0 & 0 & 1 & 0 & 0 & O & O\\
	0 & 0 & 0 & 0 & 0 & 1+[\frac{-3}{\widetilde{n}}] & r_{3} & O\\
	0 & 0 & 1 & 0 & 0 & 0 & O & O\\
	O & O & O & O & O & r_{3}^T & D_{-3}+r_{-1}^Tr_{-1} & A\\
	O & O & O & O & O & O & A+D_{-1}+r_{-1}^Tr_{-1} & O
\end{pmatrix} \\ &\text{  by column transformations.} \\
	&=5+r\begin{pmatrix}
		[\frac{3}{\widetilde{n}}]  & r_{3} & O\\
		r_{3}^T & D_{-3}+r_{-1}^Tr_{-1} & A\\
		O & A^T & O
	\end{pmatrix} \\
	&\text{  by } A+A^T+r_{-1}^Tr_{-1}+D_{-1}=O. 
\end{align*}

If $[\frac{3}{\widetilde{n}}]=0$, i.e. $n\equiv11,23(\text{mod }24)$, $$r(A_3)=5+r\begin{pmatrix}
	0  & r_{3} & O\\
	r_{3}^T & D_{-3}+r_{-1}^Tr_{-1}+r_{3}^Tr_{3} & A\\
	O & A^T & O
\end{pmatrix} $$ is odd; and if $[\frac{3}{\widetilde{n}}]=1$, i.e. $n\equiv7,19(\text{mod }24)$, $$r(A_3)=6+r\begin{pmatrix}
 D_{-3}+r_{-1}^Tr_{-1}+r_{3}^Tr_{3} & A\\
 A^T & O
\end{pmatrix} $$ is even. $\Box$

Now we prove a result about non $\pi/3$-congruent numbers:

\begin{corollary}
	For $n=\widetilde{n}=p_1...p_t\equiv7,19\pmod{24}$ and $r_4(-n)=0$, then $s_2(E_n)=2$. 
\end{corollary}

\emph{\bf Proof:} 
For $n\equiv7,19(\text{mod }24)$, we have

$$ s_2(E_n)=2t+6-r(M_n)=2t-r\begin{pmatrix}
	D_{-3}+r_{-1}^Tr_{-1}+r_{3}^Tr_{3} & A\\
	A^T & O
\end{pmatrix}.$$
 Then $s_2\leq2t-2r(A)=2t-2r(R(-n))=2+2r_4(-n)=2.$

Note that $s_2(E_n)\geq2$ by the fact that $E_n$ has $4$ elements of $2$-torsion, we have $s_2(E_n)=2$. $\Box$

\begin{theorem}
	Let $n=\widetilde{n}=p_1...p_t\equiv19\pmod{24},p_i\neq2,3$ and $r_4(-n)=1$, i.e., ${\rm Cl}^+(K)[2]\cap2{\rm Cl}^+(K)=\{1,[(d,\frac{-n+\sqrt{-n}}{2})]\}$ for some $d|n$. Then we have $s_2(E_n)=4$. Moreover, let $v_c=(x_1(D_{-2}+r_{-1}^Tr_{-1}+r_{2}^T r_{2}+r_{6}^T r_{6})+r_c,x_1(D_{-2}+r_{-1}^Tr_{-1}+r_{2}^T r_{2}))\in \mathbb{F}_2^{2t}$, $M^{\prime}=$ $\begin{pmatrix}
		D_{-3}+r_{-1}^Tr_{-1}+r_{3}^Tr_{3} & A\\
		A^T & O
	\end{pmatrix}$, where $x_1=(x_{1,1},x_{1,2}...,x_{1,t})$ such that $d=\prod\limits_{i=1}^t p_i^{x_{1,i}}$, and $c$ satisfies $4c^2=da^2+n/d b^2$ for some $a,b\in \mathbb{Z}$. If $M^{\prime} {v^{\prime}}^T=v_c^T$ no solution, then $rank_{\mathbb{Z}}E_n(\mathbb{Q})=0$ and $\Sha(E_n)=(\mathbb{Z}/2\mathbb{Z})^2$ .
\end{theorem}

\emph{\bf Proof:} 
We define $\mathcal{p}^v :=\prod\limits_{i=1}^t p_i^{v_i}$ where $v=(v_1,v_2,...,v_t)\in \mathbb{F}_2^t$. 

It is easy to see that $\{x_1,e=(1,1,...,1)\}$ are the basis of the solution set to $Au^T=0, u\in \mathbb{F}_2^t$ for $r_4=1$. Obviously, the system of linear equations $M^{\prime}{v^{\prime}}^T=0, v^{\prime}\in \mathbb{F}_2^{2t}$ have $3$ linearly independent solutions $(e,O),(O,e), $ $(O,x_1)$ $\in \mathbb{F}_2^{2t}$ and $s_2\leq2+2r_4(-n)=4$. This state $s_2=4$ by $M^{\prime}$ is skew-symmetric matrix. Then we denote the last linearly independent solution be $(y_2,x_2)$. We assume that $r_1y_2^T=0$ and $r_{-3}y_2^T+r_{-1}x_2^T=0$, if not, we can replace $y_2$ with $y_2+e$ and replace $x_2$ with $x_2+e$.

The basis of the solution set to $Mv^T=0$ are $\{(1,0,0,1,0,1,e,O),$ $(0,0,0,,$ $0,0,O,e),$ $(0,0,0,0,0,0,O,x_1),$ $(0,0,0,r_{-1}y_2^T,0,r_3y_2^T,y_2,x_2)\}$ which correspond to $\Pi_0=(-3,-n)$, $\Pi_1=(n,1)$, $\Lambda_0=(d,1)$ and $\Lambda_1=(b_1^{\prime},b_2^{\prime})$ $=(3^{r_3y_2^T}\mathcal{p}^{x_2},\mathcal{p}^{y_2})$ in ${\rm Sel}_2(E_n)$, respectively. It is worth noting that $\Pi_i,$ $i=1,2$ generate the 2-torsion points of $E_n$.

Now we calculate the Cassels pairing $\left \langle \Lambda_0,\Lambda_1 \right \rangle$. 

Note that $r_4=1$, $dx^2+n/dy^2=4z^2$ have a primitive integer solution denoted by $(a,b,c)$. We  can assume $a\equiv b\equiv 1(\text{mod }4)$ for $n\equiv 3\pmod{8}$. For $C_{\Lambda_0}$, we take $Q_1=(t,u_2,u_3)=(b,-ad,4c)$, $Q_2=(t,u_1,u_3)=(0,1,-1)$ and $Q_3=(t,u_1,u_2)=(b,-2c,-ad)$ such that $Q_i$ is a solution of $H_i$ for $i=1,2,3$. Then $L_i$ linear form  in three variables such that $L_i=0$ defined the tangent plane of $H_i$ at $Q_i$ is 

$$ \left\{
\begin{aligned}
	L_1:&\ 2\frac{n}{d}bt-au_2-2cu_3=0, \\
	L_2:&\ u_1-u_3=0, \\
	L_3:&\ b\frac{n}{d}t-au_2+2cu_1=0.
\end{aligned}
\right.
$$

Let $b_3^{\prime}=b_1^{\prime}b_2^{\prime}$, then the Cassels pairing $\left \langle \Lambda_0,\Lambda_1 \right \rangle=\sum\limits_{p|24n\infty} \left \langle \Lambda_0,\Lambda_1 \right \rangle_p=\sum\limits_{p|24n\infty} \sum\limits_{i=1}^3 \left[ L_i(P_p),b_i^{\prime} \right]_p$ where $P_p=(t,u_1,u_2,u_3)\in C_{\Lambda_0}(\mathbb{Q}_p)$ and $\left[ \cdot,\cdot \right]_p$ is additional Hilbert symbol by Lemma 7.2 and 7.4 of \cite{C98}. 

$1)$ For $p=\infty$.  $\left[ L_i(P_\infty),b_i^{\prime} \right]_\infty=0$ for $i=1,2,3$ by $b_i^{\prime}>0$.

$2)$ For $p\neq2,3 \mid d$. We take $P_p=(1,\alpha_p,0,-2\alpha_p)$ where $\alpha_p^2\equiv4(\text{mod }p)$ and  $p\mid 2c-\alpha_p b$, i.e., $2c+\alpha_p b\equiv4c(\text{mod }p)$.

\begin{align*}
	\left \langle \Lambda_0,\Lambda_1 \right \rangle_p=&\left[2\frac{n}{d}b+4c\alpha_p,b_1^{\prime} \right]_p+\left[3\alpha_p,b_2^{\prime} \right]_p+\left[\frac{n}{d}b+2c\alpha_p,b_1^{\prime}b_2^{\prime} \right]_p\\
	=&\left[2,b_1^{\prime} \right]_p+\left[3(2c+\alpha_p b),b_2^{\prime} \right]_p \\
	=&\left[2,b_1^{\prime} \right]_p+\left[3c,b_2^{\prime} \right]_p.
\end{align*}

$3)$ For $p\neq2,3 \mid \frac{n}{d}$. We take $P_p=(0,1,\beta_p,-1)$ where $\beta_p^2\equiv4(\text{mod }p)$ and  $p\mid 2c+\beta_p a$, i.e., $2c-\beta_p a\equiv4c(\text{mod }p)$.

\begin{align*}
	\left \langle \Lambda_0,\Lambda_1 \right \rangle_p=&\left[-a\beta_p+2c,b_1^{\prime} \right]_p+\left[2,b_2^{\prime} \right]_p+\left[-a\beta_p+2c,b_1^{\prime}b_2^{\prime} \right]_p\\
	=&\left[2(2c-\beta_p a),b_2^{\prime} \right]_p \\
	=&\left[2c,b_2^{\prime} \right]_p.
\end{align*}

$4)$ For $p=3$. When $d\equiv1(\text{mod }3)$, we take $P_3=(0,2,1,1)$ if $3\mid a$ and $P_p=(1,2,0,1)$ if $3\mid b$. When $d\equiv2(\text{mod }3)$, we take $P_p=(t,2,u_2,1)$ with $t,u_2\nmid 3$ and such that $3\mid bt+au_2$. Then we have $\left \langle \Lambda_0,\Lambda_1 \right \rangle_3=[c,b_2^{\prime}]_3=0.$
 
$5)$ For $p=2$. 

$i)$ When $d\equiv1(\text{mod }8)$, we take $P_2=(0,1,1,-1)$, then $\left \langle \Lambda_0,\Lambda_1 \right \rangle_2=[2,b_2^{\prime}]_2=\left[\frac{2}{b_2^{\prime}}\right].$

$ii)$ When $d\equiv3(\text{mod }8)$, we take $P_2=(1,1,0,2)$, then $\left \langle \Lambda_0,\Lambda_1 \right \rangle_2=[2,b_2^{\prime}]_2=\left[\frac{2}{b_1^{\prime}}\right].$ 

$iii)$ When $d\equiv5(\text{mod }8)$, we take $P_2=(1,1,0,2)$, then $\left \langle \Lambda_0,\Lambda_1 \right \rangle_2=[2,b_2^{\prime}]_2=\left[\frac{2}{b_1^{\prime}b_2^{\prime}}\right].$ 

$iv)$ When $d\equiv7(\text{mod }8)$, we take $P_2=(1,1,0,2)$, then $\left \langle \Lambda_0,\Lambda_1 \right \rangle_2=[2,b_2^{\prime}]_2=0.$

Then we have $\left \langle \Lambda_0,\Lambda_1 \right \rangle_2=(1-\left[\frac{-1}{d}\right])\left[\frac{2}{b_2^{\prime}}\right]+\left[\frac{2}{d}\right]\left[\frac{2}{b_1^{\prime}}\right].$

In summary, we get 
$$\left \langle \Lambda_0,\Lambda_1 \right \rangle=\sum\limits_{p|24n\infty} \left \langle \Lambda_0,\Lambda_1 \right \rangle_p=x_1(D_6+r_{-1}^Tr_2+r_2^Tr_3)y_2^T+x_1(D_2+r_2^Tr_2)x_2+r_cy_2^T.$$ It is worth noting that $(x_1,O)M^{\prime}(y_2,x_2)^T=0$ and $\left \langle \Lambda_0,\Pi_1 \right \rangle=r_{-3}x_1^T=0$ which is following from the definition of Cassels pairing, therefore we have
$$\left \langle \Lambda_0,\Lambda_1 \right \rangle=x_1(D_{-2}+r_{-1}^Tr_{-1}+r_{2}^T r_{2})(y_2+x_2)^T+r_6 x_1^T r_6y_2^T+r_cy_2^T(\text{mod }2).$$

We know for skew-symmetric matrix $N$, $v$ lies in ${\rm Im}N$ if and only if $v^Tu=0$ for each $u\in {\rm Ker}N$. Thus the fact $(e,O)v_c^T=(O,e)v_c^T=(O,x_1)v_c^T=0$ tells us $\left \langle \Lambda_0,\Lambda_1 \right \rangle=v_c(y_2,x_2)^T=1$ if and only if $M^{\prime} {v^{\prime}}^T=v_c^T$ has no solution. $\Box$

\begin{corollary}
	Let $n=\widetilde{n}=p_1...p_t\equiv19\pmod{24},p_i\neq2,3$ and $r_4(-n)=1$, i.e., ${\rm Cl}^+(K)[2]\cap2{\rm Cl}^+(K)=\{1,[(d,\frac{-n+\sqrt{-n}}{2})]\}$ for some $d|n$. If 
	
	$1)$ $n$ have only one prime factor $p_t\equiv3\pmod{4}$ and $r_8=1-\left[\frac{6}{d}\right]$; or
	
	$2)$ all prime factors of $d$ congruent to $1\pmod{8}$and $r_8=0$,
	
	\noindent then $n$ is a non-$\pi/3$-congruent.
\end{corollary}

\emph{\bf Proof:} 
First, we show that $r_cy_2^T=1$ if and only if $r_8=0$. For $r_4=1$, according to Theorem 4.1 of Kolster\cite{K05}, $r_8=0$ if and only if $rk(\begin{array}{c|c}
	A & r_c^T
\end{array})=rkA+1$. In other words, $r_8=0$ if and only if linear equations $Au=r_c^T$ have no solutions which means $M^{\prime}{v^{\prime}}^T=(r_c,O)^T$ have no solutions. This states that $r_8=0$ if and only if $r_cy_2^T=1$ or $r_c e^T=1$. Then note that $r_c e^T=\left[\frac{c}{n}\right]=\left[\frac{-n}{c}\right]=0$ by the definition of $c$, we have $r_cy_2^T=1$ if and only if $r_8=0$.

For condition $1)$, we have $x_1=y_2$, $x_2=0$, then $\left \langle \Lambda_0,\Lambda_1 \right \rangle=r_cy_2^T+r_6x_1^T=1$.

For condition $2)$, we have $\left \langle \Lambda_0,\Lambda_1 \right \rangle=r_cy_2^T=1$.

Then $\Lambda_0,\Lambda_1 \in \Sha(E_n)$ and $rank_\mathbb{Z} E_n (\mathbb{Q})=0$. $\Box$

\subsection{\bf $\eta=2$} 
Let 
$$B_1=\begin{pmatrix}
	1 & 0 & 0 & 1 & 0 & 0 & O & O\\
	0 & 0 & 0 & 1 & [\frac{-1}{\widetilde{n}}] & 1 & O & r_{-1}\\
	0 & 1 & 0 & 0 & [\frac{2}{\widetilde{n}}] & 1 & O & r_{2}\\
	1 & 1 & 0 & 0 & 0 & [\frac{-3}{\widetilde{n}}] & r_{-3} & O\\
	0 & 0 & 1 & 0 & 0 & 0 & O & O\\
	O & O & O & r_{-1}^T & r_{2}^T & r_{3}^T & D_{-3} & A+D_{2}\\
	r_{-1}^T & r_{2}^T & r_{3}^T & O & O & O & A+D_{-2} & O
\end{pmatrix}.$$

\begin{theorem}
	For $n=2\widetilde{n}=2p_1...p_t,p_i \neq 2,3$, then we have Monsky matrix $M_n=B_1$.
\end{theorem}

\emph{\bf Proof:} We list all the equivalent conditions and transform them to additional Legendre symbols to build system of linear equations. 

We know that $C_{\Lambda}(\mathbb{A})\neq \emptyset$ if and only if $C_{\Lambda}(\mathbb{Q}_S)\neq \emptyset$, i.e., 

1) $C_{\Lambda}(\mathbb{R})\neq \emptyset$ if and only if 		$b_1b_2>0$ that is equivalent to $\gamma_1+\xi_1=0$;

2) $C_{\Lambda}(\mathbb{Q}_{p_i})\neq \emptyset$ for $i=1,..,t$ if and only if 
\begin{equation*}
	\begin{cases}
		\left( \frac{\dot{b}_1}{p_i}\right)=\left( \frac{\dot{b}_2}{p_i}\right)=1, &\text{if}\ p_i \nmid b_1b_2\\
		\left( \frac{\dot{n}\dot{b}_1}{p_i}\right)=\left( \frac{\dot{b}_2}{p_i}\right)=1, &\text{if}\ p_i\ |\ b_1,p_i\nmid b_2 \\
		\left( \frac{-3\dot{b}_1}{p_i}\right)=\left( \frac{-\dot{n}\dot{b}_2}{p_i}\right)=1, &\text{if}\ p_i\nmid b_1,p_i\ |\  b_2\\
		\left( \frac{-3\dot{n}\dot{b}_1}{p_i}\right)=\left( \frac{-\dot{n}\dot{b}_2}{p_i}\right)=1, &\text{if}\ p_i\ |\ b_1,p_i\ |\  b_2 
	\end{cases}
\end{equation*}

by Hensel's Lemma. The conditions are equivalent to
\begin{equation*}
	\begin{cases}
		[\frac{\dot{b}_1}{p_i}]+[\frac{-3}{p_i}]y_i+[\frac{2\widetilde{n}p_i}{p_i}]x_i=0 \\
		[\frac{\dot{b}_2}{p_i}]+[\frac{-2\widetilde{n}p_i}{p_i}]y_i=0 
	\end{cases};
\end{equation*}
In other words, we have

\begin{equation*}
	\begin{cases}
		[\frac{-1}{p_i}]\gamma_1+[\frac{2}{p_i}]\gamma_2+[\frac{3}{p_i}]\gamma_3+[\frac{-3}{p_i}]y_i+\sum\limits_{j\neq i}[\frac{p_j}{p_i}]x_j+[\frac{2\widetilde{n}p_i}{p_i}]x_i=0 \\
		[\frac{-1}{p_i}]\xi_1+[\frac{2}{p_i}]\xi_2+[\frac{3}{p_i}]\xi_3+\sum\limits_{j\neq i}[\frac{p_j}{p_i}]y_j+[\frac{-2\widetilde{n}p_i}{p_i}]y_i=0 
	\end{cases};
\end{equation*}

3) $C_{\Lambda}(\mathbb{Q}_3)\neq \emptyset$ if and only if $3\nmid b_2$ and
\begin{equation*}
	\begin{cases}
		\left( \frac{\dot{b}_2}{3}\right)=1, &\text{if}\ 3 \nmid b_1\\
		\left( \frac{-\dot{n}\dot{b}_2}{3}\right)=1, &\text{if}\ 3\ |\ b_1 
	\end{cases}
\end{equation*}

by Hessel's Lemma. The conditions are equivalent to
\begin{equation*}
	\begin{cases}
		[\frac{-3}{\widetilde{n}}]\gamma_3+[\frac{-3}{\dot{b}_2}]=0\\
		\xi_3=0
	\end{cases}.
\end{equation*}
In other words, we have 
\begin{equation*}
	\begin{cases}
		\xi_1+\xi_2+[\frac{-3}{\widetilde{n}}]\gamma_3+\sum\limits_{i}[\frac{-3}{p_i}]y_i=0\\
		\xi_3=0
	\end{cases};
\end{equation*}

4) $C_{\Lambda}(\mathbb{Q}_2)\neq \emptyset$. If  $\widetilde{n}\equiv 1({\rm mod }\ 8)$, $C_{\Lambda}(\mathbb{Q}_2)\neq \emptyset$ if and only if 

\begin{equation*}
	\begin{cases}
		(\dot{b}_1,\dot{b}_2)\equiv(1,1),(1,7)(\text{mod }8,\text{mod }8), &\text{if}\ 2 \nmid b_1b_2\\
		\dot{b}_1\equiv1(\text{mod }8), &\text{if}\ 2\ |\ b_1,2\nmid b_2 \\
		(\dot{b}_1,\dot{b}_2)\equiv(5,1),(5,7)(\text{mod }8,\text{mod }8), &\text{if}\ 2\nmid b_1,2\ |\  b_2\\
		\dot{b}_1\equiv5(\text{mod }8), &\text{if}\ 2\ |\ b_1,2\ |\  b_2 
	\end{cases}.
\end{equation*}

If  $\widetilde{n}\equiv 3({\rm mod }\ 8)$, $C_{\Lambda}(\mathbb{Q}_2)\neq \emptyset$ if and only if 

\begin{equation*}
	\begin{cases}
		(\dot{b}_1,\dot{b}_2)\equiv(1,1),(1,3)(\text{mod }8,\text{mod }8), &\text{if}\ 2 \nmid b_1b_2\\
		\dot{b}_1\equiv3(\text{mod }8), &\text{if}\ 2\ |\ b_1,2\nmid b_2 \\
		(\dot{b}_1,\dot{b}_2)\equiv(5,5),(5,7)(\text{mod }8,\text{mod }8), &\text{if}\ 2\nmid b_1,2\ |\  b_2\\
		\dot{b}_1\equiv7(\text{mod }8), &\text{if}\ 2\ |\ b_1,2\ |\  b_2 
	\end{cases}.
\end{equation*}

If  $\widetilde{n}\equiv 5({\rm mod }\ 8)$, $C_{\Lambda}(\mathbb{Q}_2)\neq \emptyset$ if and only if 

\begin{equation*}
	\begin{cases}
		(\dot{b}_1,\dot{b}_2)\equiv(1,1),(1,7)(\text{mod }8,\text{mod }8), &\text{if}\ 2 \nmid b_1b_2\\
		\dot{b}_1\equiv5(\text{mod }8), &\text{if}\ 2\ |\ b_1,2\nmid b_2 \\
		(\dot{b}_1,\dot{b}_2)\equiv(5,3),(5,5)(\text{mod }8,\text{mod }8), &\text{if}\ 2\nmid b_1,2\ |\  b_2\\
		\dot{b}_1\equiv1(\text{mod }8), &\text{if}\ 2\ |\ b_1,2\ |\  b_2 
	\end{cases}.
\end{equation*}

If  $\widetilde{n}\equiv 7({\rm mod }\ 8)$, $C_{\Lambda}(\mathbb{Q}_2)\neq \emptyset$ if and only if 

\begin{equation*}
	\begin{cases}
		(\dot{b}_1,\dot{b}_2)\equiv(1,1),(1,3)(\text{mod }8,\text{mod }8), &\text{if}\ 2 \nmid b_1b_2\\
		\dot{b}_1\equiv7(\text{mod }8), &\text{if}\ 2\ |\ b_1,2\nmid b_2 \\
		(\dot{b}_1,\dot{b}_2)\equiv(5,1),(5,3)(\text{mod }8,\text{mod }8), &\text{if}\ 2\nmid b_1,2\ |\  b_2\\
		\dot{b}_1\equiv3(\text{mod }8), &\text{if}\ 2\ |\ b_1,2\ |\  b_2 
	\end{cases}.
\end{equation*}
 The conditions are equivalent to

\begin{equation*}
	\begin{cases}
		[\frac{-1}{\widetilde{n}}]\gamma_2+[\frac{-1}{\dot{b}_1}]=0,\\
		[\frac{2}{\widetilde{n}}]\gamma_2+\xi_2+[\frac{2}{\dot{b}_1}]=0,\\
		[\frac{2}{\widetilde{n}}]\xi_2+[\frac{2}{\dot{b}_2}]+[\frac{-1}{\widetilde{n}}][\frac{-1}{\dot{b}_2}]=0,\ \ \text{if } \eta_2=0
	\end{cases}.
\end{equation*}

Note that if $C_{\Lambda}(\mathbb{Q}_3)$ and $C_{\Lambda}(\mathbb{Q}_{p_i})\neq \emptyset$, we have $[\frac{2}{\widetilde{n}}]\xi_2+[\frac{2}{\dot{b}_2}]+[\frac{-1}{\widetilde{n}}][\frac{-1}{\dot{b}_2}]=0$ by $\xi_3=0$ and $[\frac{\dot{b}_2}{p_i}]+[\frac{-2\widetilde{n}p_i}{p_i}]y_i=0$. In other words, if we assume $C_{\Lambda}(\mathbb{Q}_3)$ and $C_{\Lambda}(\mathbb{Q}_{p_i})\neq \emptyset$, then $C_{\Lambda}(\mathbb{Q}_2)\neq \emptyset$ if and only if

\begin{equation*}
	\begin{cases}
		\gamma_1+[\frac{-1}{\widetilde{n}}]\gamma_2+\gamma_3+\sum\limits_{i}[\frac{-1}{p_i}]x_i=0,\\
		[\frac{2}{\widetilde{n}}]\gamma_2+\gamma_3+\xi_2+\sum\limits_{i}[\frac{2}{p_i}]x_i=0
	\end{cases}.
\end{equation*}
The Monsky matrix can be naturally derived from the above system of linear
equations.
$\Box$

Then we have a corollary about pairty of 2-rank as following:

\begin{corollary}
	For $n=2\widetilde{n}=2p_1...p_t,p_i\neq2,3$, the rank of $2$-Selmer group $s_2(E_n)$ is even (resp. odd) if and only if $n\equiv2,14(resp.\ 10,22)\pmod{24}$, i.e., $\widetilde{n}\equiv1(resp.\ 2)\pmod{3}$. 
\end{corollary}

\emph{\bf Proof:}
Note that $s_2(E_n)=2t+6-r(M_n)$, then we only need consider
the rank of $M_n$ which is denoted by $r(M_n)$.
We have:

\begin{align*}
	r(B_1)&=r\begin{pmatrix}
		1 & 0 & 0 & 1 & 0 & 0 & O & O\\
		0 & 0 & 0 & 1 & [\frac{-1}{\widetilde{n}}] & 1 & O & r_{-1}\\
		0 & 1 & 0 & 0 & [\frac{2}{\widetilde{n}}] & 1 & O & r_{2}\\
		1 & 1 & 0 & 0 & 0 & [\frac{-3}{\widetilde{n}}] & r_{-3} & O\\
		0 & 0 & 1 & 0 & 0 & 0 & O & O\\
		O & O & O & r_{-1}^T & r_{2}^T & r_{3}^T & D_{-3} & A+D_{2}\\
		r_{-1}^T & r_{2}^T & r_{3}^T & O & O & O & A+D_{-2} & O
	\end{pmatrix} \\
	&=r\begin{pmatrix}
		1 & 0 & 0 & 1 & 0 & 0 & O & O\\
		0 & 0 & 0 & 1 & [\frac{-1}{\widetilde{n}}] & 1 & O & r_{-1}\\
		0 & 1 & 0 & 0 & [\frac{2}{\widetilde{n}}] & 1 & O & r_{2}\\
		0 & 0 & 0 & 0 & [\frac{2}{\widetilde{n}}] & [\frac{-3}{\widetilde{n}}] & r_{-3} & r_{-2}\\
		0 & 0 & 1 & 0 & 0 & 0 & O & O\\
		O & O & O & O & [\frac{-1}{\widetilde{n}}]r_{-1}^T+r_{-2}^T & r_{-3}^T & D_{-3} & A+D_{2}+r_{-1}^Tr_{-1}\\
		O & O & O & O & [\frac{-1}{\widetilde{n}}]r_{-1}^T+[\frac{2}{\widetilde{n}}]r_{-2}^T & r_{-2}^T & A+D_{-2} & r_{-1}^Tr_{-1}+r_{2}^Tr_{2}
	\end{pmatrix} \\ &\text{  by row transformations.} \\
	&=r\begin{pmatrix}
		1 & 0 & 0 & 0 & 0 & 0 & O & O\\
		0 & 0 & 0 & 1 & 0 & 0 & O & O\\
		0 & 1 & 0 & 0 & 0 & 0 & O & O\\
		0 & 0 & 0 & 0 & 0 & [\frac{-3}{\widetilde{n}}] & r_{-3} & r_{-2}\\
		0 & 0 & 1 & 0 & 0 & 0 & O & O\\
		O & O & O & O & O & r_{-3}^T & D_{-3} & A+D_{2}+r_{-1}^Tr_{-1}\\
		O & O & O & O & O & r_{-2}^T & A+D_{-2} & r_{-1}^Tr_{-1}+r_{2}^Tr_{2}
	\end{pmatrix} \\ &\text{  by column transformations.} \\
	&=4+r\begin{pmatrix}
		[\frac{-3}{\widetilde{n}}]  & r_{-3} & r_{-2}\\
		r_{-3}^T & D_{-3} & A^T+D_{-2}\\
		r_{-2}^T & A+D_{-2} & r_{-1}^Tr_{-1}+r_{2}^Tr_{2}
	\end{pmatrix} \\
	&\text{  by } A+A^T+r_{-1}^Tr_{-1}+D_{-1}=O. 
\end{align*}

If $[\frac{-3}{\widetilde{n}}]=0$, i.e. $n\equiv2,14(\text{mod }24)$, $$r(B_1)=4+r\begin{pmatrix}
	0  & r_{-3} & r_{-2}\\
	r_{-3}^T & D_{-3}+r_{-3}^Tr_{-3} & A^T+D_{-2}+r_{-3}^Tr_{-2}\\
	r_{-2}^T & A+D_{-2}+r_{-2}^Tr_{-3} & r_{-1}^Tr_{-1}+r_{2}^Tr_{2}+r_{-2}^Tr_{-2}
\end{pmatrix} $$ is even; and if $[\frac{-3}{\widetilde{n}}]=1$, i.e. $n\equiv10,22(\text{mod }24)$, $$r(B_1)=5+r\begin{pmatrix}
 D_{-3}+r_{-3}^Tr_{-3} & A^T+D_{-2}+r_{-3}^Tr_{-2}\\
 A+D_{-2}+r_{-2}^Tr_{-3} & r_{-1}^Tr_{-1}+r_{2}^Tr_{2}+r_{-2}^Tr_{-2}
\end{pmatrix} $$ is odd. $\Box$

\subsection{\bf $\eta=3$} Let 
$$C_1=\begin{pmatrix}
	1 & 0 & 0 & 1 & 0 & 0 & O & O\\
	0 & 1 & 0 & 1 & 0 & 1 & O & r_{-1}\\
	0 & 0 & 0 & 0 & 1 & 0 & O & O\\
	1 & 1 & 1 & 0 & 0 & 1 & r_{-1} & r_2\\
	1 & 1 & 1+[\frac{-3}{\widetilde{n}}] & 0 & 0 & 0 & r_{-3} & O\\
	0 & 0 & 1+[\frac{-3}{\widetilde{n}}] & 1 & 1 & [\frac{-3}{\widetilde{n}}] & O & r_{-3}\\
	O & O & O & r_{-1}^T & r_{2}^T & r_{3}^T & D_{-3} & A+D_{3}\\
	r_{-1}^T & r_{2}^T & r_{3}^T & O & O & O & A+D_{-3} & O
\end{pmatrix},$$

and 
$$C_2=\begin{pmatrix}
	1 & 0 & 0 & 1 & 0 & 0 & O & O\\
	0 & 0 & 0 & 0 & 1 & 0 & O & O\\
	0 & 1 & 0 & 0 & 0 & 0 & O & O\\
	0 & 0 & 0 & 1 & 0 & 1 & O & r_{-1}\\
	1 & 1 & 1+[\frac{-3}{\widetilde{n}}] & 0 & 0 & 0 & r_{-3} & O\\
0 & 0 & 1+[\frac{-3}{\widetilde{n}}] & 1 & 1 & [\frac{-3}{\widetilde{n}}] & O & r_{-3}\\
O & O & O & r_{-1}^T & r_{2}^T & r_{3}^T & D_{-3} & A+D_{3}\\
r_{-1}^T & r_{2}^T & r_{3}^T & O & O & O & A+D_{-3} & O
\end{pmatrix},$$

and
$$C_3=\begin{pmatrix}
	1 & 0 & 0 & 1 & 0 & 0 & O & O\\
	0 & 1 & 0 & 0 & 0 & 0 & O & O\\
	0 & 0 & 0 & 0 & 1 & 0 & O & O\\
	1 & 0 & 1 & 0 & 0 & 0 & r_{-1} & O\\
	1 & 1 & 1+[\frac{-3}{\widetilde{n}}] & 0 & 0 & 0 & r_{-3} & O\\
O & O & 1+[\frac{-3}{\widetilde{n}}] & 1 & 1 & [\frac{-3}{\widetilde{n}}] & O & r_{-3}\\
O & O & O & r_{-1}^T & r_{2}^T & r_{3}^T & D_{-3} & A+D_{3}\\
r_{-1}^T & r_{2}^T & r_{3}^T & O & O & O & A+D_{-3} & O
\end{pmatrix}.$$

\begin{theorem}
	For $n=3\widetilde{n}=3p_1...p_t,p_i \neq 2,3$, then we have Monsky matrix 
	\begin{equation*}
		M_n=\begin{cases}
			C_1, &\text{if } \widetilde{n}\equiv 3({\rm mod }\ 8),\\
			C_2, &\text{if } \widetilde{n}\equiv 7({\rm mod }\ 8), \\
			C_3, &\text{if } \widetilde{n}\equiv 1({\rm mod }\ 4).
		\end{cases}
	\end{equation*} 
\end{theorem}

\emph{\bf Proof:}We list all the equivalent conditions and transform them to additional Legendre symbols to build system of linear equations. 

We know that $C_{\Lambda}(\mathbb{A})\neq \emptyset$ if and only if $C_{\Lambda}(\mathbb{Q}_S)\neq \emptyset$, i.e.,

1) $C_{\Lambda}(\mathbb{R})\neq \emptyset$ if and only if 		$b_1b_2>0$ that is equivalent to $\gamma_1+\xi_1=0$;

2) $C_{\Lambda}(\mathbb{Q}_{p_i})\neq \emptyset$ for $i=1,..,t$ if and only if 
\begin{equation*}
	\begin{cases}
		\left( \frac{\dot{b}_1}{p_i}\right)=\left( \frac{\dot{b}_2}{p_i}\right)=1, &\text{if}\ p_i \nmid b_1b_2\\
		\left( \frac{\dot{n}\dot{b}_1}{p_i}\right)=\left( \frac{\dot{b}_2}{p_i}\right)=1, &\text{if}\ p_i\ |\ b_1,p_i\nmid b_2 \\
		\left( \frac{-3\dot{b}_1}{p_i}\right)=\left( \frac{-\dot{n}\dot{b}_2}{p_i}\right)=1, &\text{if}\ p_i\nmid b_1,p_i\ |\  b_2\\
		\left( \frac{-3\dot{n}\dot{b}_1}{p_i}\right)=\left( \frac{-\dot{n}\dot{b}_2}{p_i}\right)=1, &\text{if}\ p_i\ |\ b_1,p_i\ |\  b_2 
	\end{cases}
\end{equation*}

by Hensel's Lemma. The conditions are equivalent to
\begin{equation*}
	\begin{cases}
		[\frac{\dot{b}_1}{p_i}]+[\frac{-3}{p_i}]y_i+[\frac{3\widetilde{n}p_i}{p_i}]x_i=0 \\
		[\frac{\dot{b}_2}{p_i}]+[\frac{-3\widetilde{n}p_i}{p_i}]y_i=0 
	\end{cases};
\end{equation*}
In other words, we have
\begin{equation*}
	\begin{cases}
		[\frac{-1}{p_i}]\gamma_1+[\frac{2}{p_i}]\gamma_2+[\frac{3}{p_i}]\gamma_3+[\frac{-3}{p_i}]y_i+\sum\limits_{j\neq i}[\frac{p_j}{p_i}]x_j+[\frac{3\widetilde{n}p_i}{p_i}]x_i=0 \\
		[\frac{-1}{p_i}]\xi_1+[\frac{2}{p_i}]\xi_2+[\frac{3}{p_i}]\xi_3+\sum\limits_{j\neq i}[\frac{p_j}{p_i}]y_j+[\frac{-3\widetilde{n}p_i}{p_i}]y_i=0 
	\end{cases};
\end{equation*}

3) $C_{\Lambda}(\mathbb{Q}_3)\neq \emptyset$ if and only if
 \begin{equation*}
	\begin{cases}
		\left( \frac{\dot{b}_1}{3}\right)=\left( \frac{\dot{b}_2}{3}\right)=1, &\text{if}\ 3 \nmid b_1b_2\\
		\left( \frac{\dot{n}\dot{b}_1}{3}\right)=\left( \frac{\dot{b}_2}{3}\right)=1, &\text{if}\ 3\ |\ b_1,3\nmid b_2 \\
		\left( \frac{-\dot{n}\dot{b}_1}{3}\right)=\left( \frac{-\dot{n}\dot{b}_2}{3}\right)=1, &\text{if}\ 3\nmid b_1,3\ |\  b_2\\
		\left( \frac{-\dot{b}_1}{3}\right)=\left( \frac{-\dot{n}\dot{b}_2}{3}\right)=1, &\text{if}\ 3\ |\ b_1,3\ |\  b_2 
	\end{cases}
\end{equation*}

by Hensel's Lemma. The conditions are equivalent to

\begin{equation*}
	\begin{cases}		
		[\frac{-3}{\widetilde{n}}]\gamma_3+(1+[\frac{-3}{\widetilde{n}}])\xi_3+[\frac{-3}{\dot{b}_1}]=0\\
		(1+[\frac{-3}{\widetilde{n}}])\xi_3+[\frac{-3}{\dot{b}_2}]=0
	\end{cases};
\end{equation*}
In other words, we have 
\begin{equation*}
	\begin{cases}
		\gamma_1+\gamma_2+[\frac{-3}{\widetilde{n}}]\gamma_3+(1+[\frac{-3}{\widetilde{n}}])\xi_3+\sum\limits_{i}[\frac{-3}{p_i}]x_i=0\\
		\xi_1+\xi_2+(1+[\frac{-3}{\widetilde{n}}])\xi_3+\sum\limits_{i}[\frac{-3}{p_i}]y_i=0
	\end{cases};
\end{equation*}
4) $C_{\Lambda}(\mathbb{Q}_2)\neq \emptyset$. If  $\widetilde{n}\equiv 3({\rm mod }\ 8)$, $C_{\Lambda}(\mathbb{Q}_2)\neq \emptyset$ if and only if (i):$2\nmid b_1b_2$ and $(\dot{b}_1,\dot{b}_2)\equiv (1,1),(1,5),(5,7),(5,3)({\rm mod }\ 8,{\rm mod }\ 8)$ or (ii):$2\nmid b_1,2\ |\ b_2$ and $(\dot{b}_1,\dot{b}_2)\equiv (7,3),$ $(3,1)({\rm mod }\ 8,{\rm mod }\ 4)$;
if $\widetilde{n}\equiv 7({\rm mod }\ 8)$, $C_{\Lambda}(\mathbb{Q}_2)\neq \emptyset$ if and only if $2\nmid b_1b_2$ and $\dot{b}_1\equiv 1({\rm mod }\ 4)$; 
if $\widetilde{n}\equiv 1({\rm mod }\ 4)$, $C_{\Lambda}(\mathbb{Q}_2)\neq \emptyset$ if and only if $2\nmid b_1b_2$ and $\dot{b}_2\equiv 1({\rm mod }\ 4)$ by Hensel's Lemma. The conditions are equivalent to

\begin{equation*}
	\begin{cases}
		\xi_2+\gamma_1+\gamma_3+\sum\limits_{i}[\frac{-1}{p_i}]x_i=0, \ \ \ \ \gamma_2=0,\\
		\xi_1+\xi_2+\xi_3+\gamma_3+\sum\limits_{i}[\frac{-1}{p_i}]y_i+\sum\limits_{i}[\frac{2}{p_i}]x_i=0
	\end{cases} \text{if } \widetilde{n}\equiv 3({\rm mod }\ 8),
\end{equation*}
and
\begin{equation*}
	\begin{cases}
		\xi_2=0, \ \ \ \ \gamma_2=0,\\
		\gamma_1+\gamma_3+\sum\limits_{i}[\frac{-1}{p_i}]x_i=0
	\end{cases} \text{if } \widetilde{n}\equiv 7({\rm mod }\ 8),
\end{equation*}
and
\begin{equation*}
	\begin{cases}
		\xi_2=0, \ \ \ \ \gamma_2=0,\\
		\xi_1+\xi_3+\sum\limits_{i}[\frac{-1}{p_i}]y_i=0
	\end{cases} \text{if } \widetilde{n}\equiv 1({\rm mod }\ 4).
\end{equation*}
The Monsky matrix can be naturally derived from the above system of linear
equations.
$\Box$

Then we have a corollary about pairty of 2-rank as following:

\begin{corollary}
	For $n=3\widetilde{n}=3p_1...p_t,p_i\neq2,3$, the rank of $2$-Selmer group $s_2(E_n)$ is even (resp. odd) if and only if $n\equiv3,9,15(resp.\ 21)\pmod{24}$,i.e. $\widetilde{n}\equiv1,3,5(resp.\ 7)\pmod{8}$. 
\end{corollary}

\emph{\bf Proof:}
Note that $s_2(E_n)=2t+6-r(M_n)$, then we only need consider
the rank of $M_n$ which is denoted by $r(M_n)$.

For $\widetilde{n}\equiv3(\text{mod 8})$, we have:
\begin{align*}
	r(C_1)&=r\begin{pmatrix}
		1 & 0 & 0 & 1 & 0 & 0 & O & O\\
		0 & 1 & 0 & 1 & 0 & 1 & O & r_{-1}\\
		0 & 0 & 0 & 0 & 1 & 0 & O & O\\
		1 & 1 & 1 & 0 & 0 & 1 & r_{-1} & r_2\\
		1 & 1 & 1+[\frac{-3}{\widetilde{n}}] & 0 & 0 & 0 & r_{-3} & O\\
		0 & 0 & 1+[\frac{-3}{\widetilde{n}}] & 1 & 1 & [\frac{-3}{\widetilde{n}}] & O & r_{-3}\\
		O & O & O & r_{-1}^T & r_{2}^T & r_{3}^T & D_{-3} & A+D_{3}\\
		r_{-1}^T & r_{2}^T & r_{3}^T & O & O & O & A+D_{-3} & O
	\end{pmatrix} \\
	&=r\begin{pmatrix}
		1 & 0 & 0 & 1 & 0 & 0 & O & O\\
		0 & 1 & 0 & 1 & 0 & 1 & O & r_{-1}\\
		0 & 0 & 0 & 0 & 1 & 0 & O & O\\
		0 & 0 & 1 & 0 & 0 & 0 & r_{-1} & r_{-2}\\
		0 & 0 & 1+[\frac{-3}{\widetilde{n}}] & 0 & 0 & 1 & r_{-3} & r_{-1}\\
		0 & 0 & 0 & 0 & 0 & 0 & O & O\\
		O & r_{-1}^T & O & O & O & r_{-3}^T & D_{-3} & A+D_{3}+r_{-1}^Tr_{-1}\\
		O & r_{-2}^T & r_{3}^T & O & O & r_{-1}^T & A+D_{-3} & r_{-1}^Tr_{-1}
	\end{pmatrix} \\ &\text{  by row transformations.} \\
	&=r\begin{pmatrix}
		1 & 0 & 0 & 0 & 0 & 0 & O & O\\
		0 & 0 & 0 & 1 & 0 & 0 & O & O\\
		0 & 0 & 0 & 0 & 1 & 0 & O & O\\
		0 & 0 & 0 & 0 & 0 & 0 & r_{-1} & r_{-2}\\
		0 & 0 & 0 & 0 & 0 & 1 & r_{-3} & r_{-1}\\
		0 & 0 & 0 & 0 & 0 & 0 & O & O\\
		O & r_{-1}^T & O & O & O & r_{-3}^T & D_{-3} & A+D_{3}+r_{-1}^Tr_{-1}\\
		O & r_{-2}^T & O & O & O & r_{-1}^T & A+D_{-3} & r_{-1}^Tr_{-1}
	\end{pmatrix} \\ &\text{  by column transformations.} \\
	&=3+r\begin{pmatrix}
		1 & 0  & r_{-3} & r_{-1}\\
		0 & 0  & r_{-1} & r_{-2}\\
		r_{-3}^T & r_{-1}^T & D_{-3} & A^T+D_{-3}\\
		r_{-1}^T &r_{-2}^T & A+D_{-3} & r_{-1}^Tr_{-1}
	\end{pmatrix} \\
	&\text{  by } A+A^T+r_{-1}^Tr_{-1}+D_{-1}=O. \\
	&=4+r\begin{pmatrix}
	 0  & r_{-1} & r_{-2}\\
	 r_{-1}^T & D_{-3}+r_{-3}^Tr_{-3} & A^T+D_{-3}+r_{-3}^Tr_{-1}\\
	r_{-2}^T & A+D_{-3}+r_{-1}^Tr_{-3} & O
	\end{pmatrix}
\end{align*}

Then we have $r(C_1)$ is even for $n\equiv9(\text{mod 24})$.

For $\widetilde{n}\equiv7(\text{mod 8})$, we have:
\begin{align*}
	r(C_2)&=r\begin{pmatrix}
		1 & 0 & 0 & 1 & 0 & 0 & O & O\\
		0 & 0 & 0 & 0 & 1 & 0 & O & O\\
		0 & 1 & 0 & 0 & 0 & 0 & O & O\\
		0 & 0 & 0 & 1 & 0 & 1 & O & r_{-1}\\
		1 & 1 & 1+[\frac{-3}{\widetilde{n}}] & 0 & 0 & 0 & r_{-3} & O\\
		0 & 0 & 1+[\frac{-3}{\widetilde{n}}] & 1 & 1 & [\frac{-3}{\widetilde{n}}] & O & r_{-3}\\
		O & O & O & r_{-1}^T & r_{2}^T & r_{3}^T & D_{-3} & A+D_{3}\\
		r_{-1}^T & r_{2}^T & r_{3}^T & O & O & O & A+D_{-3} & O
	\end{pmatrix} \\
	&=r\begin{pmatrix}
		1 & 0 & 0 & 1 & 0 & 0 & O & O\\
		0 & 0 & 0 & 0 & 1 & 0 & O & O\\
		0 & 1 & 0 & 0 & 0 & 0 & O & O\\
		0 & 0 & 0 & 1 & 0 & 1 & O & r_{-1}\\
		0 & 0 & 0 & 0 & 0 & 0 & O & O\\
		0 & 0 & 1+[\frac{-3}{\widetilde{n}}] & 0 & 0 & 1 & r_{-3} & r_{-1}\\
		O & O & O & O & O & r_{-3}^T & D_{-3} & A+D_{3}+r_{-1}^Tr_{-1}\\
		O & O & r_{3}^T & O & O & r_{-1}^T & A+D_{-3} & r_{-1}^Tr_{-1}
	\end{pmatrix} \\ &\text{  by row transformations.} \\
	&=r\begin{pmatrix}
		1 & 0 & 0 & 0 & 0 & 0 & O & O\\
		0 & 0 & 0 & 0 & 1 & 0 & O & O\\
		0 & 1 & 0 & 0 & 0 & 0 & O & O\\
		0 & 0 & 0 & 1 & 0 & 0 & O & O\\
		0 & 0 & 0 & 0 & 0 & 0 & O & O\\
		0 & 0 & 0 & 0 & 0 & 1 & r_{-3} & r_{-1}\\
		O & O & O & O & O & r_{-3}^T & D_{-3} & A+D_{3}+r_{-1}^Tr_{-1}\\
		O & O & O & O & O & r_{-1}^T & A+D_{-3} & r_{-1}^Tr_{-1}
	\end{pmatrix} \\ &\text{  by column transformations.} \\
	&=4+r\begin{pmatrix}
		1  & r_{-3} & r_{-1}\\
		r_{-3}^T & D_{-3} & A^T+D_{-3}\\
		r_{-1}^T & A+D_{-3} & r_{-1}^Tr_{-1}
	\end{pmatrix} \\
	&\text{  by } A+A^T+r_{-1}^Tr_{-1}+D_{-1}=O. \\
		&=5+r\begin{pmatrix}
 D_{-3}+r_{-3}^Tr_{-3} & A^T+D_{-3}+r_{-1}^Tr_{-3}\\
 A+D_{-3}+r_{-1}^Tr_{-3} & O
	\end{pmatrix}
\end{align*}

Then we have $r(C_2)$ is odd for $n\equiv21(\text{mod 24})$. 

For $\widetilde{n}\equiv1(\text{mod 4})$, we have:
\begin{align*}
	r(C_3)&=r\begin{pmatrix}
		1 & 0 & 0 & 1 & 0 & 0 & O & O\\
		0 & 1 & 0 & 0 & 0 & 0 & O & O\\
		0 & 0 & 0 & 0 & 1 & 0 & O & O\\
		1 & 0 & 1 & 0 & 0 & 0 & r_{-1} & O\\
		1 & 1 & 1+[\frac{-3}{\widetilde{n}}] & 0 & 0 & 0 & r_{-3} & O\\
		O & O & 1+[\frac{-3}{\widetilde{n}}] & 1 & 1 & [\frac{-3}{\widetilde{n}}] & O & r_{-3}\\
		O & O & O & r_{-1}^T & r_{2}^T & r_{3}^T & D_{-3} & A+D_{3}\\
		r_{-1}^T & r_{2}^T & r_{3}^T & O & O & O & A+D_{-3} & O
	\end{pmatrix} \\
	&=r\begin{pmatrix}
		1 & 0 & 0 & 1 & 0 & 0 & O & O\\
		0 & 0 & 0 & 0 & 1 & 0 & O & O\\
		0 & 1 & 0 & 0 & 0 & 0 & O & O\\
		0 & 0 & 1 & 1 & 0 & 0 & r_{-1} & O\\
		0 & 0 & [\frac{-3}{\widetilde{n}}] & 0 & 0 & 0 & r_{3} & O\\
		0 & 0 & 0 & 0 & 0 & 0 & O & O\\
		O & O & r_{-1}^T & O & O & r_{3}^T & D_{-3}+r_{-1}^Tr_{-1} & A+D_3\\
		O & O & r_{-3}^T & O & O & O & A+D_{-3}+r_{-1}^Tr_{-1} & O
	\end{pmatrix} \\ &\text{  by row transformations.} \\
	&=r\begin{pmatrix}
		1 & 0 & 0 & 0 & 0 & 0 & O & O\\
		0 & 0 & 0 & 0 & 1 & 0 & O & O\\
		0 & 1 & 0 & 0 & 0 & 0 & O & O\\
		0 & 0 & 0 & 1 & 0 & 0 & O & O\\
		0 & 0 & 0 & 0 & 0 & 0 & r_{3} & O\\
		0 & 0 & 0 & 0 & 0 & 0 & O & O\\
		O & O & O & O & O & r_{3}^T & D_{-3}+r_{-1}^Tr_{-1} & A+D_3\\
		O & O & O & O & O & O & A+D_{-3}+r_{-1}^Tr_{-1} & O
	\end{pmatrix} \\ &\text{  by column transformations.} \\
	&=4+r\begin{pmatrix}
		0  & r_{3} & O\\
		r_{3}^T & D_{-3}+r_{-1}^Tr_{-1}+r_{3}^Tr_{3} & A+D_3\\
		O & A^T+D_3 & O
	\end{pmatrix} \\
	&\text{  by } A+A^T+r_{-1}^Tr_{-1}+D_{-1}=O. 
\end{align*}

 Then we have $r(C_3)$ is even for $n\equiv3,15(\text{mod 24})$.  $\Box$

Now we prove a result about non $\pi/3$-congruent numbers:

\begin{corollary}
	For $n=3\widetilde{n}=3p_1...p_t\equiv3,15\pmod{24},p_i\neq2,3$, and $r_4(-n)=0$, then $s_2(E_n)=2$. 
\end{corollary}

\emph{\bf Proof:} 
For $n\equiv3,15(\text{mod }24)$, we have

$$ s_2(E_n)=2t+6-r(M_n)=2t+2-r\begin{pmatrix}
	0  & r_{3} & O\\
	r_{3}^T & D_{-3}+r_{-1}^Tr_{-1}+r_{3}^Tr_{3} & A+D_{3}\\
	O & A^T+D_{3} & 
\end{pmatrix}.$$
Then $s_2\leq2t+2-2r(r_{3}^T\ A+D_3)=2t-2r(R(-n))=2+2r_4(-n)=2.$

Note that $s_2(E_n)\geq2$ by the fact that $E_n$ has $4$ elements of $2$-torsion, we have $s_2(E_n)=2$. $\Box$

\subsection{\bf $\eta=6$} 
Let 
$$D_1=\begin{pmatrix}
	1 & 0 & 0 & 1 & 0 & 0 & O & O\\
	0 & 0 & 0 & 1 & 1+[\frac{-1}{\widetilde{n}}] & 1 & O & r_{-1}\\
	0 & 1 & 0 & 0 & 1+[\frac{2}{\widetilde{n}}] & 1 & O & r_{2}\\
	1 & 1 & [\frac{-3}{\widetilde{n}}] & 0 & 0 & 0 & r_{-3} & O\\
0 & 0 & [\frac{-3}{\widetilde{n}}] & 1 & 1 & 1+[\frac{-3}{\widetilde{n}}] & O & r_{-3}\\
	0 & 0 & 0 & r_{-1}^T & r_{2}^T & r_{3}^T & D_{-3} & A+D_{6}\\
	r_{-1}^T & r_{2}^T & r_{3}^T & 0 & 0 & 0 & A+D_{-6} & O
\end{pmatrix}.$$

\begin{theorem}
	For $n=6\widetilde{n}=6p_1...p_t,p_i \neq 2,3$, then we have Monsky matrix $M_n=D_1$.
\end{theorem}

\emph{\bf Proof:} We list all the equivalent conditions and transform them to additional Legendre symbols to build system of linear equations. 

$C_{\Lambda}(\mathbb{A})\neq \emptyset$ if and only if $C_{\Lambda}(\mathbb{Q}_S)\neq \emptyset$ , i.e.,

1) $C_{\Lambda}(\mathbb{R})\neq \emptyset$ if and only if 		$b_1b_2>0$ that is equivalent to $\gamma_1+\xi_1=0$;

2) $C_{\Lambda}(\mathbb{Q}_{p_i})\neq \emptyset$ for $i=1,..,t$ if and only if 
\begin{equation*}
	\begin{cases}
		\left( \frac{\dot{b}_1}{p_i}\right)=\left( \frac{\dot{b}_2}{p_i}\right)=1, &\text{if}\ p_i \nmid b_1b_2\\
		\left( \frac{\dot{n}\dot{b}_1}{p_i}\right)=\left( \frac{\dot{b}_2}{p_i}\right)=1, &\text{if}\ p_i\ |\ b_1,p_i\nmid b_2 \\
		\left( \frac{-3\dot{b}_1}{p_i}\right)=\left( \frac{-\dot{n}\dot{b}_2}{p_i}\right)=1, &\text{if}\ p_i\nmid b_1,p_i\ |\  b_2\\
		\left( \frac{-3\dot{n}\dot{b}_1}{p_i}\right)=\left( \frac{-\dot{n}\dot{b}_2}{p_i}\right)=1, &\text{if}\ p_i\ |\ b_1,p_i\ |\  b_2 
	\end{cases}
\end{equation*}

by Hensel's Lemma. The conditions are equivalent to
\begin{equation*}
	\begin{cases}
		[\frac{\dot{b}_1}{p_i}]+[\frac{-3}{p_i}]y_i+[\frac{6\widetilde{n}p_i}{p_i}]x_i=0 \\
		[\frac{\dot{b}_2}{p_i}]+[\frac{-6\widetilde{n}p_i}{p_i}]y_i=0 
	\end{cases};
\end{equation*}
In other words, we have

\begin{equation*}
	\begin{cases}
		[\frac{-1}{p_i}]\gamma_1+[\frac{2}{p_i}]\gamma_2+[\frac{3}{p_i}]\gamma_3+[\frac{-3}{p_i}]y_i+\sum\limits_{j\neq i}[\frac{p_j}{p_i}]x_j+[\frac{6\widetilde{n}p_i}{p_i}]x_i=0 \\
		[\frac{-1}{p_i}]\xi_1+[\frac{2}{p_i}]\xi_2+[\frac{3}{p_i}]\xi_3+\sum\limits_{j\neq i}[\frac{p_j}{p_i}]y_j+[\frac{-6\widetilde{n}p_i}{p_i}]y_i=0 
	\end{cases};
\end{equation*}

3) $C_{\Lambda}(\mathbb{Q}_3)\neq \emptyset$ if and only if
\begin{equation*}
	\begin{cases}
		\left( \frac{\dot{b}_1}{3}\right)=\left( \frac{\dot{b}_2}{3}\right)=1, &\text{if}\ 3 \nmid b_1b_2\\
		\left( \frac{\dot{n}\dot{b}_1}{3}\right)=\left( \frac{\dot{b}_2}{3}\right)=1, &\text{if}\ 3\ |\ b_1,3\nmid b_2 \\
		\left( \frac{-\dot{n}\dot{b}_1}{3}\right)=\left( \frac{-\dot{n}\dot{b}_2}{3}\right)=1, &\text{if}\ 3\nmid b_1,3\ |\  b_2\\
		\left( \frac{-\dot{b}_1}{3}\right)=\left( \frac{-\dot{n}\dot{b}_2}{3}\right)=1, &\text{if}\ 3\ |\ b_1,3\ |\  b_2 
	\end{cases}
\end{equation*}

by Hensel's Lemma. The conditions are equivalent to

\begin{equation*}
	\begin{cases}		
		(1+[\frac{-3}{\widetilde{n}}])\gamma_3+[\frac{-3}{\widetilde{n}}]\xi_3+[\frac{-3}{\dot{b}_1}]=0\\
		[\frac{-3}{\widetilde{n}}]\xi_3+[\frac{-3}{\dot{b}_2}]=0
	\end{cases};
\end{equation*}
In other words, we have 
\begin{equation*}
	\begin{cases}
		\gamma_1+\gamma_2+(1+[\frac{-3}{\widetilde{n}}])\gamma_3+[\frac{-3}{\widetilde{n}}]\xi_3+\sum\limits_{i}[\frac{-3}{p_i}]x_i=0\\
		\xi_1+\xi_2+[\frac{-3}{\widetilde{n}}]\xi_3+\sum\limits_{i}[\frac{-3}{p_i}]y_i=0
	\end{cases};
\end{equation*}

4) $C_{\Lambda}(\mathbb{Q}_2)\neq \emptyset$.
If  $\widetilde{n}\equiv 1({\rm mod }\ 8)$, $C_{\Lambda}(\mathbb{Q}_2)\neq \emptyset$ if and only if 

\begin{equation*}
	\begin{cases}
		(\dot{b}_1,\dot{b}_2)\equiv(1,1),(1,3)(\text{mod }8,\text{mod }8), &\text{if}\ 2 \nmid b_1b_2\\
		\dot{b}_1\equiv3(\text{mod }8), &\text{if}\ 2\ |\ b_1,2\nmid b_2 \\
		(\dot{b}_1,\dot{b}_2)\equiv(5,5),(5,7)(\text{mod }8,\text{mod }8), &\text{if}\ 2\nmid b_1,2\ |\  b_2\\
		\dot{b}_1\equiv7(\text{mod }8), &\text{if}\ 2\ |\ b_1,2\ |\  b_2 
	\end{cases}.
\end{equation*}

 If  $\widetilde{n}\equiv 3({\rm mod }\ 8)$, $C_{\Lambda}(\mathbb{Q}_2)\neq \emptyset$ if and only if 

\begin{equation*}
	\begin{cases}
		(\dot{b}_1,\dot{b}_2)\equiv(1,1),(1,7)(\text{mod }8,\text{mod }8), &\text{if}\ 2 \nmid b_1b_2\\
		\dot{b}_1\equiv1(\text{mod }8), &\text{if}\ 2\ |\ b_1,2\nmid b_2 \\
		(\dot{b}_1,\dot{b}_2)\equiv(5,1),(5,7)(\text{mod }8,\text{mod }8), &\text{if}\ 2\nmid b_1,2\ |\  b_2\\
		\dot{b}_1\equiv5(\text{mod }8), &\text{if}\ 2\ |\ b_1,2\ |\  b_2 
	\end{cases}.
\end{equation*}

If  $\widetilde{n}\equiv 5({\rm mod }\ 8)$, $C_{\Lambda}(\mathbb{Q}_2)\neq \emptyset$ if and only if 

\begin{equation*}
	\begin{cases}
		(\dot{b}_1,\dot{b}_2)\equiv(1,1),(1,3)(\text{mod }8,\text{mod }8), &\text{if}\ 2 \nmid b_1b_2\\
		\dot{b}_1\equiv7(\text{mod }8), &\text{if}\ 2\ |\ b_1,2\nmid b_2 \\
		(\dot{b}_1,\dot{b}_2)\equiv(5,1),(5,3)(\text{mod }8,\text{mod }8), &\text{if}\ 2\nmid b_1,2\ |\  b_2\\
		\dot{b}_1\equiv3(\text{mod }8), &\text{if}\ 2\ |\ b_1,2\ |\  b_2 
	\end{cases}.
\end{equation*}

If  $\widetilde{n}\equiv 7({\rm mod }\ 8)$, $C_{\Lambda}(\mathbb{Q}_2)\neq \emptyset$ if and only if 

\begin{equation*}
	\begin{cases}
		(\dot{b}_1,\dot{b}_2)\equiv(1,1),(1,7)(\text{mod }8,\text{mod }8), &\text{if}\ 2 \nmid b_1b_2\\
		\dot{b}_1\equiv5(\text{mod }8), &\text{if}\ 2\ |\ b_1,2\nmid b_2 \\
		(\dot{b}_1,\dot{b}_2)\equiv(5,3),(5,5)(\text{mod }8,\text{mod }8), &\text{if}\ 2\nmid b_1,2\ |\  b_2\\
		\dot{b}_1\equiv1(\text{mod }8), &\text{if}\ 2\ |\ b_1,2\ |\  b_2 
	\end{cases}.
\end{equation*}

The conditions are equivalent to

\begin{equation*}
	\begin{cases}
		(1+[\frac{-1}{\widetilde{n}}])\gamma_2+[\frac{-1}{\dot{b}_1}]=0,\\
		(1+[\frac{2}{\widetilde{n}}])\gamma_2+\xi_2+[\frac{2}{\dot{b}_1}]=0,\\
		(1+[\frac{2}{\widetilde{n}}])\xi_2+[\frac{2}{\dot{b}_2}]+(1+[\frac{-1}{\widetilde{n}}])[\frac{-1}{\dot{b}_2}]=0,\ \ \text{if } \eta_2=0
	\end{cases}.
\end{equation*}

Note if $C_{\Lambda}(\mathbb{Q}_3)$ and $C_{\Lambda}(\mathbb{Q}_{p_i})\neq \emptyset$, we have $(1+[\frac{2}{\widetilde{n}}])\xi_2+[\frac{2}{\dot{b}_2}]+(1+[\frac{-1}{\widetilde{n}}])[\frac{-1}{\dot{b}_2}]=0$ by $[\frac{-3}{\widetilde{n}}]\xi_3+[\frac{-3}{\dot{b}_2}]=0$ and $[\frac{\dot{b}_2}{p_i}]+[\frac{-6\widetilde{n}p_i}{p_i}]y_i=0$. In other words, if we assume $C_{\Lambda}(\mathbb{Q}_3)$ and $C_{\Lambda}(\mathbb{Q}_{p_i})\neq \emptyset$, then $C_{\Lambda}(\mathbb{Q}_2)\neq \emptyset$ if and only if

\begin{equation*}
	\begin{cases}
		\gamma_1+(1+[\frac{-1}{\widetilde{n}}])\gamma_2+\gamma_3+\sum\limits_{i}[\frac{-1}{p_i}]x_i=0,\\
		(1+[\frac{2}{\widetilde{n}}])\gamma_2+\gamma_3+\xi_2+\sum\limits_{i}[\frac{2}{p_i}]x_i=0
	\end{cases}.
\end{equation*}
The Monsky matrix can be naturally derived from the above system of linear equations.
$\Box$

Then we have a corollary about pairty of 2-rank as following:

\begin{corollary}
	For $n=6\widetilde{n}=6p_1...p_t,p_i\neq2,3$, the rank of $2$-Selmer group $s_2(E_n)$ is odd. 
\end{corollary}

\emph{\bf Proof:}
Note that $s_2(E_n)=2t+6-r(M_n)$, then we only need consider the rank of $M_n$ which is denoted by $r(M_n)$. We have:

\begin{align*}
	r(D_1)&=r\begin{pmatrix}
		1 & 0 & 0 & 1 & 0 & 0 & O & O\\
		0 & 0 & 0 & 1 & 1+[\frac{-1}{\widetilde{n}}] & 1 & O & r_{-1}\\
		0 & 1 & 0 & 0 & 1+[\frac{2}{\widetilde{n}}] & 1 & O & r_{2}\\
		1 & 1 & [\frac{-3}{\widetilde{n}}] & 0 & 0 & 0 & r_{-3} & O\\
		0 & 0 & [\frac{-3}{\widetilde{n}}] & 1 & 1 & 1+[\frac{-3}{\widetilde{n}}] & O & r_{-3}\\
		0 & 0 & 0 & r_{-1}^T & r_{2}^T & r_{3}^T & D_{-3} & A+D_{6}\\
		r_{-1}^T & r_{2}^T & r_{3}^T & 0 & 0 & 0 & A+D_{-6} & O
	\end{pmatrix} \\
	&=r\begin{pmatrix}
		1 & 0 & 0 & 0 & 0 & 0 & O & O\\
		0 & 0 & 0 & 1 & 0 & 0 & O & O\\
		0 & 1 & 0 & 0 & 0 & 0 & O & O\\
		0 & 0 & 0 & 0 & 0 & 0 & r_{-3} & r_{-2}\\
		0 & 0 & 0 & 0 & 0 & 0 & O & O\\
		O & O & O & O & O & r_{-3}^T & D_{-3} & A+D_{6}+r_{-1}^Tr_{-1}\\
		O & O & O & O & O & r_{-2}^T & A+D_{-6} & r_{-1}^Tr_{-1}+r_{2}^Tr_{2}
	\end{pmatrix} \\ &\text{  by column and row transformations.} \\
	&=3+r\begin{pmatrix}
		0  & r_{-3} & r_{-2}\\
		r_{-3}^T & D_{-3}+r_{-3}^Tr_{-3} & A^T+D_{-6}+r_{-3}^Tr_{-2}\\
		r_{-2}^T & A+D_{-6}+r_{-2}^Tr_{-3} & r_{-1}^Tr_{-1}+r_{2}^Tr_{2}+r_{-2}^Tr_{-2}
	\end{pmatrix} \\
	&\text{  by } A+A^T+r_{-1}^Tr_{-1}+D_{-1}=O. 
\end{align*}

 Then we have $r(D_1)$ is odd. $\Box$

\subsection{\bf $\eta=-1$} 
Let 
$$A_4=\begin{pmatrix}
	1 & 0 & 0 & 0 & 0 & 0 & O & O\\
	0 & 0 & 0 & 0 & 1 & 0 & O & O\\
	0 & 1 & 0 & 0 & 0 & 0 & O & O\\
	0 & 0 & 0 & 1 & 0 & 1 & O & r_{-1}\\
	1 & 1 & 0 & 0 & 0 & [\frac{-3}{\widetilde{n}}] & r_{-3} & O\\
	0 & 0 & 1 & 0 & 0 & 0 & O & O\\
	O & O & O & r_{-1}^T & r_{2}^T & r_{3}^T & D_{-3} & A+D_{-1}\\
	r_{-1}^T & r_{2}^T & r_{3}^T & O & O & O & A & O
\end{pmatrix},$$

and
$$A_5=\begin{pmatrix}
	1 & 0 & 0 & 0 & 0 & 0 & O & O\\
	0 & 1 & 0 & 1 & 0 & 1 & O & r_{-1}\\
	0 & 0 & 0 & 0 & 1 & 0 & O & O\\
	1 & 1 & 1 & 0 & 0 & 1 & r_{-1} & r_2\\
	1 & 1 & 0 & 0 & 0 & [\frac{-3}{\widetilde{n}}] & r_{-3} & O\\
	0 & 0 & 1 & 0 & 0 & 0 & O & O\\
	O & O & O & r_{-1}^T & r_{2}^T & r_{3}^T & D_{-3} & A+D_{-1}\\
	r_{-1}^T & r_{2}^T & r_{3}^T & O & O & O & A & O
\end{pmatrix},$$

and 

$$A_6=\begin{pmatrix}
	1 & 0 & 0 & 0 & 0 & 0 & O & O\\
	0 & 1 & 0 & 0 & 0 & 0 & O & O\\
	0 & 0 & 0 & 0 & 1 & 0 & O & O\\
	1 & 0 & 1 & 0 & 0 & 0 & r_{-1} & O\\
	1 & 1 & 0 & 0 & 0 & [\frac{-3}{\widetilde{n}}] & r_{-3} & O\\
	0 & 0 & 1 & 0 & 0 & 0 & O & O\\
	O & O & O & r_{-1}^T & r_{2}^T & r_{3}^T & D_{-3} & A+D_{-1}\\
	r_{-1}^T & r_{2}^T & r_{3}^T & O & O & O & A & O
\end{pmatrix}.$$

\begin{theorem}
	For $n=-\widetilde{n}=-p_1...p_t,p_i \neq 2,3$, then we have Monsky matrix 
	\begin{equation*}
		M_n=\begin{cases}
			A_4, &\text{if } \widetilde{n}\equiv 3({\rm mod }\ 8),\\
			A_5, &\text{if } \widetilde{n}\equiv 7({\rm mod }\ 8), \\
			A_6, &\text{if } \widetilde{n}\equiv 1({\rm mod }\ 4).
		\end{cases}
	\end{equation*} 
\end{theorem}

\emph{\bf Proof:}We list all the equivalent conditions and transform them to additional Legendre symbols to build system of linear equations. 

We know that $C_{\Lambda}(\mathbb{A})\neq \emptyset$ if and only if $C_{\Lambda}(\mathbb{Q}_S)\neq \emptyset$, i.e.,

1) $C_{\Lambda}(\mathbb{R})\neq \emptyset$ if and only if 		$b_2>0$ that is equivalent to $\xi_1=0$;

2) $C_{\Lambda}(\mathbb{Q}_{p_i})\neq \emptyset$ for $i=1,..,t$ if and only if 
\begin{equation*}
	\begin{cases}
		\left( \frac{\dot{b}_1}{p_i}\right)=\left( \frac{\dot{b}_2}{p_i}\right)=1, &\text{if}\ p_i \nmid b_1b_2\\
		\left( \frac{\dot{n}\dot{b}_1}{p_i}\right)=\left( \frac{\dot{b}_2}{p_i}\right)=1, &\text{if}\ p_i\ |\ b_1,p_i\nmid b_2 \\
		\left( \frac{-3\dot{b}_1}{p_i}\right)=\left( \frac{-\dot{n}\dot{b}_2}{p_i}\right)=1, &\text{if}\ p_i\nmid b_1,p_i\ |\  b_2\\
		\left( \frac{-3\dot{n}\dot{b}_1}{p_i}\right)=\left( \frac{-\dot{n}\dot{b}_2}{p_i}\right)=1, &\text{if}\ p_i\ |\ b_1,p_i\ |\  b_2 
	\end{cases}
\end{equation*}

by Hensel's Lemma. The conditions are equivalent to
\begin{equation*}
	\begin{cases}
		[\frac{\dot{b}_1}{p_i}]+[\frac{-3}{p_i}]y_i+[\frac{-\widetilde{n}p_i}{p_i}]x_i=0 \\
		[\frac{\dot{b}_2}{p_i}]+[\frac{\widetilde{n}p_i}{p_i}]y_i=0 
	\end{cases};
\end{equation*}
In other words, we have
\begin{equation*}
	\begin{cases}
		[\frac{-1}{p_i}]\gamma_1+[\frac{2}{p_i}]\gamma_2+[\frac{3}{p_i}]\gamma_3+[\frac{-3}{p_i}]y_i+\sum\limits_{j\neq i}[\frac{p_j}{p_i}]x_j+[\frac{-\widetilde{n}p_i}{p_i}]x_i=0 \\
		[\frac{-1}{p_i}]\xi_1+[\frac{2}{p_i}]\xi_2+[\frac{3}{p_i}]\xi_3+\sum\limits_{j\neq i}[\frac{p_j}{p_i}]y_j+[\frac{\widetilde{n}p_i}{p_i}]y_i=0 
	\end{cases};
\end{equation*}

3) $C_{\Lambda}(\mathbb{Q}_3)\neq \emptyset$ if and only if $3\nmid b_2$ and
\begin{equation*}
	\begin{cases}
		\left( \frac{\dot{b}_2}{3}\right)=1, &\text{if}\ 3 \nmid b_1\\
		\left( \frac{-\dot{n}\dot{b}_2}{3}\right)=1, &\text{if}\ 3\ |\ b_1 
	\end{cases}
\end{equation*}

by Hensel's Lemma. The conditions are equivalent to

\begin{equation*}
	\begin{cases}
		[\frac{-3}{\widetilde{n}}]\gamma_3+[\frac{-3}{\dot{b}_2}]=0\\
		\xi_3=0
	\end{cases};
\end{equation*}
In other words, we have 
\begin{equation*}
	\begin{cases}
		\xi_1+\xi_2+[\frac{-3}{\widetilde{n}}]\gamma_3+\sum\limits_{i}[\frac{-3}{p_i}]y_i=0\\
		\xi_3=0
	\end{cases};
\end{equation*}
4) $C_{\Lambda}(\mathbb{Q}_2)\neq \emptyset$. 
If $\widetilde{n}\equiv 3({\rm mod }\ 8)$, $C_{\Lambda}(\mathbb{Q}_2)\neq \emptyset$ if and only if $2\nmid b_1b_2$ and $\dot{b}_1\equiv 1({\rm mod }\ 4)$; 
if  $\widetilde{n}\equiv 7({\rm mod }\ 8)$, $C_{\Lambda}(\mathbb{Q}_2)\neq \emptyset$ if and only if (i):$2\nmid b_1b_2$ and $(\dot{b}_1,\dot{b}_2)\equiv (1,1),(1,5),(5,7),(5,3)({\rm mod }\ 8,{\rm mod }\ 8)$ or (ii):$2\nmid b_1,2\ |\ b_2$ and $(\dot{b}_1,\dot{b}_2)\equiv (7,3),$ $(3,1)({\rm mod }\ 8,{\rm mod }\ 4)$;
if $\widetilde{n}\equiv 1({\rm mod }\ 4)$, $C_{\Lambda}(\mathbb{Q}_2)\neq \emptyset$ if and only if $2\nmid b_1b_2$ and $\dot{b}_2\equiv 1({\rm mod }\ 4)$ by Hensel's Lemma. The conditions are  equivalent to

\begin{equation*}
	\begin{cases}
		\xi_2=0, \ \ \ \ \gamma_2=0,\\
		\gamma_1+\gamma_3+\sum\limits_{i}[\frac{-1}{p_i}]x_i=0
	\end{cases} \text{if } \widetilde{n}\equiv 3({\rm mod }\ 8),
\end{equation*}
and
\begin{equation*}
	\begin{cases}
		\xi_2+\gamma_1+\gamma_3+\sum\limits_{i}[\frac{-1}{p_i}]x_i=0, \ \ \ \ \gamma_2=0,\\
		\xi_1+\xi_2+\xi_3+\gamma_3+\sum\limits_{i}[\frac{-1}{p_i}]y_i+\sum\limits_{i}[\frac{2}{p_i}]x_i=0
	\end{cases} \text{if } \widetilde{n}\equiv 7({\rm mod }\ 8),
\end{equation*}
and
\begin{equation*}
	\begin{cases}
		\xi_2=0, \ \ \ \ \gamma_2=0,\\
		\xi_1+\xi_3+\sum\limits_{i}[\frac{-1}{p_i}]y_i=0
	\end{cases} \text{if } \widetilde{n}\equiv 1({\rm mod }\ 4).
\end{equation*}
The Monsky matrix can be naturally derived from the above system of linear equations.
$\Box$

Then we have a corollary about pairty of 2-rank as following:

\begin{corollary}
	For $n=-\widetilde{n}=-p_1...p_t,p_i\neq2,3$, the rank of $2$-Selmer group $s_2(E_n)$ is even (resp. odd) if and only if $n\equiv11,13,17,23$ $(resp.\ 1,5,7$ $,19)\pmod{24}$. 
\end{corollary}

\emph{\bf Proof:}
Note that $s_2(E_n)=2t+6-r(M_n)$, then we only need consider the rank of $M_n$ which is denoted by $r(M_n)$.

For $\widetilde{n}\equiv3(\text{mod 8})$, we have:
\begin{align*}
	r(A_4)&=r\begin{pmatrix}
		1 & 0 & 0 & 0 & 0 & 0 & O & O\\
		0 & 0 & 0 & 0 & 1 & 0 & O & O\\
		0 & 1 & 0 & 0 & 0 & 0 & O & O\\
		0 & 0 & 0 & 1 & 0 & 1 & O & r_{-1}\\
		1 & 1 & 0 & 0 & 0 & [\frac{-3}{\widetilde{n}}] & r_{-3} & O\\
		0 & 0 & 1 & 0 & 0 & 0 & O & O\\
		O & O & O & r_{-1}^T & r_{2}^T & r_{3}^T & D_{-3} & A+D_{-1}\\
		r_{-1}^T & r_{2}^T & r_{3}^T & O & O & O & A & O
	\end{pmatrix} \\
	&=r\begin{pmatrix}
		1 & 0 & 0 & 0 & 0 & 0 & O & O\\
		0 & 1 & 0 & 0 & 0 & 0 & O & O\\
		0 & 0 & 0 & 0 & 1 & 0 & O & O\\
		0 & 0 & 0 & 1 & 0 & 1 & O & r_{-1}\\
		0 & 0 & 0 & 0 & 0 & [\frac{-3}{\widetilde{n}}] & r_{-3} & O\\
		0 & 0 & 1 & 0 & 0 & 0 & O & O\\
		O & O & O & O & O & r_{-3}^T & D_{-3} & A+D_{-1}+r_{-1}^Tr_{-1}\\
		O & O & O & O & O & O & A & O
	\end{pmatrix} \\ &\text{  by row transformations.} \\
	&=r\begin{pmatrix}
		1 & 0 & 0 & 0 & 0 & 0 & O & O\\
		0 & 1 & 0 & 0 & 0 & 0 & O & O\\
		0 & 0 & 0 & 0 & 1 & 0 & O & O\\
		0 & 0 & 0 & 1 & 0 & 0 & O & O\\
		0 & 0 & 1 & 0 & 0 & 0 & O & O\\
		0 & 0 & 0 & 0 & 0 & [\frac{-3}{\widetilde{n}}] & r_{-3} & O\\
		O & O & O & O & O & r_{-3}^T & D_{-3} & A+D_{-1}+r_{-1}^Tr_{-1}\\
		O & O & O & O & O & O & A & O
	\end{pmatrix} \\ &\text{  by column transformations.} \\
	&=5+r\begin{pmatrix}
		[\frac{-3}{\widetilde{n}}]  & r_{-3} & O\\
		r_{-3}^T  & D_{-3} & A^T\\
		O  & A & O
	\end{pmatrix} \\
	&\text{  by } A+A^T+r_{-1}^Tr_{-1}+D_{-1}=O. 
\end{align*}

If $[\frac{-3}{\widetilde{n}}]=0$, i.e. $n\equiv5(\text{mod }24)$, $$r(A_4)=5+r\begin{pmatrix}
	0  & r_{-3} & O\\
	r_{-3}^T  & D_{-3}+r_{-3}^Tr_{-3} & A^T\\
	O  & A & O
\end{pmatrix} $$ is odd; and if $[\frac{-3}{\widetilde{n}}]=1$, i.e. $n\equiv13(\text{mod }24)$, $$r(A_4)=6+r\begin{pmatrix}
 D_{-3}+r_{-3}^Tr_{-3} & A^T\\
 A & O
\end{pmatrix} $$ is even.

For $\widetilde{n}\equiv7(\text{mod 8})$, we have:
\begin{align*}
	r(A_5)&=r\begin{pmatrix}
		1 & 0 & 0 & 0 & 0 & 0 & O & O\\
		0 & 1 & 0 & 1 & 0 & 1 & O & r_{-1}\\
		0 & 0 & 0 & 0 & 1 & 0 & O & O\\
		1 & 1 & 1 & 0 & 0 & 1 & r_{-1} & r_2\\
		1 & 1 & 0 & 0 & 0 & [\frac{-3}{\widetilde{n}}] & r_{-3} & O\\
		0 & 0 & 1 & 0 & 0 & 0 & O & O\\
		O & O & O & r_{-1}^T & r_{2}^T & r_{3}^T & D_{-3} & A+D_{-1}\\
		r_{-1}^T & r_{2}^T & r_{3}^T & O & O & O & A & O
	\end{pmatrix} \\
	&=r\begin{pmatrix}
		1 & 0 & 0 & 0 & 0 & 0 & O & O\\
		0 & 0 & 0 & 0 & 1 & 0 & O & O\\
		0 & 1 & 0 & 1 & 0 & 1 & O & r_{-1}\\
		0 & 1 & 0 & 0 & 0 & 1 & r_{-1} & r_{2}\\
		0 & 0 & 0 & 0 & 0 & [\frac{3}{\widetilde{n}}] & r_{3} & r_{2}\\
		0 & 0 & 1 & 0 & 0 & 0 & O & O\\
		O & O & O & O & O & r_{3}^T & D_{-3}+r_{-1}^Tr_{-1} & A+D_{-1}+r_{-1}^Tr_{-2}\\
		O & O & O & O & O & r_{2}^T & A+r_{2}^Tr_{-1} & r_{2}^Tr_{2}
	\end{pmatrix} \\ &\text{  by row transformations.} \\
	&=r\begin{pmatrix}
		1 & 0 & 0 & 0 & 0 & 0 & O & O\\
		0 & 0 & 0 & 0 & 1 & 0 & O & O\\
		0 & 0 & 0 & 1 & 0 & 0 & O & O\\
		0 & 1 & 0 & 0 & 0 & 0 & O & O\\
		0 & 0 & 1 & 0 & 0 & 0 & O & O\\
		0 & 0 & 0 & 0 & 0 & [\frac{3}{\widetilde{n}}] & r_{3} & r_{2}\\
		O & O & O & O & O & r_{3}^T & D_{-3}+r_{-1}^Tr_{-1} & A^T+r_{-1}^Tr_{2}\\
		O & O & O & O & O & r_{2}^T & A+r_{2}^Tr_{-1} & r_{2}^Tr_{2}
	\end{pmatrix} \\ &\text{  by column transformations.} \\
	&=5+r\begin{pmatrix}
 [\frac{3}{\widetilde{n}}] & r_{3} & r_{2}\\
 r_{3}^T & D_{-3}+r_{-1}^Tr_{-1} & A^T+r_{-1}^Tr_{2}\\
r_{2}^T & A+r_{2}^Tr_{-1} & r_{2}^Tr_{2}
	\end{pmatrix} \\
	&\text{  by } A+A^T+r_{-1}^Tr_{-1}+D_{-1}=O. 
\end{align*}

If $[\frac{3}{\widetilde{n}}]=0$, i.e. $n\equiv1(\text{mod }24)$, $$r(A_5)=5+r\begin{pmatrix}
	0 & r_{3} & r_{2}\\
	r_{3}^T & D_{-3}+r_{-1}^Tr_{-1}+r_{3}^Tr_{3} & A^T+r_{-3}^Tr_{2}\\
	r_{2}^T & A+r_{2}^Tr_{-3} & O
\end{pmatrix} $$ is odd; and if $[\frac{3}{\widetilde{n}}]=1$, i.e. $n\equiv17(\text{mod }24)$, $$r(A_5)=6+r\begin{pmatrix} 	
 D_{-3}+r_{-1}^Tr_{-1}+r_{3}^Tr_{3} & A^T+r_{-3}^Tr_{2}\\
A+r_{2}^Tr_{-3} & O
\end{pmatrix} $$ is even.

For $\widetilde{n}\equiv1(\text{mod 4})$, we have:
\begin{align*}
	r(A_6)&=r\begin{pmatrix}
		1 & 0 & 0 & 0 & 0 & 0 & O & O\\
		0 & 1 & 0 & 0 & 0 & 0 & O & O\\
		0 & 0 & 0 & 0 & 1 & 0 & O & O\\
		1 & 0 & 1 & 0 & 0 & 0 & r_{-1} & O\\
		1 & 1 & 0 & 0 & 0 & [\frac{-3}{\widetilde{n}}] & r_{-3} & O\\
		0 & 0 & 1 & 0 & 0 & 0 & O & O\\
		O & O & O & r_{-1}^T & r_{2}^T & r_{3}^T & D_{-3} & A+D_{-1}\\
		r_{-1}^T & r_{2}^T & r_{3}^T & 0 & 0 & 0 & A & O
	\end{pmatrix} \\
	&=r\begin{pmatrix}
		1 & 0 & 0 & 0 & 0 & 0 & O & O\\
		0 & 1 & 0 & 0 & 0 & 0 & O & O\\
		0 & 0 & 0 & 0 & 1 & 0 & O & O\\
		0 & 0 & 0 & 0 & 0 & 0 & O & O\\
		0 & 0 & 0 & 0 & 0 & [\frac{3}{\widetilde{n}}] & r_{-3} & O\\
		0 & 0 & 1 & 0 & 0 & 0 & O & O\\
		O & O & O & r_{-1}^T & O & r_{3}^T & D_{-3} & A+D_{-1}\\
		O & O & O & O & O & O & A+r_{-1}^Tr_{-1} & O
	\end{pmatrix} \\ &\text{  by row transformations.} \\
	&=4+r\begin{pmatrix}
 [\frac{3}{\widetilde{n}}] & r_{-3} & O\\
 r_{-3}^T & D_{-3} & A+D_{-1}\\
 O & A^T+D_{-1} & O
	\end{pmatrix} \\ 
&\text{  by column transformations and }
A+A^T+r_{-1}^Tr_{-1}+D_{-1}=O.
\end{align*}

If $[\frac{3}{\widetilde{n}}]=0$, i.e. $n\equiv11,23(\text{mod }24)$, $$r(A_6)=4+r\begin{pmatrix}
	0  & r_{-3} & O\\
	r_{-3}^T & D_{-3}+r_{-3}^Tr_{-3} & A+D_{-1}\\
	O & A^T+D_{-1} & O
\end{pmatrix} $$ is even; and if $[\frac{3}{\widetilde{n}}]=1$, i.e. $n\equiv7,19(\text{mod }24)$, $$r(A_6)=5+r\begin{pmatrix}
	D_{-3}+r_{-3}^Tr_{-3} & A+D_{-1}\\
	A^T+D_{-1} & O
\end{pmatrix} $$ is odd. $\Box$

Now we prove a result about non $2\pi/3$-congruent numbers:

\begin{corollary}
	For $n=-\widetilde{n}=-p_1...p_t\equiv13(\text{mod }24),p_i\neq2,3$, and $r_4(n)=0$, then $s_2(E_n)=2$. 
\end{corollary}

\emph{\bf Proof:} 
For $n\equiv13(\text{mod }24)$, we have

$$ s_2(E_n)=2t+6-r(M_n)=2t-r\begin{pmatrix}
	D_{-3}+r_{-3}^Tr_{-3} & A^T\\
	A & O
\end{pmatrix}.$$ 
Then $s_2\leq2t-2r(A)=2t-2r(R(n))=2+2r_4(n)=2.$

Note that $s_2(E_n)\geq2$ by the fact that $E_n$ have $4$ elements of $2$-torsion, we have $s_2(E_n)=2$. $\Box$

\subsection{\bf $\eta=-2$} 
Let 
$$B_2=\begin{pmatrix}
	1 & 0 & 0 & 0 & 0 & 0 & O & O\\
	0 & 0 & 0 & 1 & 1+[\frac{-1}{\widetilde{n}}] & 1 & O & r_{-1}\\
	0 & 1 & 0 & 0 & [\frac{2}{\widetilde{n}}] & 1 & O & r_{2}\\
	1 & 1 & 0 & 0 & 0 & 1+[\frac{-3}{\widetilde{n}}] & r_{-3} & O\\
	0 & 0 & 1 & 0 & 0 & 0 & O & O\\
	O & O & O & r_{-1}^T & r_{2}^T & r_{3}^T & D_{-3} & A+D_{-2}\\
	r_{-1}^T & r_{2}^T & r_{3}^T & O & O & O & A+D_{2} & O
\end{pmatrix}.$$

\begin{theorem}
	For $n=-2\widetilde{n}=-2p_1...p_t,p_i \neq 2,3$, then we have Monsky matrix $M_n=B_2$.
\end{theorem}

\emph{\bf Proof:}We list all the equivalent conditions and transform them to additional Legendre symbols to build system of linear equations. 

We know that $C_{\Lambda}(\mathbb{A})\neq \emptyset$ if and only if $C_{\Lambda}(\mathbb{Q}_S)\neq \emptyset$, i.e.,

1) $C_{\Lambda}(\mathbb{R})\neq \emptyset$ if and only if 		$b_2>0$ that is equivalent to $\xi_1=0$;

2) $C_{\Lambda}(\mathbb{Q}_{p_i})\neq \emptyset$ for $i=1,..,t$ if and only if 
\begin{equation*}
	\begin{cases}
		\left( \frac{\dot{b}_1}{p_i}\right)=\left( \frac{\dot{b}_2}{p_i}\right)=1, &\text{if}\ p_i \nmid b_1b_2\\
		\left( \frac{\dot{n}\dot{b}_1}{p_i}\right)=\left( \frac{\dot{b}_2}{p_i}\right)=1, &\text{if}\ p_i\ |\ b_1,p_i\nmid b_2 \\
		\left( \frac{-3\dot{b}_1}{p_i}\right)=\left( \frac{-\dot{n}\dot{b}_2}{p_i}\right)=1, &\text{if}\ p_i\nmid b_1,p_i\ |\  b_2\\
		\left( \frac{-3\dot{n}\dot{b}_1}{p_i}\right)=\left( \frac{-\dot{n}\dot{b}_2}{p_i}\right)=1, &\text{if}\ p_i\ |\ b_1,p_i\ |\  b_2 
	\end{cases}
\end{equation*}

by Hensel's Lemma. The conditions are equivalent to
\begin{equation*}
	\begin{cases}
		[\frac{\dot{b}_1}{p_i}]+[\frac{-3}{p_i}]y_i+[-\frac{2\widetilde{n}p_i}{p_i}]x_i=0 \\
		[\frac{\dot{b}_2}{p_i}]+[\frac{2\widetilde{n}p_i}{p_i}]y_i=0 
	\end{cases};
\end{equation*}
In other words, we have

\begin{equation*}
	\begin{cases}
		[\frac{-1}{p_i}]\gamma_1+[\frac{2}{p_i}]\gamma_2+[\frac{3}{p_i}]\gamma_3+[\frac{-3}{p_i}]y_i+\sum\limits_{j\neq i}[\frac{p_j}{p_i}]x_j+[\frac{-2\widetilde{n}p_i}{p_i}]x_i=0 \\
		[\frac{-1}{p_i}]\xi_1+[\frac{2}{p_i}]\xi_2+[\frac{3}{p_i}]\xi_3+\sum\limits_{j\neq i}[\frac{p_j}{p_i}]y_j+[\frac{2\widetilde{n}p_i}{p_i}]y_i=0 
	\end{cases};
\end{equation*}

3) $C_{\Lambda}(\mathbb{Q}_3)\neq \emptyset$ if and only if $3\nmid b_2$ and
\begin{equation*}
	\begin{cases}
		\left( \frac{\dot{b}_2}{3}\right)=1, &\text{if}\ 3 \nmid b_1\\
		\left( \frac{-\dot{n}\dot{b}_2}{3}\right)=1, &\text{if}\ 3\ |\ b_1 
	\end{cases}
\end{equation*}

by Hensel's Lemma. The conditions are equivalent to
\begin{equation*}
	\begin{cases}
		(1+[\frac{-3}{\widetilde{n}}])\gamma_3+[\frac{-3}{\dot{b}_2}]=0\\
		\xi_3=0
	\end{cases}.
\end{equation*}
In other words, we have 
\begin{equation*}
	\begin{cases}
		\xi_1+\xi_2+(1+[\frac{-3}{\widetilde{n}}])\gamma_3+\sum\limits_{i}[\frac{-3}{p_i}]y_i=0\\
		\xi_3=0
	\end{cases};
\end{equation*}

4) $C_{\Lambda}(\mathbb{Q}_2)\neq \emptyset$. If  $\widetilde{n}\equiv 1({\rm mod }\ 8)$, $C_{\Lambda}(\mathbb{Q}_2)\neq \emptyset$ if and only if 

\begin{equation*}
	\begin{cases}
		(\dot{b}_1,\dot{b}_2)\equiv(1,1),(1,3)(\text{mod }8,\text{mod }8), &\text{if}\ 2 \nmid b_1b_2\\
		\dot{b}_1\equiv7(\text{mod }8), &\text{if}\ 2\ |\ b_1,2\nmid b_2 \\
		(\dot{b}_1,\dot{b}_2)\equiv(5,1),(5,3)(\text{mod }8,\text{mod }8), &\text{if}\ 2\nmid b_1,2\ |\  b_2\\
		\dot{b}_1\equiv3(\text{mod }8), &\text{if}\ 2\ |\ b_1,2\ |\  b_2 
	\end{cases}.
\end{equation*}

If  $\widetilde{n}\equiv 3({\rm mod }\ 8)$, $C_{\Lambda}(\mathbb{Q}_2)\neq \emptyset$ if and only if 

\begin{equation*}
	\begin{cases}
		(\dot{b}_1,\dot{b}_2)\equiv(1,1),(1,7)(\text{mod }8,\text{mod }8), &\text{if}\ 2 \nmid b_1b_2\\
		\dot{b}_1\equiv5(\text{mod }8), &\text{if}\ 2\ |\ b_1,2\nmid b_2 \\
		(\dot{b}_1,\dot{b}_2)\equiv(5,3),(5,5)(\text{mod }8,\text{mod }8), &\text{if}\ 2\nmid b_1,2\ |\  b_2\\
		\dot{b}_1\equiv1(\text{mod }8), &\text{if}\ 2\ |\ b_1,2\ |\  b_2 
	\end{cases}.
\end{equation*}

If  $\widetilde{n}\equiv 5({\rm mod }\ 8)$, $C_{\Lambda}(\mathbb{Q}_2)\neq \emptyset$ if and only if 

\begin{equation*}
	\begin{cases}
		(\dot{b}_1,\dot{b}_2)\equiv(1,1),(1,3)(\text{mod }8,\text{mod }8), &\text{if}\ 2 \nmid b_1b_2\\
		\dot{b}_1\equiv3(\text{mod }8), &\text{if}\ 2\ |\ b_1,2\nmid b_2 \\
		(\dot{b}_1,\dot{b}_2)\equiv(5,5),(5,7)(\text{mod }8,\text{mod }8), &\text{if}\ 2\nmid b_1,2\ |\  b_2\\
		\dot{b}_1\equiv7(\text{mod }8), &\text{if}\ 2\ |\ b_1,2\ |\  b_2 
	\end{cases}.
\end{equation*}

If  $\widetilde{n}\equiv 7({\rm mod }\ 8)$, $C_{\Lambda}(\mathbb{Q}_2)\neq \emptyset$ if and only if 

\begin{equation*}
	\begin{cases}
		(\dot{b}_1,\dot{b}_2)\equiv(1,1),(1,7)(\text{mod }8,\text{mod }8), &\text{if}\ 2 \nmid b_1b_2\\
		\dot{b}_1\equiv1(\text{mod }8), &\text{if}\ 2\ |\ b_1,2\nmid b_2 \\
		(\dot{b}_1,\dot{b}_2)\equiv(5,1),(5,7)(\text{mod }8,\text{mod }8), &\text{if}\ 2\nmid b_1,2\ |\  b_2\\
		\dot{b}_1\equiv5(\text{mod }8), &\text{if}\ 2\ |\ b_1,2\ |\  b_2 
	\end{cases}.
\end{equation*}

The conditions are equivalent to

\begin{equation*}
	\begin{cases}
		(1+[\frac{-1}{\widetilde{n}}])\gamma_2+[\frac{-1}{\dot{b}_1}]=0,\\
		[\frac{2}{\widetilde{n}}]\gamma_2+\xi_2+[\frac{2}{\dot{b}_1}]=0,\\
		[\frac{2}{\widetilde{n}}]\xi_2+[\frac{2}{\dot{b}_2}]+(1+[\frac{-1}{\widetilde{n}}])[\frac{-1}{\dot{b}_2}]=0,\ \ \text{if } \eta_2=0
	\end{cases}.
\end{equation*}

Note if $C_{\Lambda}(\mathbb{Q}_3)$ and $C_{\Lambda}(\mathbb{Q}_{p_i})\neq \emptyset$, we have $[\frac{2}{\widetilde{n}}]\xi_2+[\frac{2}{\dot{b}_2}]+(1+[\frac{-1}{\widetilde{n}}])[\frac{-1}{\dot{b}_2}]=0$ by $\xi_3=0$ and $[\frac{\dot{b}_2}{p_i}]+[\frac{2\widetilde{n}p_i}{p_i}]y_i=0$. In other words, if we assume $C_{\Lambda}(\mathbb{Q}_3)$ and $C_{\Lambda}(\mathbb{Q}_{p_i})\neq \emptyset$, then $C_{\Lambda}(\mathbb{Q}_2)\neq \emptyset$ if and only if

\begin{equation*}
	\begin{cases}
		\gamma_1+(1+[\frac{-1}{\widetilde{n}}])\gamma_2+\gamma_3+\sum\limits_{i}[\frac{-1}{p_i}]x_i=0,\\
		[\frac{2}{\widetilde{n}}]\gamma_2+\gamma_3+\xi_2+\sum\limits_{i}[\frac{2}{p_i}]x_i=0
	\end{cases}.
\end{equation*}
The Monsky matrix can be naturally derived from the above system of linear equations.
$\Box$

Then we have a corollary about pairty of 2-rank as following:

\begin{corollary}
	For $n=-2\widetilde{n}=-2p_1...p_t,p_i\neq2,3$, the rank of $2$-Selmer group $s_2(E_n)$ is even (resp. odd) if and only if $n\equiv10,22(resp.\ 2,14)\pmod{24}$, i.e., $\widetilde{n}\equiv1(resp.\ 2)\pmod{3}$. 
\end{corollary}

\emph{\bf Proof:}
Note that $s_2(E_n)=2t+6-r(M_n)$, then we only need consider the rank of $M_n$ which is denoted by $r(M_n)$. We have:

\begin{align*}
	r(B_2)&=r\begin{pmatrix}
		1 & 0 & 0 & 0 & 0 & 0 & O & O\\
		0 & 0 & 0 & 1 & 1+[\frac{-1}{\widetilde{n}}] & 1 & O & r_{-1}\\
		0 & 1 & 0 & 0 & [\frac{2}{\widetilde{n}}] & 1 & O & r_{2}\\
		1 & 1 & 0 & 0 & 0 & 1+[\frac{-3}{\widetilde{n}}] & r_{-3} & O\\
		0 & 0 & 1 & 0 & 0 & 0 & O & O\\
		O & O & O & r_{-1}^T & r_{2}^T & r_{3}^T & D_{-3} & A+D_{-2}\\
		r_{-1}^T & r_{2}^T & r_{3}^T & O & O & O & A+D_{2} & O
	\end{pmatrix} \\
	&=r\begin{pmatrix}
		1 & 0 & 0 & 1 & 0 & 0 & O & O\\
		0 & 0 & 0 & 1 & 1+[\frac{-1}{\widetilde{n}}] & 1 & O & r_{-1}\\
		0 & 1 & 0 & 0 & [\frac{2}{\widetilde{n}}] & 1 & O & r_{2}\\
		0 & 0 & 0 & 0 & [\frac{2}{\widetilde{n}}] & [\frac{-3}{\widetilde{n}}] & r_{-3} & r_{2}\\
		0 & 0 & 1 & 0 & 0 & 0 & O & O\\
		O & O & O & O & [\frac{-1}{\widetilde{n}}]r_{-1}^T+r_{-2}^T & r_{-3}^T & D_{-3} & A+D_{-2}+r_{-1}^Tr_{-1}\\
		O & O & O & O & [\frac{2}{\widetilde{n}}]r_{-2}^T & r_{2}^T & A+D_{2} & r_{2}^Tr_{2}
	\end{pmatrix} \\ &\text{  by row transformations.} \\
	&=r\begin{pmatrix}
		1 & 0 & 0 & 0 & 0 & 0 & O & O\\
		0 & 0 & 0 & 1 & 0 & 0 & O & O\\
		0 & 1 & 0 & 0 & 0 & 0 & O & O\\
		0 & 0 & 0 & 0 & 0 & [\frac{-3}{\widetilde{n}}] & r_{-3} & r_{2}\\
		0 & 0 & 1 & 0 & 0 & 0 & O & O\\
		O & O & O & O & O & r_{-3}^T & D_{-3} & A+D_{-2}+r_{-1}^Tr_{-1}\\
		O & O & O & O & O & r_{2}^T & A+D_{2} & r_{2}^Tr_{2}
	\end{pmatrix} \\ &\text{  by column transformations.} \\
	&=4+r\begin{pmatrix}
		[\frac{-3}{\widetilde{n}}]  & r_{-3} & r_{2}\\
		r_{-3}^T & D_{-3} & A^T+D_{2}\\
		r_{2}^T & A+D_{2} & r_{2}^Tr_{2}
	\end{pmatrix} \\
	&\text{  by } A+A^T+r_{-1}^Tr_{-1}+D_{-1}=O. 
\end{align*}

If $[\frac{-3}{\widetilde{n}}]=0$, i.e. $n\equiv10,22(\text{mod }24)$, $$r(B_2)=4+r\begin{pmatrix}
	0  & r_{-3} & r_{2}\\
	r_{-3}^T & D_{-3}+r_{-3}^Tr_{-3} & A^T+D_{2}+r_{-3}^Tr_{2}\\
	r_{2}^T & A+D_{2}+r_{2}^Tr_{-3} & O
\end{pmatrix} $$ is even; and if $[\frac{-3}{\widetilde{n}}]=1$, i.e. $n\equiv2,14(\text{mod }24)$, $$r(B_2)=5+r\begin{pmatrix}
	D_{-3}+r_{-3}^Tr_{-3} & A^T+D_{2}+r_{-3}^Tr_{2}\\
	A+D_{2}+r_{2}^Tr_{-3} & O
\end{pmatrix} $$ is odd. $\Box$

Now we prove a result about non $2\pi/3$-congruent numbers:

\begin{corollary}
	For $n=-2\widetilde{n}=-2p_1...p_t\equiv10,22\pmod{24},p_i\neq2,3$, if $r_4(n)=0$, then $s_2(E_n)=2$. 
\end{corollary}

\emph{\bf Proof:} 
For $n\equiv10,22(\text{mod }24)$, we have

$$ s_2(E_n)=2t+6-r(M_n)=2t+2-r\begin{pmatrix}
	0  & r_{-3} & r_{2}\\
	r_{-3}^T & D_{-3}+r_{-3}^Tr_{-3} & A^T+D_{2}+r_{-3}^Tr_{2}\\
	r_{2}^T & A+D_{2}+r_{2}^Tr_{-3} & O
\end{pmatrix}.$$ Then $s_2\leq2t+2-2r(	r_{2}^T \ A+D_{2}+r_{2}^Tr_{-3})=2t+2-2r(r_{2}^T \ A+D_{2})=2t+2-2r(R(n))=2+2r_4(n)=2.$

Note that $s_2(E_n)\geq2$ by the fact $E_n$ have $4$ elements of $2$-torsion, we have $s_2(E_n)=2$. $\Box$

\subsection{\bf $\eta=-3$} Let 
$$C_4=\begin{pmatrix}
	1 & 0 & 0 & 0 & 0 & 0 & O & O\\
	0 & 0 & 0 & 0 & 1 & 0 & O & O\\
	0 & 1 & 0 & 0 & 0 & 0 & O & O\\
	0 & 0 & 0 & 1 & 0 & 1 & O & r_{-1}\\
	1 & 1 & [\frac{-3}{\widetilde{n}}] & 0 & 0 & 0 & r_{-3} & O\\
	0 & 0 & [\frac{-3}{\widetilde{n}}] & 1 & 1 & 1+[\frac{-3}{\widetilde{n}}] & O & r_{-3}\\
	O & O & O & r_{-1}^T & r_{2}^T & r_{3}^T & D_{-3} & A+D_{-3}\\
	r_{-1}^T & r_{2}^T & r_{3}^T & O & O & O & A+D_{3} & O
\end{pmatrix},$$

and
$$C_5=\begin{pmatrix}
	1 & 0 & 0 & 0 & 0 & 0 & O & O\\
	0 & 1 & 0 & 1 & 0 & 1 & O & r_{-1}\\
	0 & 0 & 0 & 0 & 1 & 0 & O & O\\
	1 & 1 & 1 & 0 & 0 & 1 & r_{-1} & r_2\\
	1 & 1 & [\frac{-3}{\widetilde{n}}] & 0 & 0 & 0 & r_{-3} & O\\
	0 & 0 & [\frac{-3}{\widetilde{n}}] & 1 & 1 & 1+[\frac{-3}{\widetilde{n}}] & O & r_{-3}\\
	O & O & O & r_{-1}^T & r_{2}^T & r_{3}^T & D_{-3} & A+D_{-3}\\
	r_{-1}^T & r_{2}^T & r_{3}^T & O & O & O & A+D_{3} & O
\end{pmatrix},$$

and 

$$C_6=\begin{pmatrix}
	1 & 0 & 0 & 0 & 0 & 0 & O & O\\
	0 & 1 & 0 & 0 & 0 & 0 & O & O\\
	0 & 0 & 0 & 0 & 1 & 0 & O & O\\
	1 & 0 & 1 & 0 & 0 & 0 & r_{-1} & O\\
	1 & 1 & [\frac{-3}{\widetilde{n}}] & 0 & 0 & 0 & r_{-3} & O\\
	0 & 0 & [\frac{-3}{\widetilde{n}}] & 1 & 1 & 1+[\frac{-3}{\widetilde{n}}] & O & r_{-3}\\
	O & O & O & r_{-1}^T & r_{2}^T & r_{3}^T & D_{-3} & A+D_{-3}\\
	r_{-1}^T & r_{2}^T & r_{3}^T & O & O & O & A+D_{3} & O
\end{pmatrix}.$$

\begin{theorem}
	For $n=-3\widetilde{n}=-3p_1...p_t,p_i \neq 2,3$, then we have Monsky matrix 
	\begin{equation*}
		M_n=\begin{cases}
			C_4, &\text{if } \widetilde{n}\equiv 1({\rm mod }\ 8),\\
			C_5, &\text{if } \widetilde{n}\equiv 5({\rm mod }\ 8), \\
			C_6, &\text{if } \widetilde{n}\equiv 3({\rm mod }\ 4).
		\end{cases}
	\end{equation*} 
\end{theorem}

\emph{\bf Proof:} We list all the equivalent conditions and transform them to additional Legendre symbols to build system of linear equations. 

We know that $C_{\Lambda}(\mathbb{A})\neq \emptyset$ if and only if $C_{\Lambda}(\mathbb{Q}_S)\neq \emptyset$, i.e.,

1) $C_{\Lambda}(\mathbb{R})\neq \emptyset$ if and only if $b_2>0$ that is equivalent to $\xi_1=0$;

2) $C_{\Lambda}(\mathbb{Q}_{p_i})\neq \emptyset$ for $i=1,..,t$ if and only if 
\begin{equation*}
	\begin{cases}
		\left( \frac{\dot{b}_1}{p_i}\right)=\left( \frac{\dot{b}_2}{p_i}\right)=1, &\text{if}\ p_i \nmid b_1b_2\\
		\left( \frac{\dot{n}\dot{b}_1}{p_i}\right)=\left( \frac{\dot{b}_2}{p_i}\right)=1, &\text{if}\ p_i\ |\ b_1,p_i\nmid b_2 \\
		\left( \frac{-3\dot{b}_1}{p_i}\right)=\left( \frac{-\dot{n}\dot{b}_2}{p_i}\right)=1, &\text{if}\ p_i\nmid b_1,p_i\ |\  b_2\\
		\left( \frac{-3\dot{n}\dot{b}_1}{p_i}\right)=\left( \frac{-\dot{n}\dot{b}_2}{p_i}\right)=1, &\text{if}\ p_i\ |\ b_1,p_i\ |\  b_2 
	\end{cases}
\end{equation*}

by Hensel's Lemma. The conditions are equivalent to
\begin{equation*}
	\begin{cases}
		[\frac{\dot{b}_1}{p_i}]+[\frac{-3}{p_i}]y_i+[\frac{-3\widetilde{n}p_i}{p_i}]x_i=0 \\
		[\frac{\dot{b}_2}{p_i}]+[\frac{3\widetilde{n}p_i}{p_i}]y_i=0 
	\end{cases};
\end{equation*}
In other words, we have
\begin{equation*}
	\begin{cases}
		[\frac{-1}{p_i}]\gamma_1+[\frac{2}{p_i}]\gamma_2+[\frac{3}{p_i}]\gamma_3+[\frac{-3}{p_i}]y_i+\sum\limits_{j\neq i}[\frac{p_j}{p_i}]x_j+[\frac{-3\widetilde{n}p_i}{p_i}]x_i=0 \\
		[\frac{-1}{p_i}]\xi_1+[\frac{2}{p_i}]\xi_2+[\frac{3}{p_i}]\xi_3+\sum\limits_{j\neq i}[\frac{p_j}{p_i}]y_j+[\frac{3\widetilde{n}p_i}{p_i}]y_i=0 
	\end{cases};
\end{equation*}

3) $C_{\Lambda}(\mathbb{Q}_3)\neq \emptyset$ if and only if
\begin{equation*}
	\begin{cases}
		\left( \frac{\dot{b}_1}{3}\right)=\left( \frac{\dot{b}_2}{3}\right)=1, &\text{if}\ 3 \nmid b_1b_2\\
		\left( \frac{\dot{n}\dot{b}_1}{3}\right)=\left( \frac{\dot{b}_2}{3}\right)=1, &\text{if}\ 3\ |\ b_1,3\nmid b_2 \\
		\left( \frac{-\dot{n}\dot{b}_1}{3}\right)=\left( \frac{-\dot{n}\dot{b}_2}{3}\right)=1, &\text{if}\ 3\nmid b_1,3\ |\  b_2\\
		\left( \frac{-\dot{b}_1}{3}\right)=\left( \frac{-\dot{n}\dot{b}_2}{3}\right)=1, &\text{if}\ 3\ |\ b_1,3\ |\  b_2 
	\end{cases}
\end{equation*}

by Hensel's Lemma. The conditions are equivalent to

\begin{equation*}
	\begin{cases}		
		(1+[\frac{-3}{\widetilde{n}}])\gamma_3+[\frac{-3}{\widetilde{n}}]\xi_3+[\frac{-3}{\dot{b}_1}]=0\\
		[\frac{-3}{\widetilde{n}}]\xi_3+[\frac{-3}{\dot{b}_2}]=0
	\end{cases};
\end{equation*}
In other words, we have 
\begin{equation*}
	\begin{cases}
		\gamma_1+\gamma_2+(1+[\frac{-3}{\widetilde{n}}])\gamma_3+[\frac{-3}{\widetilde{n}}]\xi_3+\sum\limits_{i}[\frac{-3}{p_i}]x_i=0\\
		\xi_1+\xi_2+[\frac{-3}{\widetilde{n}}]\xi_3+\sum\limits_{i}[\frac{-3}{p_i}]y_i=0
	\end{cases};
\end{equation*}
4) $C_{\Lambda}(\mathbb{Q}_2)\neq \emptyset$. If $\widetilde{n}\equiv 1({\rm mod }\ 8)$, $C_{\Lambda}(\mathbb{Q}_2)\neq \emptyset$ if and only if $2\nmid b_1b_2$ and $\dot{b}_1\equiv 1({\rm mod }\ 4)$; if  $\widetilde{n}\equiv 5({\rm mod }\ 8)$, $C_{\Lambda}(\mathbb{Q}_2)\neq \emptyset$ if and only if (i):$2\nmid b_1b_2$ and $(\dot{b}_1,\dot{b}_2)\equiv (1,1),(1,5),(5,7),(5,3)({\rm mod }\ 8,{\rm mod }\ 8)$ or (ii):$2\nmid b_1,2\ |\ b_2$ and $(\dot{b}_1,\dot{b}_2)\equiv (7,3),$ $(3,1)({\rm mod }\ 8,{\rm mod }\ 4)$;
if $\widetilde{n}\equiv 3({\rm mod }\ 4)$, $C_{\Lambda}(\mathbb{Q}_2)\neq \emptyset$ if and only if $2\nmid b_1b_2$ and $\dot{b}_2\equiv 1({\rm mod }\ 4)$ by Hensel's Lemma. The conditions are  equivalent to

\begin{equation*}
	\begin{cases}
		\xi_2=0, \ \ \ \ \gamma_2=0,\\
		\gamma_1+\gamma_3+\sum\limits_{i}[\frac{-1}{p_i}]x_i=0
	\end{cases} \text{if } \widetilde{n}\equiv 1({\rm mod }\ 8),
\end{equation*}
and
\begin{equation*}
	\begin{cases}
		\xi_2+\gamma_1+\gamma_3+\sum\limits_{i}[\frac{-1}{p_i}]x_i=0, \ \ \ \ \gamma_2=0,\\
		\xi_1+\xi_2+\xi_3+\gamma_3+\sum\limits_{i}[\frac{-1}{p_i}]y_i+\sum\limits_{i}[\frac{2}{p_i}]x_i=0
	\end{cases} \text{if } \widetilde{n}\equiv 5({\rm mod }\ 8),
\end{equation*}
and
\begin{equation*}
	\begin{cases}
		\xi_2=0, \ \ \ \ \gamma_2=0,\\
		\xi_1+\xi_3+\sum\limits_{i}[\frac{-1}{p_i}]y_i=0
	\end{cases} \text{if } \widetilde{n}\equiv 3({\rm mod }\ 4).
\end{equation*}
The Monsky matrix can be naturally derived from the above system of linear equations.
$\Box$

Then we have a corollary about pairty of 2-rank as following:

\begin{corollary}
	For $n=-3\widetilde{n}=-3p_1...p_t,p_i\neq2,3$, the rank of $2$-Selmer group $s_2(E_n)$ is even (resp. odd) if and only if $n\equiv21(resp.\ 3,9,15)\pmod{24}$. 
\end{corollary}

\emph{\bf Proof:}
Note that $s_2(E_n)=2t+6-r(M_n)$, then we only need consider the rank of $M_n$ which is denoted by $r(M_n)$.

For $\widetilde{n}\equiv1(\text{mod 8})$, we have:
\begin{align*}
	r(C_4)&=r\begin{pmatrix}
		1 & 0 & 0 & 0 & 0 & 0 & O & O\\
		0 & 0 & 0 & 0 & 1 & 0 & O & O\\
		0 & 1 & 0 & 0 & 0 & 0 & O & O\\
		0 & 0 & 0 & 1 & 0 & 1 & O & r_{-1}\\
		1 & 1 & [\frac{-3}{\widetilde{n}}] & 0 & 0 & 0 & r_{-3} & O\\
		0 & 0 & [\frac{-3}{\widetilde{n}}] & 1 & 1 & 1+[\frac{-3}{\widetilde{n}}] & O & r_{-3}\\
		O & O & O & r_{-1}^T & r_{2}^T & r_{3}^T & D_{-3} & A+D_{-3}\\
		r_{-1}^T & r_{2}^T & r_{3}^T & O & O & O & A+D_{3} & O
	\end{pmatrix} \\
	&=r\begin{pmatrix}
		1 & 0 & 0 & 0 & 0 & 0 & O & O\\
		0 & 0 & 0 & 0 & 1 & 0 & O & O\\
		0 & 1 & 0 & 0 & 0 & 0 & O & O\\
		0 & 0 & 0 & 1 & 0 & 1 & O & r_{-1}\\
		0 & 0 & [\frac{-3}{\widetilde{n}}] & 0 & 0 & 0 & r_{-3} & O\\
		O & O & O & O & O & r_{-3}^T & D_{-3} & A+D_{-3}+r_{-1}^Tr_{-1}\\
		O & O & r_{3}^T & O & O & O & A+D_{3} & O
	\end{pmatrix} \\ &\text{  by row transformations.} \\
	&=r\begin{pmatrix}
		1 & 0 & 0 & 0 & 0 & 0 & O & O\\
		0 & 0 & 0 & 0 & 1 & 0 & O & O\\
		0 & 1 & 0 & 0 & 0 & 0 & O & O\\
		0 & 0 & 0 & 1 & 0 & 0 & O & O\\
		0 & 0 & 0 & 0 & 0 & 0 & r_{-3} & O\\
		O & O & O & O & O & r_{-3}^T & D_{-3} & A^T+D_{3}\\
		O & O & O & O & O & O & A+D_{3} & O
	\end{pmatrix} \\ 
&\text{  by column transformations, and }A+A^T+r_{-1}^Tr_{-1}+D_{-1}=O \\
	&=4+r\begin{pmatrix}
         0  & r_{-3} & O\\
		r_{-3}^T & D_{-3}+r_{-3}^Tr_{-3} & A^T+D_{3}\\
        O  & A+D_{3} & O
	\end{pmatrix}
\end{align*}

Then we have $r(C_4)$ is even for $n\equiv21(\text{mod 24})$.

For $\widetilde{n}\equiv5(\text{mod 8})$, we have:
\begin{align*}
	r(C_5)&=r\begin{pmatrix}
		1 & 0 & 0 & 0 & 0 & 0 & O & O\\
		0 & 1 & 0 & 1 & 0 & 1 & O & r_{-1}\\
		0 & 0 & 0 & 0 & 1 & 0 & O & O\\
		1 & 1 & 1 & 0 & 0 & 1 & r_{-1} & r_2\\
		1 & 1 & [\frac{-3}{\widetilde{n}}] & 0 & 0 & 0 & r_{-3} & O\\
		0 & 0 & [\frac{-3}{\widetilde{n}}] & 1 & 1 & 1+[\frac{-3}{\widetilde{n}}] & O & r_{-3}\\
		O & O & O & r_{-1}^T & r_{2}^T & r_{3}^T & D_{-3} & A+D_{-3}\\
		r_{-1}^T & r_{2}^T & r_{3}^T & O & O & O & A+D_{3} & O
	\end{pmatrix} \\ 
	&=r\begin{pmatrix}
		1 & 0 & 0 & 0 & 0 & 0 & O & O\\
		0 & 0 & 0 & 0 & 1 & 0 & O & O\\
		0 & 1 & 0 & 1 & 0 & 1 & O & r_{-1}\\
		0 & 1 & 1 & 0 & 0 & 1 & r_{-1} & r_{2}\\
		0 & 0 & 1+[\frac{-3}{\widetilde{n}}] & 0 & 0 & 1 & r_{3} & r_{2}\\
		0 & 0 & 0 & 0 & 0 & 0 & O & O\\
		O & O & r_{-3}^T+[\frac{-3}{\widetilde{n}}]r_{3}^T & O & O & O & D_{-3}+r_{-1}^Tr_{-1}+r_{3}^Tr_{3} & A+D_{-3}+r_{-1}^Tr_{-1}+r_{-3}^Tr_{2}\\
		O & O & r_{3}^T+[\frac{-3}{\widetilde{n}}]r_{2}^T & O & O & O & A+D_{3}+r_{2}^Tr_{-3} & O
	\end{pmatrix} \\ 
&\text{  by row transformations.} \\
	&=r\begin{pmatrix}
		1 & 0 & 0 & 0 & 0 & 0 & O & O\\
		0 & 0 & 0 & 0 & 1 & 0 & O & O\\
		0 & 0 & 0 & 1 & 0 & 0 & O & O\\
		0 & 1 & 0 & 0 & 0 & 0 & O & O\\
		0 & 0 & 0 & 0 & 0 & 1 & O & O\\
		0 & 0 & 0 & 0 & 0 & 0 & O & O\\
		O & O & O & O & O & O & D_{-3}+r_{-1}^Tr_{-1}+r_{3}^Tr_{3} & A^T+D_{3}+r_{-3}^Tr_{2}\\
		O & O & O & O & O & O & A+D_{3}+r_{2}^Tr_{-3} & O
	\end{pmatrix} \\ 
&\text{  by column transformations, and }A+A^T+r_{-1}^Tr_{-1}+D_{-1}=O \\
	&=5+r\begin{pmatrix}
		D_{-3}+r_{-1}^Tr_{-1}+r_{3}^Tr_{3} & A^T+D_{3}+r_{-3}^Tr_{2}\\
        A+D_{3}+r_{2}^Tr_{-3} & O
	\end{pmatrix} 
\end{align*}

Then we have $r(C_5)$ is odd for $n\equiv9(\text{mod 24})$. 

For $\widetilde{n}\equiv3(\text{mod 4})$, we have:
\begin{align*}
	r(C_6)&=r\begin{pmatrix}
		1 & 0 & 0 & 0 & 0 & 0 & O & O\\
		0 & 1 & 0 & 0 & 0 & 0 & O & O\\
		0 & 0 & 0 & 0 & 1 & 0 & O & O\\
		1 & 0 & 1 & 0 & 0 & 0 & r_{-1} & O\\
		1 & 1 & [\frac{-3}{\widetilde{n}}] & 0 & 0 & 0 & r_{-3} & O\\
		0 & 0 & [\frac{-3}{\widetilde{n}}] & 1 & 1 & 1+[\frac{-3}{\widetilde{n}}] & O & r_{-3}\\
		O & O & O & r_{-1}^T & r_{2}^T & r_{3}^T & D_{-3} & A+D_{-3}\\
		r_{-1}^T & r_{2}^T & r_{3}^T & O & O & O & A+D_{3} & O
	\end{pmatrix} \\
	&=r\begin{pmatrix}
		1 & 0 & 0 & 0 & 0 & 0 & O & O\\
		0 & 0 & 0 & 0 & 1 & 0 & O & O\\
		0 & 1 & 0 & 0 & 0 & 0 & O & O\\
		0 & 0 & 1 & 0 & 0 & 0 & r_{-1} & O\\
		0 & 0 & [\frac{3}{\widetilde{n}}] & 0 & 0 & 0 & r_{3} & O\\
		0 & 0 & 0 & 0 & 0 & 0 & O & O\\
		O & O & O & r_{-1}^T & O & r_{3}^T & D_{-3} & A+D_{-3}\\
		O & O & r_{3}^T & O & O & O & A+D_{-3} & O
	\end{pmatrix} \\ &\text{  by row transformations.} \\
	&=r\begin{pmatrix}
		1 & 0 & 0 & 0 & 0 & 0 & O & O\\
		0 & 0 & 0 & 0 & 1 & 0 & O & O\\
		0 & 1 & 0 & 0 & 0 & 0 & O & O\\
		0 & 0 & 0 & 0 & 0 & 0 & r_{-1} & O\\
		0 & 0 & 0 & 0 & 0 & 0 & r_{3} & O\\
		0 & 0 & 0 & 0 & 0 & 0 & O & O\\
		O & O & O & r_{-1}^T & O & r_{3}^T & D_{-3}+r_{-1}^Tr_{-1}+r_{3}^Tr_{3} & A^T+D_{3}\\
		O & O & O & O & O & O & A+D_{3} & O
	\end{pmatrix} \\ 
&\text{  by column transformations, and }A+A^T+r_{-1}^Tr_{-1}+D_{-1}=O \\
	&=3+r\begin{pmatrix}
 0 & 0 & r_{-1} & O\\
 0 & 0 & r_{3} & O\\
 r_{-1}^T & r_{3}^T & D_{-3}+r_{-1}^Tr_{-1}+r_{3}^Tr_{3} & A^T+D_{3}\\
 O & O & A+D_{3} & O
	\end{pmatrix} 
\end{align*}

Then we have $r(C_6)$ is odd for $n\equiv3,15(\text{mod 24})$.  $\Box$

Now we prove a result about non $2\pi/3$-congruent numbers:

\begin{corollary}
	For $n=-3\widetilde{n}=-3p_1...p_t\equiv21\pmod{24},p_i\neq2,3$, if $r_4(n)=0$, then $s_2(E_n)=2$. 
\end{corollary}

\emph{\bf Proof:} 
For $n\equiv21(\text{mod }24)$, we have

$$ s_2(E_n)=2t+6-r(M_n)=2t+2-r\begin{pmatrix}
	0  & r_{-3} & O\\
	r_{-3}^T & D_{-3}+r_{-3}^Tr_{-3} & A^T+D_{3}\\
	O & A+D_{3} & O
\end{pmatrix}.$$ Then $s_2\leq2t+2-2r(	r_{-3}^T \ A^T+D_{3})=2t+2-2r(R(n))=2+2r_4(n)=2.$

Note that $s_2(E_n)\geq2$ by the fact $E_n$ have $4$ elements of $2$-torsion, we have $s_2(E_n)=2$. $\Box$

\subsection{\bf $\eta=-6$} 
Let 
$$D_2=\begin{pmatrix}
	1 & 0 & 0 & 0 & 0 & 0 & O & O\\
	0 & 0 & 0 & 1 & [\frac{-1}{\widetilde{n}}] & 1 & O & r_{-1}\\
	0 & 1 & 0 & 0 & 1+[\frac{2}{\widetilde{n}}] & 1 & O & r_{2}\\
	1 & 1 & 1+[\frac{-3}{\widetilde{n}}] & 0 & 0 & 0 & r_{-3} & O\\
	0 & 0 & 1+[\frac{-3}{\widetilde{n}}] & 1 & 1 & [\frac{-3}{\widetilde{n}}] & O & r_{-3}\\
	O & O & O & r_{-1}^T & r_{2}^T & r_{3}^T & D_{-3} & A+D_{-6}\\
	r_{-1}^T & r_{2}^T & r_{3}^T & O & O & O & A+D_{6} & O
\end{pmatrix}.$$

\begin{theorem}
	For $n=-6\widetilde{n}=-6p_1...p_t,p_i \neq 2,3$, then we have Monsky matrix $M_n=D_2$.
\end{theorem}

\emph{\bf Proof:} We list all the equivalent conditions and transform them to additional Legendre symbols to build system of linear equations. 

$C_{\Lambda}(\mathbb{A})\neq \emptyset$ if and only if $C_{\Lambda}(\mathbb{Q}_S)\neq \emptyset$, i.e., 

1) $C_{\Lambda}(\mathbb{R})\neq \emptyset$ if and only if 		$b_2>0$ that is equivalent to $\xi_1=0$;

2) $C_{\Lambda}(\mathbb{Q}_{p_i})\neq \emptyset$ for $i=1,..,t$ if and only if 
\begin{equation*}
	\begin{cases}
		\left( \frac{\dot{b}_1}{p_i}\right)=\left( \frac{\dot{b}_2}{p_i}\right)=1, &\text{if}\ p_i \nmid b_1b_2\\
		\left( \frac{\dot{n}\dot{b}_1}{p_i}\right)=\left( \frac{\dot{b}_2}{p_i}\right)=1, &\text{if}\ p_i\ |\ b_1,p_i\nmid b_2 \\
		\left( \frac{-3\dot{b}_1}{p_i}\right)=\left( \frac{-\dot{n}\dot{b}_2}{p_i}\right)=1, &\text{if}\ p_i\nmid b_1,p_i\ |\  b_2\\
		\left( \frac{-3\dot{n}\dot{b}_1}{p_i}\right)=\left( \frac{-\dot{n}\dot{b}_2}{p_i}\right)=1, &\text{if}\ p_i\ |\ b_1,p_i\ |\  b_2 
	\end{cases}
\end{equation*}

by Hensel's Lemma. The conditions are equivalent to
\begin{equation*}
	\begin{cases}
		[\frac{\dot{b}_1}{p_i}]+[\frac{-3}{p_i}]y_i+[-\frac{6\widetilde{n}p_i}{p_i}]x_i=0 \\
		[\frac{\dot{b}_2}{p_i}]+[\frac{6\widetilde{n}p_i}{p_i}]y_i=0 
	\end{cases};
\end{equation*}
In other words, we have

\begin{equation*}
	\begin{cases}
		[\frac{-1}{p_i}]\gamma_1+[\frac{2}{p_i}]\gamma_2+[\frac{3}{p_i}]\gamma_3+[\frac{-3}{p_i}]y_i+\sum\limits_{j\neq i}[\frac{p_j}{p_i}]x_j+[\frac{-6\widetilde{n}p_i}{p_i}]x_i=0 \\
		[\frac{-1}{p_i}]\xi_1+[\frac{2}{p_i}]\xi_2+[\frac{3}{p_i}]\xi_3+\sum\limits_{j\neq i}[\frac{p_j}{p_i}]y_j+[\frac{6\widetilde{n}p_i}{p_i}]y_i=0 
	\end{cases};
\end{equation*}

3) $C_{\Lambda}(\mathbb{Q}_3)\neq \emptyset$ if and only if
\begin{equation*}
	\begin{cases}
		\left( \frac{\dot{b}_1}{3}\right)=\left( \frac{\dot{b}_2}{3}\right)=1, &\text{if}\ 3 \nmid b_1b_2\\
		\left( \frac{\dot{n}\dot{b}_1}{3}\right)=\left( \frac{\dot{b}_2}{3}\right)=1, &\text{if}\ 3\ |\ b_1,3\nmid b_2 \\
		\left( \frac{-\dot{n}\dot{b}_1}{3}\right)=\left( \frac{-\dot{n}\dot{b}_2}{3}\right)=1, &\text{if}\ 3\nmid b_1,3\ |\  b_2\\
		\left( \frac{-\dot{b}_1}{3}\right)=\left( \frac{-\dot{n}\dot{b}_2}{3}\right)=1, &\text{if}\ 3\ |\ b_1,3\ |\  b_2 
	\end{cases}
\end{equation*}

by Hensel's Lemma. The conditions are equivalent to

\begin{equation*}
	\begin{cases}		
		[\frac{-3}{\widetilde{n}}]\gamma_3+(1+[\frac{-3}{\widetilde{n}}])\xi_3+[\frac{-3}{\dot{b}_1}]=0\\
		(1+[\frac{-3}{\widetilde{n}}])\xi_3+[\frac{-3}{\dot{b}_2}]=0
	\end{cases};
\end{equation*}
In other words, we have 
\begin{equation*}
	\begin{cases}
		\gamma_1+\gamma_2+[\frac{-3}{\widetilde{n}}]\gamma_3+(1+[\frac{-3}{\widetilde{n}}])\xi_3+\sum\limits_{i}[\frac{-3}{p_i}]x_i=0\\
		\xi_1+\xi_2+(1+[\frac{-3}{\widetilde{n}}])\xi_3+\sum\limits_{i}[\frac{-3}{p_i}]y_i=0
	\end{cases};
\end{equation*}

4) $C_{\Lambda}(\mathbb{Q}_2)\neq \emptyset$.
If  $\widetilde{n}\equiv 1({\rm mod }\ 8)$, $C_{\Lambda}(\mathbb{Q}_2)\neq \emptyset$ if and only if 

\begin{equation*}
	\begin{cases}
		(\dot{b}_1,\dot{b}_2)\equiv(1,1),(1,7)(\text{mod }8,\text{mod }8), &\text{if}\ 2 \nmid b_1b_2\\
		\dot{b}_1\equiv5(\text{mod }8), &\text{if}\ 2\ |\ b_1,2\nmid b_2 \\
		(\dot{b}_1,\dot{b}_2)\equiv(5,3),(5,5)(\text{mod }8,\text{mod }8), &\text{if}\ 2\nmid b_1,2\ |\  b_2\\
		\dot{b}_1\equiv1(\text{mod }8), &\text{if}\ 2\ |\ b_1,2\ |\  b_2 
	\end{cases}.
\end{equation*}

If  $\widetilde{n}\equiv 3({\rm mod }\ 8)$, $C_{\Lambda}(\mathbb{Q}_2)\neq \emptyset$ if and only if 

\begin{equation*}
	\begin{cases}
		(\dot{b}_1,\dot{b}_2)\equiv(1,1),(1,3)(\text{mod }8,\text{mod }8), &\text{if}\ 2 \nmid b_1b_2\\
		\dot{b}_1\equiv7(\text{mod }8), &\text{if}\ 2\ |\ b_1,2\nmid b_2 \\
		(\dot{b}_1,\dot{b}_2)\equiv(5,1),(5,3)(\text{mod }8,\text{mod }8), &\text{if}\ 2\nmid b_1,2\ |\  b_2\\
		\dot{b}_1\equiv3(\text{mod }8), &\text{if}\ 2\ |\ b_1,2\ |\  b_2 
	\end{cases}.
\end{equation*}

If  $\widetilde{n}\equiv 5({\rm mod }\ 8)$, $C_{\Lambda}(\mathbb{Q}_2)\neq \emptyset$ if and only if 

\begin{equation*}
	\begin{cases}
		(\dot{b}_1,\dot{b}_2)\equiv(1,1),(1,7)(\text{mod }8,\text{mod }8), &\text{if}\ 2 \nmid b_1b_2\\
		\dot{b}_1\equiv1(\text{mod }8), &\text{if}\ 2\ |\ b_1,2\nmid b_2 \\
		(\dot{b}_1,\dot{b}_2)\equiv(5,1),(5,7)(\text{mod }8,\text{mod }8), &\text{if}\ 2\nmid b_1,2\ |\  b_2\\
		\dot{b}_1\equiv5(\text{mod }8), &\text{if}\ 2\ |\ b_1,2\ |\  b_2 
	\end{cases}.
\end{equation*}

If  $\widetilde{n}\equiv 7({\rm mod }\ 8)$, $C_{\Lambda}(\mathbb{Q}_2)\neq \emptyset$ if and only if 

\begin{equation*}
	\begin{cases}
		(\dot{b}_1,\dot{b}_2)\equiv(1,1),(1,3)(\text{mod }8,\text{mod }8), &\text{if}\ 2 \nmid b_1b_2\\
		\dot{b}_1\equiv3(\text{mod }8), &\text{if}\ 2\ |\ b_1,2\nmid b_2 \\
		(\dot{b}_1,\dot{b}_2)\equiv(5,5),(5,7)(\text{mod }8,\text{mod }8), &\text{if}\ 2\nmid b_1,2\ |\  b_2\\
		\dot{b}_1\equiv7(\text{mod }8), &\text{if}\ 2\ |\ b_1,2\ |\  b_2 
	\end{cases}.
\end{equation*}

The conditions are equivalent to

\begin{equation*}
	\begin{cases}
		[\frac{-1}{\widetilde{n}}]\gamma_2+[\frac{-1}{\dot{b}_1}]=0,\\
		(1+[\frac{2}{\widetilde{n}}])\gamma_2+\xi_2+[\frac{2}{\dot{b}_1}]=0,\\
		(1+[\frac{2}{\widetilde{n}}])\xi_2+[\frac{2}{\dot{b}_2}]+[\frac{-1}{\widetilde{n}}][\frac{-1}{\dot{b}_2}]=0,\ \ \text{if } \eta_2=0
	\end{cases}.
\end{equation*}

Note if $C_{\Lambda}(\mathbb{Q}_3)$ and $C_{\Lambda}(\mathbb{Q}_{p_i})\neq \emptyset$, we have $(1+[\frac{2}{\widetilde{n}}])\xi_2+[\frac{2}{\dot{b}_2}]+[\frac{-1}{\widetilde{n}}][\frac{-1}{\dot{b}_2}]=0$ by $(1+[\frac{-3}{\widetilde{n}}])\xi_3+[\frac{-3}{\dot{b}_2}]=0$ and $[\frac{\dot{b}_2}{p_i}]+[\frac{6\widetilde{n}p_i}{p_i}]y_i=0$. In other words, if we assume $C_{\Lambda}(\mathbb{Q}_3)$ and $C_{\Lambda}(\mathbb{Q}_{p_i})\neq \emptyset$, then $C_{\Lambda}(\mathbb{Q}_2)\neq \emptyset$ if and only if

\begin{equation*}
	\begin{cases}
		\gamma_1+[\frac{-1}{\widetilde{n}}]\gamma_2+\gamma_3+\sum\limits_{i}[\frac{-1}{p_i}]x_i=0,\\
		(1+[\frac{2}{\widetilde{n}}])\gamma_2+\gamma_3+\xi_2+\sum\limits_{i}[\frac{2}{p_i}]x_i=0
	\end{cases}.
\end{equation*}
The Monsky matrix can be naturally derived from the above system of linear equations.
$\Box$

Then we have a corollary about pairty of 2-rank as following:

\begin{corollary}
	For $n=-6\widetilde{n}=-6p_1...p_t,p_i\neq2,3$, the rank of $2$-Selmer group $s_2(E_n)$ is even. 
\end{corollary}

\emph{\bf Proof:}
Note that $s_2(E_n)=2t+6-r(M_n)$, then we only need consider the rank of $M_n$ which is denoted by $r(M_n)$. We have:

\begin{align*}
	r(D_2)&=r\begin{pmatrix}
		1 & 0 & 0 & 0 & 0 & 0 & O & O\\
		0 & 0 & 0 & 1 & [\frac{-1}{\widetilde{n}}] & 1 & O & r_{-1}\\
		0 & 1 & 0 & 0 & 1+[\frac{2}{\widetilde{n}}] & 1 & O & r_{2}\\
		1 & 1 & 1+[\frac{-3}{\widetilde{n}}] & 0 & 0 & 0 & r_{-3} & O\\
		0 & 0 & 1+[\frac{-3}{\widetilde{n}}] & 1 & 1 & [\frac{-3}{\widetilde{n}}] & O & r_{-3}\\
		O & O & O & r_{-1}^T & r_{2}^T & r_{3}^T & D_{-3} & A+D_{-6}\\
		r_{-1}^T & r_{2}^T & r_{3}^T & O & O & O & A+D_{6} & O
	\end{pmatrix} \\
	&=r\begin{pmatrix}
		1 & 0 & 0 & 0 & 0 & 0 & O & O\\
		0 & 0 & 0 & 1 & [\frac{-1}{\widetilde{n}}] & 1 & O & r_{-1}\\
		0 & 1 & 0 & 0 & [\frac{2}{\widetilde{n}}] & 1 & O & r_{2}\\
		0 & 0 & 1+[\frac{-3}{\widetilde{n}}] & 0 & 1+[\frac{2}{\widetilde{n}}] & 1 & r_{-3} & r_{2}\\
		0 & 0 & 0 & 0 & 0 & 0 & O & O\\
		O & O & O & O & [\frac{-1}{\widetilde{n}}]r_{-1}^Tr_{2}^T & r_{-3}^T & D_{-3} & A+D_{-6}+r_{-1}^Tr_{-1}\\
		O & O & r_3^T & O & (1+[\frac{2}{\widetilde{n}}])r_{2}^T & r_{2}^T & A+D_{6} & r_{2}^Tr_{2}
	\end{pmatrix} \\ &\text{  by row transformations.} \\
	&=r\begin{pmatrix}
		1 & 0 & 0 & 0 & 0 & 0 & O & O\\
		0 & 0 & 0 & 1 & 0 & 0 & O & O\\
		0 & 1 & 0 & 0 & 0 & 0 & O & O\\
		0 & 0 & 0 & 0 & 0 & 1 & r_{-3} & r_{2}\\
		0 & 0 & 0 & 0 & 0 & 0 & O & O\\
		O & O & O & O & O & r_{-3}^T & D_{-3} & A+D_{-6}+r_{-1}^Tr_{-1}\\
		O & O & O & O & O & r_{2}^T & A+D_{6} & r_{2}^Tr_{2}
	\end{pmatrix} \\ &\text{  by column transformations.} \\
	&=3+r\begin{pmatrix}
1 & r_{-3} & r_{2}\\
r_{-3}^T & D_{-3} & A^T+D_{6}\\
 r_{2}^T & A+D_{6} & r_{2}^Tr_{2}
	\end{pmatrix} \\
	&\text{  by } A+A^T+r_{-1}^Tr_{-1}+D_{-1}=O. \\
	&=4+r\begin{pmatrix}
 D_{-3}+r_{-3}^Tr_{-3} & A^T+D_{6}+r_{-3}^Tr_{2}\\
A+D_{6}+r_{2}^Tr_{-3} & O
	\end{pmatrix} 
\end{align*}

Then we have $r(D_2)$ is even. $\Box$ 

Now we prove a result about non $2\pi/3$-congruent numbers:

\begin{corollary}
	For $n=-6\widetilde{n}=-6p_1...p_t\equiv6,18\pmod{24},p_i\neq2,3$, if $r_4(n)=0$, then $s_2(E_n)=2$. 
\end{corollary}

\emph{\bf Proof:} 
For $n\equiv6,18(\text{mod }24)$, we have

$$ s_2(E_n)=2t+6-r(M_n)=2t+3-r\begin{pmatrix}
	1 & r_{-3} & r_{2}\\
	r_{-3}^T & D_{-3} & A^T+D_{6}\\
	r_{2}^T & A+D_{6} & r_{2}^Tr_{2}
\end{pmatrix}.$$ Then 
\begin{align*}
	s_2&=2t+4-2r\begin{pmatrix}
				1 & 0 & O & r_{2}\\
	0&	1 & r_{-3} & r_{2}\\
	O &	r_{-3}^T & D_{-3} & A^T+D_{6}\\
	r_{2}^T &	r_{2}^T & A+D_{6} & O
	\end{pmatrix} \\
&=2t+4-2r\begin{pmatrix}
	1 & O & 1 & r_{2}\\
	O & D_{-3} &	r_{-3}^T  & A^T+D_{6}\\
1 & r_{-3} &	0	 & O\\
	r_{2}^T &	 A+D_{6} & O  & O
\end{pmatrix} \\
&\leq 2t+4-2r\begin{pmatrix}
	1 & r_{-3} \\
	r_{2}^T &	 A+D_{6} 
\end{pmatrix} \\
&=2t+4-r(R(n))\\
&=2+c_4(n)=2.
\end{align*}

Note that $s_2(E_n)\geq2$ by the fact $E_n$ have $4$ elements of $2$-torsion, we have $s_2(E_n)=2$. $\Box$

\section{Some Examples}
In this section, we give some non $\pi/3$ or $2\pi/3$-congruent numbers examples for $n=pq$.

\begin{theorem}
	For $n=pq\equiv 5\pmod{24}$ and $p,q$ are primes, if $1)$ $[\frac{p}{q}]=1$ or $2)$ $[\frac{p}{q}]=0$ and $[\frac{\beta}{q}]=1$ where $\beta^2=p\pmod{q}$, then $n$ is a non $\pi/3$-congruent number. Moreover, the density that $n$ is non $\pi/3$-congruent number is at least $75\%$. 
\end{theorem}

\emph{\bf Proof:} Let $[\frac{-1}{p}]=[\frac{-1}{q}]=u$, $[\frac{q}{p}]=v$. Without loss of generality, we assume $[\frac{-3}{p}]=0$.

By theorem 6.2, we have system of linear equations as following:

\begin{equation*}
	\begin{cases}
		\gamma_1+\xi_1=0 \\
		\gamma_2=0\\
		\xi_2=0\\
		\xi_3=0\\
		\gamma_1=y_2\\
		\gamma_3=y_2+ux_1+ux_2\\
		(u+v)x_1+(u+v)x_2=0\\
		(u+v)y_1+(u+v)y_2=0
	\end{cases} (\star)
\end{equation*}

When $[\frac{p}{q}]=u+v=1$, the basis of the solution set to $(\star)$ is $(1,0,0,1,0,1,1,$ $1,0,0)$ and $(0,0,0,0,0,0,0,0,1,1)$ which correspond to $(-3,-n)$ and $(n,1)$ in ${\rm Sel}_2(E_{n})$, respectively. In this case, those points are torsion in $E_{n}(\mathbb{Q})$, i.e., $n$ is non $\pi/3$-congruent number.

Now we assume $[\frac{p}{q}]=u+v=0$.

 The 2-Selmer group ${\rm Sel}_2(E_n)=\{\left\langle(1,p),(-3,-n),(3^uq,1),(n,1)\right\rangle\}$. $n$ is non $\pi/3$-congruent number if the Cassels pairing $\left \langle (1,p),(3^uq,1) \right \rangle=1$. And we have $\left \langle (1,p),(3^uq,1) \right \rangle$ =$\sum\limits_{\mathcal{p}|24n\infty} \left[ L_1(P_\mathcal{p}),3^uq \right]_\mathcal{p}\left[ L_3(P_\mathcal{p}),3^uq \right]_\mathcal{p}$.

For $\Lambda=(1,p)$, we take $Q_1=(t,u_2,u_3)=(0,1,1)$, $Q_3=(t,u_1,u_2)=(b,pa,c)$ such that $Q_i$ is a solution of $H_i$ for $i=1,3$. And let $(a,b,c)$ be a primitive integer solution to $px^2-qy^2=z^2$. We assume that $c$ is even, $c\equiv1\pmod{3}$, $a$, $b$ is odd, $a\equiv1\pmod{4}$ and $3\mid a$. If not, we have a new solution $$(\widetilde{a},\widetilde{b},\widetilde{c}):=(-(p+q)a+2qb,(p+q)b-2pa,(p-q)c).$$
Then the 2-adic valuation of $\widetilde{c}$ is $2$, while those of $\widetilde{a}$ and $\widetilde{b}$ is $1$. And the 3-adic valuation of $\widetilde{a}$ is equal to that of $b$. Moreover, the 3-adic valuation of $\widetilde{b}$ is equal to that of $a$. By dividing ${\rm gcd}(\widetilde{a},\widetilde{b},\widetilde{c})$, we have a new primitive solution $(a^{\prime},b^{\prime},c^{\prime})$ such that $c^\prime$ is even and $(a^\prime,b^\prime)\equiv(b,a)(\text{mod}\ 3, \text{mod}\ 3)$. Then we find a solution $(a,b,c)$ while $c$ is even and $3\mid a$ as those process. If necessary, we can take $-a$ replace $a$ or $-c$ replace $c$ to make $a\equiv1\pmod{4}$ and  $c\equiv1\pmod{3}$. 

$L_1$ and $L_3$ are 

$$ \left\{
\begin{aligned}
	L_1:&\ u_2-u_3=0, \\
	L_3:&\ -bqt-cu_2+au_1=0.
\end{aligned}
\right.
$$

Since $3^uq>0$, $[\frac{-1}{3^uq}]=u+u=0$ and $[\frac{3^uq}{p}]=u+v=0$, we have $\left[ L_i(P_\infty),3^uq \right]_\infty=\left[ L_i(P_p),3^uq \right]_p=0$ for $i=1,3$. Let $P_3=(t,u_1,u_2,u_3)=(0,1,1,-1)$, then we have $\left[ L_1(P_3),3^uq \right]_3+\left[ L_3(P_3),3^uq \right]_3=u[\frac{c}{3}]=0$.

For $\mathcal{p}=q$, we take $P_q=(t,u_1,u_2,u_3)=(0,\beta,1,-1)$ where $\beta^2=p\pmod{q}$ and $q\mid c+\beta a$, i.e., $c+\beta a=2\beta a\pmod{q}$. Then $\left[ L_1(P_q),3^uq \right]_q+\left[ L_3(P_q),3^uq \right]_q=[\frac{\beta a}{q}]$. 

Now we consider $\mathcal{p}=2$. Let $u_2=2,u_3=0$ and $u_1^2=t^2=1$ such that $au_1-qbt\equiv2\pmod{4}$, we have $\left[ L_1(P_2),3^uq \right]_2+\left[ L_3(P_2),3^uq \right]_2=0$. 

Note that $-qb^2\equiv c^2\pmod{a}$, i.e., $[\frac{-q}{a}]=0$, then we have $[\frac{a}{q}]=0$. Therefore, $\left \langle (1,p),(3^uq,1) \right \rangle=[\frac{\beta}{q}]$.

In summary, $n=pq\equiv 5\pmod{24}$ is non $\pi/3$-congruent number if $1)$ $[\frac{p}{q}]=1$ or $2)$ $[\frac{p}{q}]=0$ and $[\frac{\beta}{q}]=1$ where $\beta$ satisfy $\beta^2=p\pmod{q}$ and $q\mid c+\beta a$. By Chebotarev's density theorem, the density is at least $ 75 \% $. $\Box$

\begin{theorem}
	For $n=pq\equiv 11\pmod{24}$, $p,q$ are primes, $1)$ $[\frac{p}{q}]=1$ or $2)$ $[\frac{p}{q}]=0$ and $[\frac{\beta}{q}]=1$ where $\beta^2=p\pmod{q}$, then $n$ is a non $2\pi/3$-congruent number. Moreover, the density that $n$ is non $2\pi/3$-congruent number is at least $75\%$. 
\end{theorem}

\emph{\bf Proof:} Note $[\frac{-1}{p}]=[\frac{-1}{q}]+1=u$, $[\frac{q}{p}]=[\frac{p}{q}]=v$. Without loss of generality, we assume $[\frac{-3}{p}]=0$.

By theorem 6.14, we have system of linear equations as following:

\begin{equation*}
	\begin{cases}
		\xi_1=0 \\
		\xi_2=0\\
		\xi_3=0\\
		\gamma_1=y_2+ux_1+(u+1)x_2\\
		\gamma_2=0\\
		\gamma_3=y_2\\
		vx_1+vx_2=0\\
		vy_1+vy_2=0
	\end{cases} (\star)
\end{equation*}

When $[\frac{p}{q}]=v=1$, the basis of the solution set to $(\star)$ is $(0,0,0,1,0,1,1,$ $1,0,0)$ and $(0,0,0,1,0,0,0,0,1,1)$ which correspond to $(-3,n)$ and $(-n,1)$ in ${\rm Sel}_2(E_{-n})$, respectively. In this case, those points are torsion in $E_{-n}(\mathbb{Q})$, i.e., $n$ is non $2\pi/3$-congruent number.

Now we assume $[\frac{p}{q}]=0$. 

The 2-Selmer group ${\rm Sel}_2(E_{-n})=\{\left\langle(1,p),(-3,-n),((-1)^uq,1),(n,1)\right\rangle\}$. $n$ is non $2\pi/3$-congruent number if the Cassels pairing $\left \langle (1,p),((-1)^uq,1) \right \rangle=1$. And we have $\left \langle (1,p),((-1)^uq,1) \right \rangle$ =$\sum\limits_{\mathcal{p}|24n\infty} \left[ L_1(P_\mathcal{p}),(-1)^uq \right]_\mathcal{p}\left[ L_3(P_\mathcal{p}),(-1)^uq \right]_\mathcal{p}$.

For $\Lambda=(1,p)$, we take $Q_1=(t,u_2,u_3)=(0,1,1)$, $Q_3=(t,u_1,u_2)=(b,pa,c)$ such that $Q_i$ is a solution of $H_i$ for $i=1,3$. Let $(a,b,c)$ is a primitive integer solution to $px^2+qy^2=z^2$. We assume $c$ is even, negative, $a$, $b$ is odd, $a\equiv1\pmod{4}$ and $3\mid a$. If not, we have a new solution $$(\widetilde{a},\widetilde{b},\widetilde{c}):=(-(p+q)a+2qb,(p+q)b-2pa,(p-q)c).$$
Then the 2-adic valuation of $\widetilde{c}$ is $2$, while those of $\widetilde{a}$ and $\widetilde{b}$ is $1$. And the 3-adic valuation of $\widetilde{a}$ is equal to that of $b$. Moreover, the 3-adic valuation of $\widetilde{b}$ is equal to that of $a$. By dividing ${\rm gcd}(\widetilde{a},\widetilde{b},\widetilde{c})$, we have a new primitive solution $(a^{\prime},b^{\prime},c^{\prime})$ such that $c^\prime$ is even and $(a^\prime,b^\prime)\equiv(b,a)(\text{mod}\ 3, \text{mod}\ 3)$. Then we can find a solution $(a,b,c)$ satisfies $c$ is even and $3\mid a$ as those process. If necessary, we can take $-a$ replace $a$ or $-c$ replace $c$ to make $a\equiv1\pmod{4}$ and  $c$ is negative. 

$L_1$ and $L_3$ is 

$$ \left\{
\begin{aligned}
	L_1:&\ u_2-u_3=0, \\
	L_3:&\ bqt-cu_2+au_1=0.
\end{aligned}
\right.
$$

For $\mathcal{p}=\infty$, let $P_\infty=(t,u_1,u_2,u_3)=(0,1,-1,\sqrt{p})$, we have $\left[ L_i(P_\infty),(-1)^uq \right]_\infty$ $=0$ for $i=1,3$.

For $\mathcal{p}=3$, let $u_2=0,u_3=1$, and $u_1^2=t^2=1$ such that $3\mid au_1-bt$, then we have $\left[ L_1(P_3),(-1)^uq \right]_3+\left[ L_3(P_3),(-1)^uq \right]_3=0$. 

For $\mathcal{p}=p$, let $P_p=(t,u_1,u_2,u_3)=(1,0,\alpha,s\alpha)$ where $s^2\equiv-3\pmod{p}$, $\alpha^2\equiv q\pmod{p}$ and $p\mid c+\alpha b$, then we have $\left[ L_i(P_p),(-1)^uq \right]_p=0$ for $i=1,3$.

For $\mathcal{p}=q$, let  $P_q=(t,u_1,u_2,u_3)=(0,\beta,1,-1)$ where $\beta^2=p\pmod{q}$ and $q\mid c+\beta a$, i.e., $c+\beta a=2\beta a\pmod{q}$, then we have $\left[ L_1(P_q),(-1)^uq \right]_q+\left[ L_3(P_q),(-1)^uq \right]_q=[\frac{\beta a}{q}]$. 

Now we consider $\mathcal{p}=2$. Let $u_2=2,u_3=0$, and $u_1^2=t^2=1$ such that $au_1+qbt\equiv2\pmod{4}$. we have $\left[ L_1(P_2),(-1)^uq \right]_2+\left[ L_3(P_2),(-1)^uq \right]_2=0$. 

Note that $qb^2\equiv c^2\pmod{a}$, i.e., $[\frac{q}{a}]=0$, then we have $[\frac{a}{q}]=0$. Therefore, $\left \langle (1,p),((-1)^uq,1) \right \rangle=[\frac{\beta}{q}]$.

In summary, $n=pq\equiv 11\pmod{24}$ is non $2\pi/3$-congruent number if $1)$ $[\frac{p}{q}]=1$ or $2)$ $[\frac{p}{q}]=0$ and $[\frac{\beta}{q}]=1$ where  $\beta^2=p\pmod{q}$ and $q\mid c+\beta a$. By Chebotarev's density theorem, the density is at least $ 75 \% $. $\Box$

\end{document}